\documentclass[12pt,reqno]{amsart}

\usepackage[T1]{fontenc}
\usepackage[utf8]{inputenc}

\usepackage{amsmath,amssymb,mathtools}
\usepackage{etoolbox}
\patchcmd{\subsection}{-.5em}{.5em}{}{\ClassError{amsart}{Unable to change subsection heading layout}{Check the definition of \string\subsection.}}
\usepackage{esvect}
\usepackage{microtype}
\usepackage{enumitem}
\usepackage{booktabs}
\usepackage{graphicx}
\graphicspath{{./}{../}}
\usepackage[dvipsnames]{xcolor}
\usepackage{tikz}
\usepackage{tikz-cd}
\usepackage{subcaption}
\usepackage{placeins}
\usepackage{float}

\definecolor{brianchonblue}{RGB}{65,130,210}
\definecolor{brianchonpurple}{RGB}{163,112,185}
\definecolor{brianchonred}{RGB}{218,48,48}

\usetikzlibrary{arrows.meta,calc,intersections,positioning}
\usetikzlibrary{calc,intersections}
\usetikzlibrary{calc,arrows.meta}
\usetikzlibrary{arrows.meta}


\usepackage{hyperref}
\usepackage[nameinlink,capitalise]{cleveref}

\allowdisplaybreaks

\hypersetup{
    colorlinks=true,
    linkcolor=blue,
    citecolor=Green,
    urlcolor=red
}

\newtheorem{theorem}{Theorem}[section]
\newtheorem{proposition}[theorem]{Proposition}
\newtheorem{lemma}[theorem]{Lemma}
\newtheorem{corollary}[theorem]{Corollary}

\theoremstyle{definition}
\newtheorem{definition}[theorem]{Definition}
\newtheorem{example}[theorem]{Example}

\theoremstyle{remark}
\newtheorem{remark}[theorem]{Remark}

\newcommand{\Ftwo}{\mathbb{F}_2}
\newcommand{\PP}{\mathbb{P}}

\newcommand{\eg}{\operatorname{eg}}

\newcommand{\cl}{\operatorname{cl}}

\newcommand{\id}{\mathrm{id}}

\newcommand{\Aut}{\operatorname{Aut}}

\newcommand{\Span}{\operatorname{span}}

\newcommand{\Hom}{\operatorname{Hom}}
\newcommand{\Hol}{\operatorname{Hol}}
\newcommand{\slodd}{\operatorname{sl}^{\mathrm{odd}}}
\newcommand{\cleven}{\operatorname{cl}^{\mathrm{even}}}
\newcommand{\tl}{\operatorname{tl}}

\newcommand{\Int}{\operatorname{Int}}

\title{Incidence Geometry with Polarity and Tiled Surfaces}

\author[Yuanqi Jiang, Peixuan Liu, Zhiyuan Liu]{%
Yuanqi Jiang$^{*}$,
Peixuan Liu$^{\dagger}$,
Zhiyuan Liu$^{\ddagger}$
}

\thanks{%
$^{*}$Nankai University, Tianjin
\href{mailto:jiangyuanqi@mail.nankai.edu.cn}{\nolinkurl{jiangyuanqi@mail.nankai.edu.cn}}.}

\thanks{%
$^{\dagger}$Southern University of Science and Technology, Shenzhen
\href{mailto:sun_don@foxmail.com}{\nolinkurl{sun_don@foxmail.com}}.}

\thanks{%
$^{\ddagger}$Southern University of Science and Technology, Shenzhen
\href{mailto:kracyliu050403@gmail.com}{\nolinkurl{kracyliu050403@gmail.com}}.}

\date{\today}

\advance\paperheight1.00cm
\advance\textheight1.00cm
\advance\vsize1.00cm
\advance\topskip-0.5cm
\advance\voffset-0.5cm
\DeclareTextFontCommand{\emph}{\color{Blue}\em}
\makeatletter
\g@addto@macro\@floatboxreset{\footnotesize}
\makeatother
\let\manuscriptmathbb\mathbb
\renewcommand{\mathbb}[1]{%
  \ifstrequal{#1}{k}{\Bbbk}{\manuscriptmathbb{#1}}%
}

\begin{document}

\begin{abstract}
In 2023, Fomin and Pylyavskyy introduced the \emph{master theorem} to find and prove linear incidence theorems using quadrilateral tilings of oriented surfaces. Each face of the surface corresponds to an incidence condition, which they call a \emph{coherence condition}. We extend their master theorem by incorporating polarity into coherence and allowing nonorientable surfaces. Our new coherence condition can be interpreted as perpendicularity from the perspective of Cayley--Klein geometry.

Moreover, a polarity incidence theorem on a nonorientable tiled surface gives rise to a self-polar linear incidence theorem on the orientation double cover. We also extend the master theorem to Hermitian polarities over division rings, with its validity depending on the topology of the surface. Finally, we determine which incidence presentations admit proofs using tilings. When such proofs exist, we express the minimum Euler genus of the corresponding tiling in terms of the fundamental group.
\end{abstract}

\maketitle

\section{Introduction}
\label{sec:introduction}

Linear incidence geometry studies configurations of points, lines, and higher-dimensional subspaces, with incidence conditions specifying containment and intersection relations among these objects. The subject can be traced back at least to ancient Greece, and generations of mathematicians have established numerous incidence theorems showing how prescribed incidence conditions imply further ones. Classical examples include the theorems of Desargues~\cite{BaralicEtAl2020}, Pappus~\cite{marchisotto}, and Möbius~\cite{moebius}, as well as the complete quadrangle~\cite{richter-gebert-mechanical} and bundle theorems~\cite{kahn}, among others.


Advances in computational commutative algebra have led to effective methods for the automated proof of incidence theorems~\cite{richter-gebert-mechanical,sturmfels1991,LiWu2001,hongbo-li,chou-gao-zhang}. These developments motivate us to seek structural explanations for incidence theorems beyond their individual proofs. Connections between incidence geometry and topology provide different perspectives on this structure, as illustrated by the work of Richter-Gebert~\cite{richter-gebert-ceva}, Glynn~\cite{Glynn2013RabbitHole}, Richter-Gebert and Li~\cite{RichterGebertLi2019}, and Barali{\'c}--Curien--Mili{\'c}evi{\'c}--Obradovi{\'c}--Petri{\'c}--Zeki{\'c}--{\v Z}ivaljevi{\'c}~\cite{BaralicEtAl2020}. In this direction, Fomin and Pylyavskyy introduced the \emph{master theorem} (\cite[Theorem~2.6]{FominPylyavskyy2023}) that provides a systematic framework for deriving linear incidence theorems from quadrilateral tilings of closed connected oriented surfaces. Within this framework, many classical incidence theorems can be recovered from suitable surface tilings, while other choices of tilings give rise to new theorems.

Subsequent work has investigated the scope and structure of tiling proofs. Pylyavskyy and Skopenkov~\cite{PylyavskyySkopenkov2026} study their dependence on the ground field and the underlying topological space, while K{\"u}hne and Larson~\cite{KuhneLarson2025} relate the master theorem to absolute incidence theorems. Katzenberger and Richter-Gebert~\cite{KatzenbergerRichterGebert2026} compare quadrilateral-tiling proofs with other methods of proving projective incidence theorems. In a related dynamical direction, Izosimov and Pylyavskyy~\cite{IzosimovPylyavskyy2025} connect incidence configurations with the dimer model.

We use polarity to extend this framework to incidences involving a conic and to tilings on nonorientable surfaces. We first state the result in the projective plane; the general construction is given in Section~\ref{sec:polarity-master}.

Let $V$ be a three-dimensional vector space over a field $k$, equipped with a nondegenerate symmetric bilinear form $\beta$. The \emph{polar line} of a point $P=[p]\in\mathbb P(V)$ is
\[
  P^\perp=\{[x]\in\mathbb P(V):\beta(p,x)=0\}.
\]
The \emph{pole} of a line $\ell$ is the unique point $\ell^\perp$ whose polar line is $\ell$. This point--line correspondence is the \emph{polarity} induced by $\beta$. Points $A=[a]$ and $B=[b]$ are \emph{conjugate} if $\beta(a,b)=0$, equivalently $A\in B^\perp$. Lines $\ell$ and $m$ are \emph{conjugate}, written $\ell\perp_\beta m$, if their poles are conjugate, equivalently $\ell^\perp\in m$. Both conjugacy relations are symmetric. Figures~\ref{fig:unitcir} and~\ref{fig:polarity-orthogonality} illustrate the unit-circle polarity and a pair of conjugate lines.

\begin{figure}[htbp]
\centering
\begin{minipage}[t]{0.41\textwidth}
\centering
\begin{minipage}[c][5.6cm][c]{\linewidth}
\centering
\begin{tikzpicture}[
    scale=1.8,
    every node/.style={font=\small},
    conic/.style={
        draw=green!55!black,
        line width=1.1pt
    },
    geom/.style={
        draw=black,
        line width=0.8pt
    },
    pt/.style={
        circle,
        fill=red!75!black,
        inner sep=1.7pt
    }
]

\coordinate (O) at (0,0);
\coordinate (P) at (2,0);

\coordinate (A) at (60:1);
\coordinate (B) at (-60:1);

\draw[conic] (O) circle (1);
\node[green!55!black, below left]
    at (225:1) {$\mathcal C$};

\draw[geom] (0.5,-1.3)--(0.5,1.3);
\node[above] at (0.5,1.3) {$P^\perp$};

\draw[geom] (A)--(P)--(B);

\node[pt] at (A) {};
\node[above left] at (A) {$A$};

\node[pt] at (B) {};
\node[below left] at (B) {$B$};

\node[pt] at (P) {};
\node[right] at (P) {$P$};

\end{tikzpicture}
\end{minipage}
\caption{Polarity with respect to the unit circle $\mathcal C$}
\label{fig:unitcir}
\end{minipage}\hfill
\begin{minipage}[t]{0.57\textwidth}
\centering
\begin{minipage}[c][5.6cm][c]{\linewidth}
\centering
\definecolor{polarityRed}{RGB}{199,45,75}
\definecolor{polarityCoral}{RGB}{239,94,94}
\definecolor{polarityTeal}{RGB}{47,133,122}
\begin{tikzpicture}[
    x=1.4cm,y=1.4cm,
    font=\small,
    red polar/.style={draw=polarityRed,line width=1pt},
    green polar/.style={draw=polarityTeal,line width=1pt},
    red join/.style={draw=polarityCoral,line width=.9pt},
    green join/.style={draw=polarityTeal,line width=.6pt},
    red label/.style={text=polarityRed,inner sep=1pt},
    green label/.style={text=polarityTeal,inner sep=1pt}
]
    \pgfmathsetmacro{\joinx}{1.35}
    \pgfmathsetmacro{\polex}{1/\joinx}
    \def\Ay{.72}
    \def\By{-.28}
    \def\Cy{-1.38}
    \def\Dy{.99}

    \begin{scope}[rotate=1.5]
        \coordinate (A) at (\joinx,\Ay);
        \coordinate (B) at (\joinx,\By);
        \coordinate (C) at (\polex,\Cy);
        \coordinate (D) at (\polex,\Dy);
        \coordinate (P) at (\polex,0); 
        \coordinate (Q) at (\joinx,0); 
    \end{scope}

    \begin{scope}
        \clip (-1.60,-1.65) rectangle (3.80,2.05);
        \draw[gray!70,line width=.7pt] (0,0) circle (1);
        \begin{scope}[rotate=1.5]
            \draw[red join] (\joinx,-5)--(\joinx,5);
            \draw[green join] (\polex,-5)--(\polex,5);

            \draw[red polar]
                (-3,{(1+3*\joinx)/\Ay})--(5,{(1-5*\joinx)/\Ay});
            \draw[red polar]
                (-3,{(1+3*\joinx)/\By})--(5,{(1-5*\joinx)/\By});
            \draw[green polar]
                (-3,{(1+3*\polex)/\Cy})--(5,{(1-5*\polex)/\Cy});
            \draw[green polar]
                (-3,{(1+3*\polex)/\Dy})--(5,{(1-5*\polex)/\Dy});
        \end{scope}
    \end{scope}

    \foreach \X in {A,B,P}
        \filldraw[fill=polarityCoral,draw=black,line width=.3pt]
            (\X) circle (1.5pt);
    \foreach \X in {C,D,Q}
        \filldraw[fill=polarityTeal,draw=black,line width=.3pt]
            (\X) circle (1.35pt);

    \node[red label,right=5pt] at (A) {$A$};
    \node[red label,below right=3pt] at (B) {$B$};
    \node[green label,below right=3pt] at (C) {$C$};
    \node[green label,left=5pt] at (D) {$D$};
    \node[red label,left=6pt] at (P) {$A^\perp\cap B^\perp$};
    \node[green label,right=6pt] at (Q) {$C^\perp\cap D^\perp$};

    \begin{scope}[rotate=1.5]
        \node[red label,above right=3pt]
            at (.13,{(1-.13*\joinx)/\Ay}) {$A^\perp$};
        \node[red label,left=5pt]
            at (.50,{(1-.50*\joinx)/\By}) {$B^\perp$};
        \node[green label,above left=3pt]
            at (2.43,{(1-2.43*\polex)/\Cy}) {$C^\perp$};
        \node[green label,above right=3pt]
            at (2.43,{(1-2.43*\polex)/\Dy}) {$D^\perp$};
    \end{scope}
\end{tikzpicture}
\end{minipage}
\caption{Conjugate lines $AB$ and $CD$, with poles $A^\perp\cap B^\perp$ and $C^\perp\cap D^\perp$.}
\label{fig:polarity-orthogonality}
\end{minipage}
\end{figure}

\begin{theorem}[Polarity master theorem]
\label{thm:polarity-master-intro}
Consider a quadrilateral tiling of a closed surface $S$, whose vertices are labeled by points of $\mathbb P(V)$. Assume that the labels at the endpoints of every edge are nonconjugate, and that a closed edge path preserves local orientation exactly when its length is even. If the two lines joining opposite labels are conjugate on every face but one, then the same holds on the remaining face. Here the condition is taken as automatic when either opposite pair of labels coincides.
\end{theorem}

We generalize the conjugacy conditions of a face to the \emph{polarity coherence} condition (Definition~\ref{def:polarity-tile}) and prove the corresponding master theorem (Theorem~\ref{thm:polarity-master}). We interpret general polarity coherence as perpendicularity (up to proper parallel translation) of geodesics in Cayley--Klein geometry (Theorem~\ref{thm:ck-orthogonality}).

The planar conjugacy condition has equivalent incidence formulations by the following \emph{polarity principle}.

\begin{proposition}[Polarity principle]
\label{prop:polarity_inv}
For every point $A$ and line $\ell$ in $\mathbb P(V)$,
\[
  A\in\ell \quad\Longleftrightarrow\quad \ell^\perp\in A^\perp.
\]
\end{proposition}

For a face with cyclic labels $A,C,B,D$, if $A\ne B$ and $C\ne D$, the coherence condition $AB\perp_\beta CD$ is equivalent to each of the following conditions:

\begin{equation}
\tag{*}\label{eq:coherence-intro}
\begin{aligned}
\text{(a)}\quad & A,\ B,\ C^\perp\cap D^\perp
  \text{ are collinear}; \\[3pt]
\text{(b)}\quad & AB,\ C^\perp,\ D^\perp
  \text{ are concurrent}; \\[3pt]
\text{(c)}\quad & C,\ D,\ A^\perp\cap B^\perp
  \text{ are collinear}; \\[3pt]
\text{(d)}\quad & CD,\ A^\perp,\ B^\perp
  \text{ are concurrent}.
\end{aligned}
\end{equation}

The Fomin--Pylyavskyy master theorem is recovered in the case of orientable surfaces (Section~\ref{subsec:polarity-tilings}). Our theorem also gives tiling proofs of the theorems of Hesse (Example~\ref{ex:hesse}), Pascal (Theorem~\ref{ex:pascal}), and Brianchon (Theorem~\ref{ex:brianchon}). It also yields a special case of the Cayley–Bacharach theorem (Theorem~\ref{ex:cb-line-conic}). The Hesse tiling below provides a simple nonorientable example.

\begin{example}[Hesse's theorem]
\label{ex:hesse}
Let $A,B,C,D$ be four points in $\mathbb{RP}^2$, equipped with the polarity induced by $\beta$. Assume that the joining lines below are well defined. Then
\[
    AB\perp_\beta CD
    \quad\text{and}\quad
    AC\perp_\beta BD
    \quad\Longrightarrow\quad
    AD\perp_\beta BC.
\]
Figures~\ref{fig:hesse-geometry} and~\ref{fig:hesse-tiling} show the geometric configuration and its tiling proof.

\begin{figure}[htbp]
\centering
\begingroup
\definecolor{hesseblue}{RGB}{36,99,154}
\definecolor{hessered}{RGB}{180,73,61}
\definecolor{hesseteal}{RGB}{45,130,115}

\begin{minipage}[t]{.49\textwidth}
\centering
\captionsetup{width=\linewidth}
\begin{minipage}[c][4.6cm][c]{\linewidth}
\centering
\begin{tikzpicture}[scale=.96, font=\small, line cap=round]
    \coordinate (A) at (2,0);
    \coordinate (B) at (0,2);
    \coordinate (C) at (-1,-1);
    \coordinate (D) at (0.25,0.25);
    \path[use as bounding box]
        (-2.15,-1.85) rectangle (2.6,2.55);

    \begin{scope}
        \clip (-2.05,-1.65) rectangle (2.35,2.25);
        \draw[gray!60, line width=.6pt]
            (0,0) circle[radius=1];

        \draw[hesseblue, line width=1pt]
            ($(A)!-3!(B)$)--($(A)!3!(B)$)
            ($(C)!-3!(D)$)--($(C)!3!(D)$);
        \draw[hessered, line width=1pt]
            ($(A)!-3!(C)$)--($(A)!3!(C)$)
            ($(B)!-3!(D)$)--($(B)!3!(D)$);
        \draw[hesseteal, line width=1pt]
            ($(A)!-3!(D)$)--($(A)!3!(D)$)
            ($(B)!-3!(C)$)--($(B)!3!(C)$);
    \end{scope}

    \node[gray!80, anchor=east, inner sep=1pt]
        at (-1.12,0.06) {$\mathcal C$};
    \foreach \V in {A,B,C,D}
        \fill (\V) circle[radius=1.5pt];
    \node[above right, inner sep=3pt] at (A) {$A$};
    \node[above right, inner sep=3pt] at (B) {$B$};
    \node[left, inner sep=3pt] at (C) {$C$};
    \node[below right, inner sep=3pt] at (D) {$D$};
\end{tikzpicture}
\end{minipage}
\caption{Hesse's configuration.}
\label{fig:hesse-geometry}
\end{minipage}\hfill
\begin{minipage}[t]{.49\textwidth}
\centering
\captionsetup{width=\linewidth}
\begin{minipage}[c][4.6cm][c]{\linewidth}
\centering
\begin{tikzpicture}[
    scale=.72,
    font=\small,
    line width=.8pt,
    every node/.style={inner sep=4pt},
    identification/.style={
        -{Latex[length=2mm]},
        line width=1.2pt
    }
]
\node (Bt) at (0,2.55) {$B$};
\node (Ct) at (-2,1.15) {$C$};
\node (Dt) at (2,1.15) {$D$};
\node (Db) at (-2,-1.15) {$D$};
\node (Cb) at (2,-1.15) {$C$};
\node (Bb) at (0,-2.55) {$B$};
\node (A)  at (0,0) {$A$};

\draw (Bt)--(Dt)--(Cb)--(Bb)--(Db)--(Ct)--(Bt);
\draw (Ct)--(A)--(Dt);
\draw (A)--(Bb);

\node at (0,1.10) {$F_1$};
\node at (-.95,-.70) {$F_2$};
\node at (.95,-.70) {$F_3$};

\draw[identification,red!75!black]
    ($(Bt)!.36!(Dt)$)--($(Bt)!.64!(Dt)$);
\draw[identification,red!75!black]
    ($(Bb)!.36!(Db)$)--($(Bb)!.64!(Db)$);

\draw[identification,blue!75!black]
    ($(Dt)!.36!(Cb)$)--($(Dt)!.64!(Cb)$);
\draw[identification,blue!75!black]
    ($(Db)!.36!(Ct)$)--($(Db)!.64!(Ct)$);

\draw[identification,green!50!black]
    ($(Ct)!.36!(Bt)$)--($(Ct)!.64!(Bt)$);
\draw[identification,green!50!black]
    ($(Cb)!.36!(Bb)$)--($(Cb)!.64!(Bb)$);
\end{tikzpicture}
\end{minipage}
\caption{Hesse's tiling of $\mathbb{RP}^2$.}
\label{fig:hesse-tiling}
\end{minipage}
\endgroup
\end{figure}
\end{example}

The orientation double cover of the Hesse tiling is the Desargues tiling (Example~\ref{ex:hesse-desargues-lift}). This illustrates how the cover relates polarity to projective duality (Section~\ref{sec:lifting-descent}). A self-dual incidence statement is invariant under projective duality, up to relabeling, while a self-polar realization assigns polar objects to vertices exchanged by the deck involution.

\begin{theorem}[Canonical lift]
\label{thm:canonical-lift}
Let $T$ be a realized polarity tiling of a quadrangulation $Q$ of a closed nonorientable surface satisfying the parity condition of Theorem~\ref{thm:polarity-master-intro}, with vertex labeling $\mathcal L$ and a marked face $f$. The orientation double cover of $Q$, marked at the two lifts of $f$, admits an FP tiling $\widetilde T$ carrying a valid self-dual incidence statement: under the coherence conditions on all other faces, the conditions on the two marked faces are equivalent. The original polarity induces a self-polar realization of $\widetilde T$ with respect to the deck involution, preserving face coherence.
\end{theorem}

We also characterize the involutions along which self-polar realizations descend (Theorem~\ref{thm:cover-characterization}).

We next study how surface topology governs the Hermitian extension over division rings (Section~\ref{sec:division-rings}). Following Izosimov~\cite{Izosimov2026}, we express face coherence as trivial holonomy and examine the boundary word after removing the conclusion face: a product of commutators for orientable surfaces, or of squares for nonorientable surfaces. Since squares can have nontrivial holonomy even over an abelian coefficient group, we use a semidirect product as the edge coefficient group, encoding orientation reversal together with Hermitian symmetry.

Let $\mathrm{D}$ be a division ring with an involutive anti-automorphism $\sigma$, and let $h$ be a nondegenerate $\sigma$-Hermitian form on a finite-dimensional right $\mathrm{D}$-vector space $V$, with $\dim_{\mathrm{D}} V\ge3$. Coherence in this setting is defined in Section~\ref{subsec:hermitian-master}. A tiling represents the polarity implication that coherence on all but one face forces coherence on the remaining face. We require the parity condition of Theorem~\ref{thm:polarity-master-intro}.

\begin{theorem}[Hermitian extensions of tiling proofs]
\label{thm:hermitian-extensions-intro}
If a polarity implication is represented by a tiling of the sphere, then it holds for arbitrary $\mathrm{D}$ and $\sigma$. If the tiling is on a closed orientable surface of positive genus, then it holds for any field $\mathrm{D}$ and arbitrary $\sigma$. If the tiling is on a closed nonorientable surface, then it holds for any field $\mathrm{D}$ with $\sigma=\id$. For a fixed tiling with a simple one-skeleton, these conditions are also necessary for universal validity when $\dim\mathbb P(V)$ satisfies the bound in Theorem~\ref{thm:hermitian-realization}.
\end{theorem}

The sufficient cases follow from Theorem~\ref{thm:hermitian-polarity-master}. The converse, proved in Section~\ref{subsec:hermitian-realization}, is obtained by realizing compatible edge pairings with distinct projective points for the fixed form $h$. This converts boundary-holonomy obstructions into geometric counterexamples.

Finally, for a formal presentation $I:Q_1,\ldots,Q_m\Longrightarrow Q_0$ (Definition~\ref{def:formal-presentation}), we determine the existence and minimum Euler genus of tiling proofs (Section~\ref{sec:tiling-surface}). These proofs must satisfy the parity condition above and use only the given labels and premise tiles. Following Zhang's construction \cite[Definition~3.1]{ZhangTiledGenus}, let $G_I$ be the fundamental group of the premise complex, based at a conclusion vertex, and $\delta_I$ the conclusion boundary class (Section~\ref{subsec:premise-complexes}).

An element of $G_I$ is even or odd according to the length modulo two of a representing edge loop. Let $\cleven_{G_I}$ denote commutator length with even factors, and $\slodd_{G_I}$ the least positive number of squares of odd elements, with value $\infty$ when the required factorization does not exist. Let $\eg(I)$ be the least Euler genus $2-\chi(S)$ of a closed surface supporting such a proof, with $\eg(I)=\infty$ if none exists. We obtain the following extension of Zhang's genus formula \cite[Theorem~1.9]{ZhangTiledGenus}.

\begin{theorem}[Euler-genus formula]\label{thm:surface-formula}
For every formal presentation $I$,
\begin{equation}\label{eq:euler-genus}
  \eg(I)
  =\min\left\{
    \slodd_{G_I}(\delta_I),
    2\cleven_{G_I}(\delta_I)
  \right\}.
\end{equation}
\end{theorem}

Even commutator length and odd square length are finite exactly when orientable and nonorientable tiling proofs exist, respectively (Theorem~\ref{thm:parity-tiling-criterion}); factorizations attaining the displayed minimum yield optimal tilings.

Section~\ref{sec:preliminaries} collects the background. Section~\ref{sec:polarity-master} proves the polarity master theorem and gives applications. Section~\ref{sec:ck-coherence} explains Cayley--Klein lengths and angles, interprets coherence as orthogonality, and discusses the Euclidean limit. Sections~\ref{sec:lifting-descent}, \ref{sec:division-rings}, and~\ref{sec:tiling-surface} treat lifting and descent, Hermitian polarities, and minimum Euler genus, respectively.


\subsection*{Acknowledgments}
This work was carried out as part of the PACE 2026 program at Peking University. We thank our mentor, Yuchong Zhang, for proposing this project and for his continuous support and guidance throughout the program. We are grateful to Yibo Gao, the organizer of PACE 2026, for creating the platform that brought us together and for his valuable advice. We also thank Jidong Wang and Tao Gui for their mentorship and insightful discussions. We are also grateful to Sergey Fomin, Anton Izosimov, Matt Larson and Mikhail Skopenkov for taking the time to read our manuscript and for sharing their thoughtful comments and suggestions.

\subsection*{AI Disclosure}
During the preparation of this paper, we used AI tools to help prepare illustrative figures, identify relevant references, and explore examples and counterexamples with us, including the Klein bottle polarity tiling in Example~\ref{ex:klein-bottle-tiling} and the failure of the complex Hermitian Hesse analogue in Example~\ref{ex:hermitian-hesse}. We also used AI tools to improve the exposition and correct grammatical and typographical errors. All AI-assisted content was reviewed and verified by the authors, who take full responsibility for the manuscript.

\section{Preliminaries}
\label{sec:preliminaries}

Since our result builds on the master theorem of Fomin and Pylyavskyy, we first recall their master theorem over a commutative field. We then fix the projective and topological conventions used throughout the paper, allowing division rings in the geometric setting.

\subsection{The Fomin--Pylyavskyy Master Theorem}
\label{subsec:fp-master}

\begin{theorem}[\cite{FominPylyavskyy2023}, Theorem~2.6]
\label{thm:fp-master}
\normalfont

Consider a tiling of a closed oriented surface by quadrilateral~tiles. Assume that the vertices of the tiling are colored black and white, so that every edge connects vertices of different color. Associate to each black (resp., white) vertex a~point (resp., a hyperplane) in a finite-dimensional projective space $\mathbb P$, with $\dim\mathbb P\geq 2$. For every edge joining labels $A$ and $h$, require $A\notin h$.

For each tile, write its labels in cyclic order as

\begin{center}
\begin{tikzpicture}[scale=1]
    \node at (0,1.5) {$A$};
    \node at (1.6,1.5) {$\ell$};
    \node at (0,0) {$m$};
    \node at (1.6,0) {$B$};

    \draw (0.28,1.5) -- (1.3,1.5);
    \draw (0,1.15) -- (0,0.28);
    \draw (0.25,0) -- (1.3,0);
    \draw (1.6,1.15) -- (1.6,0.28);
\end{tikzpicture}
\end{center}

where $A$ and $B$ are points and $\ell$ and $m$ are hyperplanes. Such a tile is called \emph{coherent} if either $A=B$, or $\ell=m$, or else the line $(AB)$ and the codimension-two subspace $\ell\cap m$ have a nonempty intersection. If all tiles but one are coherent, then the remaining tile is coherent as well. In particular, when $\mathbb P$ is a projective plane, hyperplanes are lines, and coherence means that the points $A$, $B$, and $\ell\cap m$ are collinear.

\end{theorem}

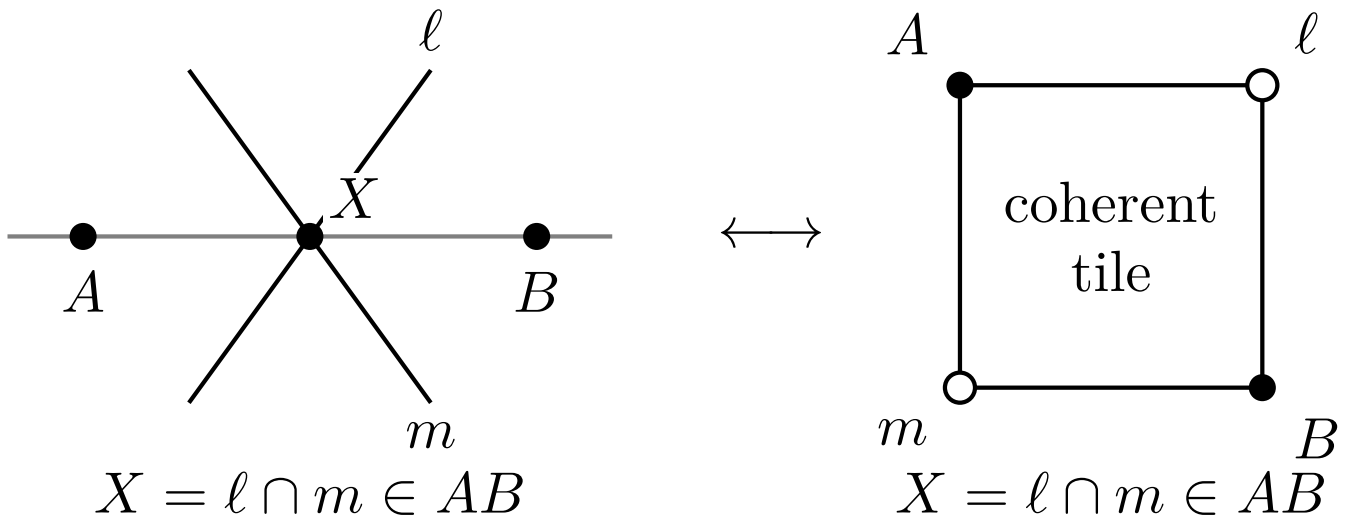
\begin{figure}[htbp]
\centering
\begin{tikzpicture}[
    font=\small,
    thick,
    pv/.style={circle,fill=black,inner sep=1.8pt},
    lv/.style={circle,draw=black,fill=white,inner sep=2pt}
]

\begin{scope}
    \draw[gray] (-2,0)--(2,0);
    \draw (-.8,-1.1)--(.8,1.1) node[above] {$\ell$};
    \draw (-.8,1.1)--(.8,-1.1) node[below] {$m$};

    \node[pv,label=below:$A$] at (-1.5,0) {};
    \node[pv,label=below:$B$] at (1.5,0) {};
    \node[
        pv,
        label={[fill=white,inner sep=1pt]above right:$X$}
    ] at (0,0) {};

    \node at (0,-1.7) {$X=\ell\cap m\in AB$};
\end{scope}

\node at (3.05,0) {$\longleftrightarrow$};

\begin{scope}[xshift=5.3cm]
    \draw (-1,1)--(1,1)--(1,-1)--(-1,-1)--cycle;

    \node[pv,label=above left:$A$] at (-1,1) {};
    \node[lv,label=above right:$\ell$] at (1,1) {};
    \node[pv,label=below right:$B$] at (1,-1) {};
    \node[lv,label=below left:$m$] at (-1,-1) {};

    \node[align=center] at (0,0) {coherent\\tile};
    \node at (0,-1.7) {$X=\ell\cap m\in AB$};
\end{scope}

\end{tikzpicture}
\caption{The left panel shows the condition $\ell\cap m\in AB$. The right panel represents the same condition by a quadrilateral tile. Black vertices represent points, and white vertices represent lines.}
\label{fig:coherent-tile}
\end{figure}

Figure~\ref{fig:coherent-tile} shows an example of a \emph{coherent tile} in dimension two, where hyperplanes are lines in the projective plane.

We express this formulation algebraically as follows. Assume that~$\mathbb{P}:=\mathbb{P}(V)$ is the projectivization of a vector space~$V$, so that points in $\PP$ are identified with one-dimensional subspaces in~$V$.
Hyperplanes in~$\PP$ correspond to codimension~1 subspaces of~$V$, so they can be identified with one-dimensional subspaces in the dual space~$V^*$.

For a vector $\mathbf{V}\!\in\!V$ and a covector $\mathbf{h}\!\in\!V^*$, let $\langle \mathbf{V},\mathbf{h}\rangle$ denote their pairing.

Here for points $A,B\!\in\!\PP$ and lines $m,\ell \in\PP^*$, let $\mathbf{A}, \mathbf{B}\!\in\!V$ be the corresponding vectors and $ \mathbf{m}, \mathbf{\boldsymbol \ell}\in V^*$ be the corresponding covectors. Then we can define the \emph{mixed cross-ratio} by
\[
    \chi(A,B;m,\ell)
    =
    \frac{\langle \mathbf{A},\mathbf{m}\rangle\langle \mathbf{B},\mathbf{\boldsymbol \ell}\rangle}
         {\langle \mathbf{B},\mathbf{m}\rangle\langle \mathbf{A},\mathbf{\boldsymbol \ell}\rangle}.
\]
We use this mixed cross-ratio as the weight of the tile and it is well-defined (see~\cite{coxeter-non-euclidean, coxeter-projective}).

\begin{proposition}
\label{prop:coherent-iff-weight-one}
Let $A,B\in\PP$ be points and $m,\ell\in\PP^*$ be hyperplanes such that neither point lies on either hyperplane. The tile whose vertices are cyclically labeled $A,\ell,B,m$ is coherent if and only if
\[
    \chi(A,B;m,\ell)=1.
\]
\end{proposition}

\begin{proof}
Assume that the tile is coherent. Then the line $AB$ passes through a point lying in $m\cap\ell$. Hence $\tfrac{\langle\mathbf A,\mathbf m\rangle} {\langle\mathbf A,\boldsymbol{\ell}\rangle} = \tfrac{\langle\mathbf B,\mathbf m\rangle} {\langle\mathbf B,\boldsymbol{\ell}\rangle}$, so $\chi(A,B;m,\ell)=1$.

Conversely, assume that the tile is not coherent. Then the points $M=AB\cap m$ and $L=AB\cap\ell$ are distinct. Hence $\tfrac{\langle\mathbf A,\mathbf m\rangle} {\langle\mathbf A,\boldsymbol{\ell}\rangle} \neq \tfrac{\langle\mathbf B,\mathbf m\rangle} {\langle\mathbf B,\boldsymbol{\ell}\rangle}$, so $\chi(A,B;m,\ell)\neq1$.
\end{proof}

\begin{proof}[Proof of Theorem~\ref{thm:fp-master}]
Fix an orientation of $S$ and give each tile the induced boundary orientation. Compute each weight $\chi(f)$ in this boundary order, starting at a point vertex.

Choose a nonzero vector $\mathbf P\in V$ for each point label $P$, and a nonzero covector $\mathbf h\in V^*$ for each hyperplane label $h$. The induced boundary orientations traverse each shared edge in opposite directions. Thus, for every edge $P-h$, the pairing $\langle\mathbf P,\mathbf h\rangle$ appears in the numerator of one adjacent tile's weight and in the denominator of the other's. Since $S$ is closed, every edge contributes exactly two such factors, which cancel. Hence
\[
    \prod_{\text{tiles }f}\chi(f)=1.
\]

Suppose that every tile except $f_0$ is coherent. By Proposition~\ref{prop:coherent-iff-weight-one}, $\chi(f)=1$ for every $f\neq f_0$. The product identity therefore gives $\chi(f_0)=1$, so $f_0$ is coherent.
\end{proof}

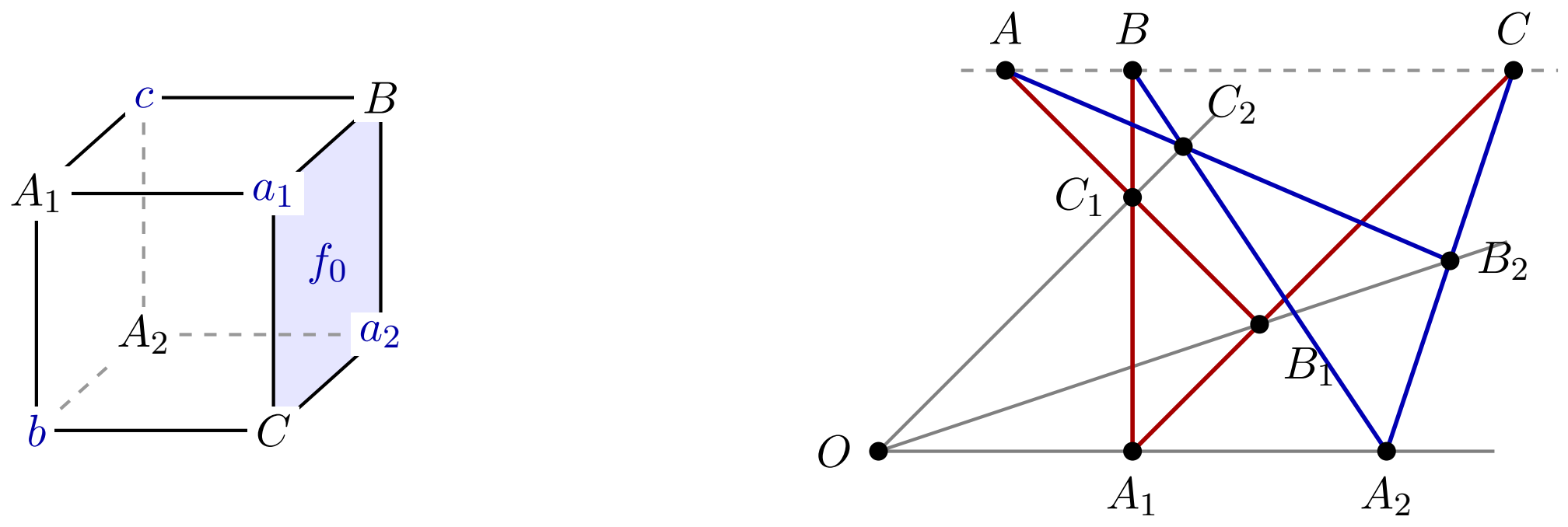
\begin{figure}[htbp]
\centering
\begin{minipage}[c]{.37\linewidth}
\centering
\begin{tikzpicture}[
    scale=.96,
    font=\small,
    line width=.8pt,
    lab/.style={fill=white,inner sep=2pt}
]
\coordinate (F1) at (0,2.1);
\coordinate (F2) at (2.1,2.1);
\coordinate (F3) at (2.1,0);
\coordinate (F4) at (0,0);

\coordinate (H1) at (.95,2.95);
\coordinate (H2) at (3.05,2.95);
\coordinate (H3) at (3.05,.85);
\coordinate (H4) at (.95,.85);

\fill[blue!10] (F2)--(H2)--(H3)--(F3)--cycle;

\draw[dashed,gray!80]
    (H1)--(H4)--(H3)
    (H4)--(F4);

\draw (F1)--(F2)--(F3)--(F4)--cycle;
\draw
    (F1)--(H1)--(H2)--(H3)--(F3)
    (F2)--(H2);

\node[lab] at (F1) {$A_1$};
\node[lab,text=blue!65!black] at (F2) {$a_1$};
\node[lab] at (F3) {$C$};
\node[lab,text=blue!65!black] at (F4) {$b$};

\node[lab,text=blue!65!black] at (H1) {$c$};
\node[lab] at (H2) {$B$};
\node[lab,text=blue!65!black] at (H3) {$a_2$};
\node[lab] at (H4) {$A_2$};

\node[text=blue!65!black] at (2.58,1.48) {$f_0$};
\end{tikzpicture}
\end{minipage}\hfill
\begin{minipage}[c]{.60\linewidth}
\centering
\begin{tikzpicture}[
    scale=1.08,
    font=\small,
    line width=.8pt,
    dot/.style={circle,fill=black,inner sep=1.6pt},
    first/.style={red!65!black,line width=1.05pt},
    second/.style={blue!70!black,line width=1.05pt}
]
\coordinate (O) at (-2,0);

\coordinate (A1) at (0,0);
\coordinate (B1) at (1,1);
\coordinate (C1) at (0,2);

\coordinate (A2) at (2,0);
\coordinate (B2) at (2.5,1.5);
\coordinate (C2) at (.4,2.4);

\coordinate (X) at (-1,3);
\coordinate (Y) at (0,3);
\coordinate (Z) at (3,3);

\draw[gray]
    (O)--(2.85,0)
    (O)--(2.95,1.65)
    (O)--(.65,2.65);

\draw[first]
    (X)--(B1)--(A1)--(Y)
    (B1)--(Z);

\draw[second]
    (X)--(B2)--(A2)--(Y)
    (B2)--(Z);

\draw[dashed,gray!85] (-1.35,3)--(3.35,3);

\node[dot,label=left:$O$] at (O) {};

\node[dot,label=below:$A_1$] at (A1) {};
\node[dot,label=below right:$B_1$] at (B1) {};
\node[dot,label=left:$C_1$] at (C1) {};

\node[dot,label=below:$A_2$] at (A2) {};
\node[dot,label=right:$B_2$] at (B2) {};
\node[dot,label=above right:$C_2$] at (C2) {};

\node[dot,label=above:$A$] at (X) {};
\node[dot,label=above:$B$] at (Y) {};
\node[dot,label=above:$C$] at (Z) {};
\end{tikzpicture}
\end{minipage}
\caption{Desargues' theorem and its cube tiling. The shaded face gives the conclusion.}
\label{fig:desargues-tiling}
\end{figure}

\begin{example}[Desargues' theorem]
\label{ex:desargues}

Let $a,b,c$ be distinct concurrent lines on the complex or real projective plane. Pick generic points $A_1,A_2\in a$, $B_1,B_2\in b$, and $C_1,C_2\in c$. Then the points
\[
A=(B_1C_1)\cap(B_2C_2),\qquad
B=(A_1C_1)\cap(A_2C_2),\qquad
C=(A_1B_1)\cap(A_2B_2)
\]
are collinear.

This is a direct application of the Fomin--Pylyavskyy master theorem. Write $a_i=(B_iC_i)$ for $i=1,2$. Figure~\ref{fig:desargues-tiling} represents the hypotheses and conclusion of Desargues' theorem as coherence conditions on the six faces of a cube tiling of the sphere. Five faces are coherent by construction, so the master theorem forces the shaded face to be coherent. Since the condition expressed by this face is $ a_1 \cap a_2 \in BC $, and since $a_1 \cap a_2 = A$, the coherence of this face is exactly the collinearity of $A,B,C$.
\end{example}

To extend this construction, we fix a polarity on the projective space. This allows us to replace a point label by its polar hyperplane, or a hyperplane label by its pole, and thus pass between point and hyperplane labels in the tiling. We next fix the projective conventions and define the polarities used in the paper.

\subsection{Projective Spaces and Polarities}
\label{subsec:projective-polarities}

Let $\mathrm{D}$ be a division ring, and let $V$ be a right $\mathrm{D}$-vector space of dimension $n+1$, where $n\ge2$. We write $\mathbb P(V)$ for the projective space of one-dimensional right subspaces of $V$, and put $\mathbb P^n(\mathrm{D})=\mathbb P(\mathrm{D}^{n+1})$.

\begin{definition}[]
\label{def:projective-points-lines-incidence}
A \emph{projective point} is a right subspace $A=[a]=a\mathrm{D}$, where $a\in V\setminus\{0\}$. Its vector representatives are $ad$, with $d\in \mathrm{D}^\times$. For a nonzero right subspace $U\subseteq V$, the corresponding projective subspace $\mathbb P(U)$ has dimension $\dim_{\mathrm{D}} U-1$. \emph{Incidence} is containment of the corresponding subspaces. In particular, a \emph{projective line} comes from a two-dimensional right subspace, and the join of distinct points $A=[a]$ and $B=[b]$ is $AB=\mathbb P(a\mathrm{D}+b\mathrm{D})$.
\end{definition}

The dual $V^*=\operatorname{Hom}_{\mathrm{D}}(V,\mathrm{D})$ consists of right-linear functionals. It is naturally a \emph{left} $\mathrm{D}$-vector space, with $(c\alpha)(v)=c\alpha(v)$. Its one-dimensional left subspaces $\mathrm{D}\alpha$ parameterize the projective hyperplanes $H_\alpha=\mathbb P(\ker\alpha)$, where $\alpha\ne0$. Thus the dual projective space is the left projective space of $V^*$. This distinction fixes the side on which representatives may be rescaled, and will matter in the division-ring setting.

For $A=[a]$ and $H=H_\alpha$, we use the pairing notation $\langle A,H\rangle=\alpha(a)$, with representatives understood. The scalar depends on these representatives, since $(c\alpha)(ad)=c\,\alpha(a)\,d$ for $c,d\in \mathrm{D}^\times$, but its vanishing is well defined: $A\in H$ if and only if $\alpha(a)=0$. When $\mathrm{D}=k$ is a field, left and right scalar actions agree, and we use the usual notation $\mathbb P(V^*)$ for the dual projective space.

To introduce a polarity, we require additional structure. Assume that $\mathrm{D}$ has an involutive anti-automorphism $\sigma$, so that $\sigma(cd)=\sigma(d)\sigma(c)$ and $\sigma^2=\id$, and fix a nondegenerate $\sigma$-Hermitian form $h:V\times V\to \mathrm{D}$. Our convention is $h(xc,yd)=\sigma(c)h(x,y)d$ and $h(y,x)=\sigma(h(x,y))$. These are the standard conventions for Hermitian forms on right vector spaces \cite[Chapter~8]{Scharlau1985}.

\begin{definition}[]
\label{def:orthogonal-complement}
\label{def:projective-polarity}
For a right subspace $U\subseteq V$, set
\[
 U^\perp=\{x\in V:h(u,x)=0\text{ for every }u\in U\}.
\]
We also use $\perp$ for the induced operation on projective subspaces. The polar of a point $A$ is the hyperplane $A^\perp$, and the pole of a hyperplane $H$ is the point $H^\perp$. Two points $[a]$ and $[b]$ are \emph{conjugate} if $h(a,b)=0$.
\end{definition}

Nondegeneracy gives $\dim_{\mathrm{D}} U^\perp=n+1-\dim_{\mathrm{D}} U$ and $(U^\perp)^\perp=U$. Consequently the polarity is involutive and satisfies reciprocity: $A\in H\Longleftrightarrow H^\perp\in A^\perp$. In particular, point conjugacy is symmetric and does not depend on the chosen vector representatives.

The polarity constructions in Sections~\ref{sec:polarity-master} and \ref{sec:lifting-descent} use $\mathrm{D}=k$, $\sigma=\id$, and $h=\beta$, where $\beta$ is a fixed nondegenerate symmetric bilinear form on $V$. In the projective plane, hyperplanes are lines. Two lines $L,M$ are \emph{$\beta$-conjugate}, written $L\perp_\beta M$, if $L^\perp\in M$, equivalently $M^\perp\in L$. For distinct points $A,B$, the pole of their join is $(AB)^\perp=A^\perp\cap B^\perp$. Thus, when both joins are defined, $AB\perp_\beta CD$ means $(AB)^\perp\in CD$. This last formulation uses the plane setting: in higher dimension $(AB)^\perp$ is a projective subspace, and the analogous intersection condition is $(AB)^\perp\cap CD\ne\varnothing$.

A \emph{projective realization} assigns geometric objects to formal labels. Distinct labels may have the same realization unless this is explicitly excluded. The division ring, ambient dimension, and, when used, the Hermitian form are part of the geometric setting; each result specifies which of these are fixed. Vectors and covectors serve only as representatives of the realized projective objects.

\subsection{Topological Conventions}
\label{subsec:topological-conventions}

Surfaces are compact and connected unless stated otherwise, and may have boundary. All cell structures are finite.

\begin{definition}[]
\label{def:quadrangular-complex}
\label{def:quadrangulation}
A \emph{quadrangular complex} is a regular two-dimensional CW complex whose closed two-cells are disks bounded by four distinct edges and four distinct vertices. Its vertices, edges, and faces are denoted by $V(Q),E(Q),F(Q)$, respectively. A \emph{quadrangulation} of a surface $S$ is such a complex with underlying space $S$.
\end{definition}

Regularity means that each characteristic map embeds its closed cell. Intersections of closed cells may be unions of common cells; two faces may therefore share more than one edge. The face structure is part of the data: a four-cycle in $Q^{(1)}$ need not bound a face. The premise complexes introduced below are more general CW complexes and need not be regular or be surfaces.

A \emph{tiling} $T$ of a quadrangular complex $Q$ consists of a formal vertex labeling $\mathcal L:V(Q)\to\Lambda$. We reserve $Q$ for the underlying complex and $T$ for the tiling with its labeling and any chosen realization. Each face of $Q$, equipped with its vertex labels, is called a \emph{tile} of $T$.

For an \emph{FP tiling}, the label set is partitioned into point labels and hyperplane labels, and every edge joins vertices of different types. The types define a two-coloring $c$, with value $0$ at point vertices and $1$ at hyperplane vertices. A \emph{realization} $R$ assigns projective objects to the formal labels; for an FP tiling, it respects their types and makes adjacent point and hyperplane labels nonincident. A tiling with a chosen realization is called \emph{realized}. We write $R(v)$ for $R(\mathcal L(v))$ when $v$ is a vertex. Face coherence always refers to a chosen realization.

\begin{definition}[]
\label{def:orientation-character}
\label{def:orientation-double-cover}
The \emph{orientation class} of a surface $S$ is $w_{1,S}\in H^1(S;\Ftwo)$. Viewed as a character of $\pi_1(S)$, it is given by
\[
 w_{1,S}([\gamma])=
 \begin{cases}
  0,&\gamma\text{ preserves local orientation},\\
  1,&\gamma\text{ reverses local orientation}.
 \end{cases}
\]
The \emph{orientation double cover} $\pi:\widetilde S\to S$ has as its fiber over each point the two possible local orientations at that point.
\end{definition}

The cover carries its tautological orientation. Its deck involution $\tau$ exchanges the two local orientations and reverses orientation on $\widetilde S$. It is connected exactly when $S$ is nonorientable; otherwise it is the disjoint union of two copies of $S$. A quadrangulation lifts to $\widetilde Q=\pi^{-1}(Q)$, with two lifts of each cell exchanged by $\tau$.

The second character we use records the parity of edge-path length. Unlike the orientation character, it depends on the chosen cell structure.

\begin{definition}[]
\label{def:parity-character}
Let $K$ be a connected CW complex of dimension at most two whose two-cell attaching maps are edge-paths of even length. Let $\lambda_K$ be the cellular $1$-cochain taking the value $1\in\Ftwo$ on every edge. Its cohomology class defines the \emph{parity character}
\[
 \varepsilon_K:\pi_1(K)\longrightarrow\Ftwo,
 \qquad \varepsilon_K([\gamma])=|\gamma|\pmod2,
\]
where $|\gamma|$ is the length of a closed edge-path representative.
\end{definition}

Indeed, $d\lambda_K=0$, since each two-cell boundary has even length, with edge occurrences counted with multiplicity. This also shows that the formula is independent of the edge-path representative. For a quadrangulated surface, the condition $[\lambda_Q]=w_{1,S}$ means that a closed edge-path has odd length exactly when it reverses local orientation. We shall use this as a compatibility condition; it is not an automatic property of a quadrangulation.

If a cellular map $F:K\to K'$ between complexes of the above kind maps each edge homeomorphically to an edge, then $F^*\lambda_{K'}=\lambda_K$ and $\varepsilon_{K'}\circ F_*=\varepsilon_K$. Thus parity behaves naturally under the maps of tiled surfaces used later in the paper.

\section{Polarity Master Theorem}
\label{sec:polarity-master}

\subsection{Polarity Tiles and Tilings}
\label{subsec:polarity-tiles}

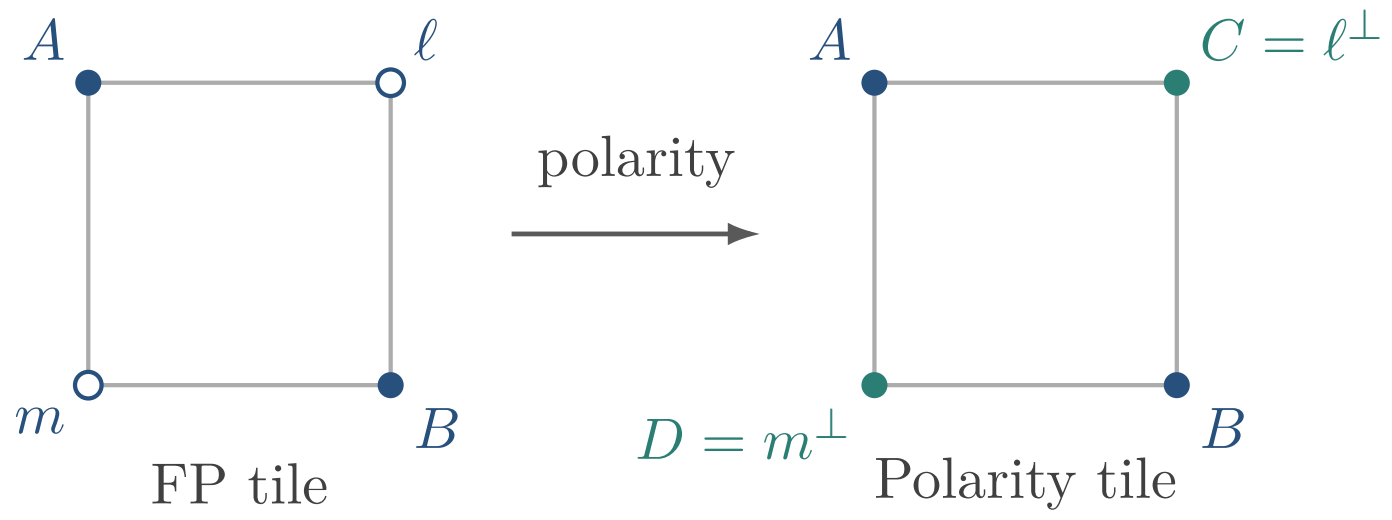
\begin{figure}[htbp]
\centering
\begingroup
\definecolor{conversionblue}{RGB}{40,80,125}
\definecolor{conversionteal}{RGB}{42,126,116}
\begin{tikzpicture}[
    font=\small,
    edge/.style={draw=gray!65,line width=.8pt},
    point/.style={circle,draw=conversionblue,fill=conversionblue,
        minimum size=4.2pt,inner sep=0pt,line width=.7pt},
    line/.style={circle,draw=conversionblue,fill=white,
        minimum size=5pt,inner sep=0pt,line width=.8pt},
    pole/.style={circle,draw=conversionteal,fill=conversionteal,
        minimum size=4.2pt,inner sep=0pt,line width=.7pt},
    label/.style={inner sep=2pt},
    title/.style={font=\small,text=black!75}
]
    \draw[edge] (-1,1)--(1,1)--(1,-1)--(-1,-1)--cycle;
    \node[point] at (-1,1) {};
    \node[line]  at (1,1) {};
    \node[point] at (1,-1) {};
    \node[line]  at (-1,-1) {};

    \node[label,above left=3pt,text=conversionblue] at (-1,1) {$A$};
    \node[label,above right=3pt,text=conversionblue] at (1,1) {$\ell$};
    \node[label,below right=3pt,text=conversionblue] at (1,-1) {$B$};
    \node[label,below left=3pt,text=conversionblue] at (-1,-1) {$m$};
    \node[title] at (0,-1.65) {FP tile};

    \draw[-{Latex[length=2.2mm,width=1.5mm]},black!65,line width=.9pt]
        (1.8,0)--(3.45,0)
        node[midway,above=5pt,text=black!75] {polarity};

    \begin{scope}[xshift=5.2cm]
        \draw[edge] (-1,1)--(1,1)--(1,-1)--(-1,-1)--cycle;
        \node[point] at (-1,1) {};
        \node[pole]  at (1,1) {};
        \node[point] at (1,-1) {};
        \node[pole]  at (-1,-1) {};

        \node[label,above left=3pt,text=conversionblue] at (-1,1) {$A$};
        \node[label,above right=3pt,text=conversionteal] at (1,1)
            {$C=\ell^\perp$};
        \node[label,below right=3pt,text=conversionblue] at (1,-1) {$B$};
        \node[label,below left=3pt,text=conversionteal] at (-1,-1)
            {$D=m^\perp$};
        \node[title] at (0,-1.65) {Polarity tile};
    \end{scope}
\end{tikzpicture}
\endgroup
\caption{Replacing line labels by their poles preserves coherence.}
\label{fig:fp-polarity-conversion}
\end{figure}

Throughout this section, we take $\mathrm{D}=k$ to be a field and fix a nondegenerate symmetric bilinear form $\beta$ on $V$. We call the quadrilateral tiles of Section~\ref{subsec:fp-master} \emph{FP tiles}. Polarity allows us to replace hyperplane labels by their poles. An FP tile with labels $A,\ell,B,m$ becomes a quadrilateral with point realizations $A,C,B,D$, where $C=\ell^\perp$ and $D=m^\perp$; see Figure~\ref{fig:fp-polarity-conversion}.

\begin{definition}
\label{def:polarity-tile}
A \emph{polarity tile} $(A,B;C,D)$ is a quadrilateral with point labels $A,C,B,D\in\mathbb P(V)$ in cyclic order. The labels need not be distinct, but adjacent labels are required to be nonconjugate. If $A\ne B$ and $C\ne D$, the tile is \emph{polarity coherent}, or simply \emph{coherent} when $AB\cap(CD)^\perp\ne\varnothing$. If $A=B$ or $C=D$, the tile is also called coherent. For $C\ne D$, the \emph{polar subspace} $(CD)^\perp=C^\perp\cap D^\perp$ is a point in the projective plane. For a polarity tile, choose nonzero representatives $a,b,c,d$ and define
\begin{equation}
\label{eq:fp-polarity-weight}
 \chi_\beta(A,B;C,D)
 :=\frac{\beta(a,c)\beta(b,d)}{\beta(a,d)\beta(b,c)}.
\end{equation}
The ratio is independent of the representatives.
\end{definition}

\begin{proposition}
\label{prop:polarity-coherence-weight}
A polarity tile is coherent if and only if $\chi_\beta(A,B;C,D)=1$.
\end{proposition}

\begin{proof}
If $A=B$ or $C=D$, the equality holds automatically. Otherwise, $\chi_\beta=1$ is equivalent to $\beta(a,c)\beta(b,d)=\beta(a,d)\beta(b,c)$. This holds exactly when there is a nonzero vector $x\in\operatorname{span}\{a,b\}$ such that $\beta(x,c)=\beta(x,d)=0$. Equivalently, the line $AB$ meets $(CD)^\perp$.
\end{proof}

For an FP tile with $C=\ell^\perp$ and $D=m^\perp$, we have $\chi_\beta(A,B;C,D)=\chi(A,B;\ell,m)$. By Propositions~\ref{prop:coherent-iff-weight-one} and~\ref{prop:polarity-coherence-weight}, the conversion preserves coherence.

\begin{definition}
\label{def:polarity-labeling}
A \emph{polarity tiling} $T$ of a quadrangular complex $Q$ has a labeling $\mathcal L:V(Q)\to\Lambda$ whose formal labels all denote point variables. A realization of $T$ is a map $R:\Lambda\to\mathbb P(V)$ such that adjacent vertices have nonconjugate realizations.

Each face of $Q$, together with the realizations of its vertex labels, forms a polarity tile. A face is \emph{coherent} in $T$ if the corresponding polarity tile is coherent.
\end{definition}

We now turn to the topology of the underlying quadrangulation. Before stating the master theorem for polarity tilings, we introduce the condition under which it applies.

\subsection{Orientation-Twisted Quadrangulation}

For a quadrangulation $Q$ of $S$, let $\lambda:=\lambda_Q$ be the cochain from Definition~\ref{def:parity-character}, and write $\ell:=[\lambda]\in H^1(S;\Ftwo)$. A function $c\in C^0(Q;\Ftwo)$ is a two-coloring of the one-skeleton if and only if $dc=\lambda$.

\begin{definition}[]\label{def:orientation-twisted}
The quadrangulation $Q$ of $S$ is \emph{orientation-twisted} if $\ell=w_{1,S}.$
\end{definition}

Intuitively, a quadrangulation is orientation-twisted if and only if every orientation-preserving closed edge path has even length and every orientation-reversing closed edge path has odd length. In particular, when $S$ is orientable, orientation-twisted means precisely that $Q^{(1)}$ is bipartite. We use the standard characteristic-class and orientation-cover facts recorded in \cite{MilnorStasheff1974,Hatcher2002}.

Recall that the first Stiefel--Whitney class $w_{1}$ is the obstruction to orientability, which evaluates to zero on orientation-preserving loops and to one on orientation-reversing loops. We can construct it in $Q$ as follows. Let $S=|Q|$, and let $p:\operatorname{Or}(S)\to S$ be the orientation double cover with deck involution $\tau$. Choose one lift $\widetilde v\in p^{-1}(v)$ for every vertex $v$. For an oriented edge $e:u\to v$, lift $e$ from $\widetilde u$ and define $\varepsilon(e)\in\mathbb F_2$ by requiring its endpoint to be $\tau^{\varepsilon(e)}\widetilde v$. Since the boundary of every face is null-homotopic, $\mathrm d\varepsilon=0$. Replacing $\widetilde v$ by $\tau^{s(v)}\widetilde v$ changes $\varepsilon$ to $\varepsilon+\mathrm ds$; hence $[\varepsilon]$ is independent of all choices. Its evaluation on a closed edge path records whether parallel transport preserves or reverses local orientation. Therefore $[\varepsilon]=w_1(S)$.

\begin{proposition}[]
\label{prop:cover-coloring}
Let $Q$ be a quadrangulation of a closed connected nonorientable surface $S$. Let $\pi:\widetilde S\to S$ be the orientation double cover, let $\widetilde Q=\pi^{-1}(Q)$, and let $\tau$ be the deck involution. Then $Q$ is orientation-twisted if and only if $\widetilde Q$ is bipartite and admits a two-coloring $c:V(\widetilde Q)\to\Ftwo$ satisfying $c\circ\tau=c+1$.
\end{proposition}

\begin{proof}
Let $\lambda\in Z^1(Q;\Ftwo)$ be the edge-length cocycle, $\lambda(e)=1$, and write $\ell=[\lambda]$.

Suppose first that $Q$ is orientation-twisted, so that $\ell=w_{1,S}$. Naturality of the first Stiefel--Whitney class gives $\pi^*\ell=\pi^*w_{1,S}=w_{1,\widetilde S}=0$. Hence $\pi^*\lambda=dc$ for some $c\in C^0(\widetilde Q;\Ftwo)$. Since $\pi^*\lambda(\widetilde e)=1$ for every edge $\widetilde e$ of $\widetilde Q$, the function $c$ changes value across every edge and therefore defines a two-coloring of $\widetilde Q$.

Choose a path $\widetilde\gamma$ from a vertex $\widetilde v$ to $\tau\widetilde v$. Its projection $\gamma=\pi\circ\widetilde\gamma$ is an orientation-reversing closed edge path in $S$. Therefore $|\widetilde\gamma|=|\gamma| =\langle\ell,[\gamma]\rangle =\langle w_{1,S},[\gamma]\rangle=1\pmod 2$. Since $c$ changes along every edge, it follows that $c(\tau\widetilde v)=c(\widetilde v)+1$. Moreover, $d(c\circ\tau+c)=\tau^*(dc)+dc =\tau^*(\pi^*\lambda)+\pi^*\lambda=0$. Thus $c\circ\tau+c$ is constant on the connected graph $\widetilde Q^{(1)}$. Since its value at $\widetilde v$ is $1$, we obtain $c\circ\tau=c+1$ at every vertex.

Conversely, suppose that $\widetilde Q$ admits a two-coloring $c$ satisfying $c\circ\tau=c+1$. Let $\gamma$ be any closed edge path in $Q$, and lift it to a path $\widetilde\gamma$ beginning at $\widetilde v$. If $\gamma$ is orientation-preserving, then $\widetilde\gamma$ ends at $\widetilde v$, so its length is even. If $\gamma$ is orientation-reversing, then $\widetilde\gamma$ ends at $\tau\widetilde v$; since these endpoints have opposite colors, its length is odd. Consequently, $\langle\ell,[\gamma]\rangle=|\gamma| =\langle w_{1,S},[\gamma]\rangle\pmod 2$ for every closed edge path $\gamma$. Hence $\ell=w_{1,S}$, and $Q$ is orientation-twisted.
\end{proof}

\subsection{Polarity Master Theorem}

\begin{lemma}[]
\label{lem:quotient-cross-ratio}
Let $T$ be a realized orientation-twisted polarity tiling of a quadrangulation $Q$ of a nonorientable surface $S$, and let $(\widetilde Q,\tau)$ be the orientation-cover lift of $Q$. Choose a two-coloring $c$ of $\widetilde Q^{(1)}$. Then $c\circ\tau=c+1$.

For $f\in F(Q)$, choose a lift $\widetilde f$ and write its positively oriented boundary, beginning at a color-zero vertex, as $\widetilde A-\widetilde C-\widetilde B-\widetilde D- \widetilde A$. Then $\chi_\beta(f):=\chi_\beta(A,B;C,D)$ is independent of the chosen lift and of the chosen color-zero initial vertex.
\end{lemma}

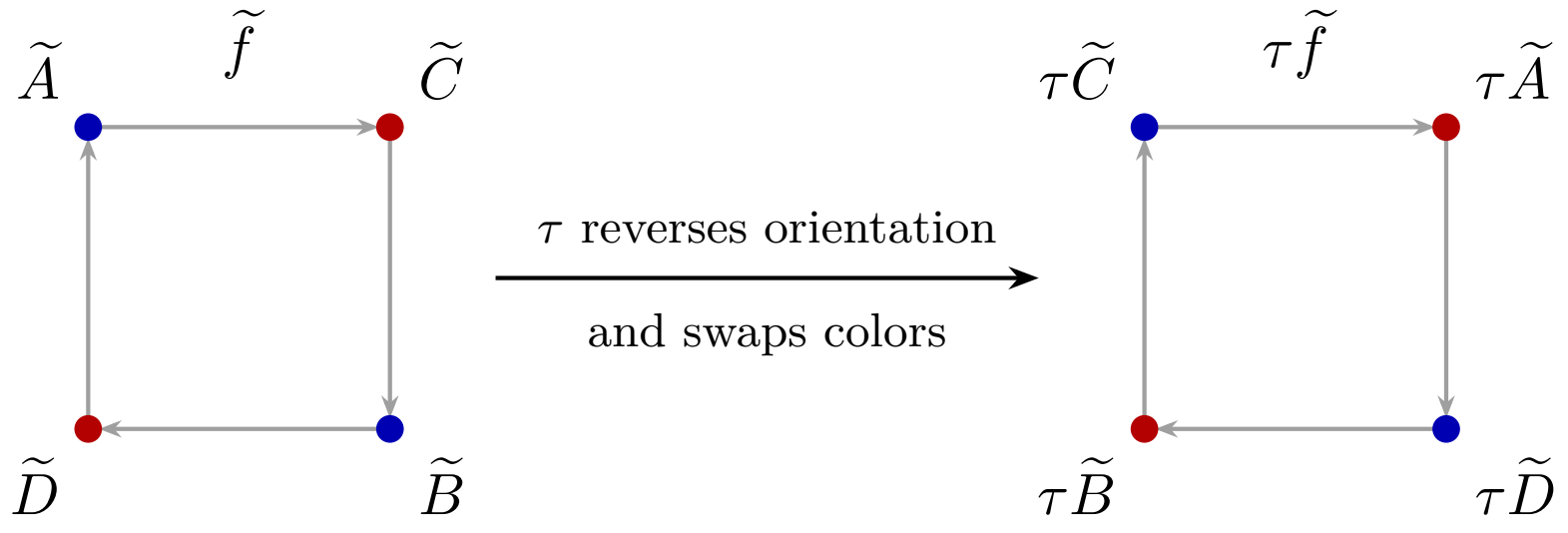
\begin{figure}[htbp]
\centering
\begin{tikzpicture}[
    font=\small,
    zero/.style={
        circle,fill=blue!70!black,draw=blue!70!black,
        inner sep=1.7pt
    },
    one/.style={
        circle,fill=red!70!black,draw=red!70!black,
        inner sep=1.7pt
    },
    boundary/.style={
        draw=gray!75,line width=.8pt,
        -{Stealth[length=1.5mm]},
        shorten <=2pt,shorten >=2pt
    }
]

\coordinate (A) at (-4.5,1);
\coordinate (C) at (-2.5,1);
\coordinate (B) at (-2.5,-1);
\coordinate (D) at (-4.5,-1);

\draw[boundary] (A)--(C);
\draw[boundary] (C)--(B);
\draw[boundary] (B)--(D);
\draw[boundary] (D)--(A);

\node[zero,label=above left:$\widetilde A$] at (A) {};
\node[one,label=above right:$\widetilde C$] at (C) {};
\node[zero,label=below right:$\widetilde B$] at (B) {};
\node[one,label=below left:$\widetilde D$] at (D) {};
\node at (-3.5,1.55) {$\widetilde f$};

\coordinate (tC) at (2.5,1);
\coordinate (tA) at (4.5,1);
\coordinate (tD) at (4.5,-1);
\coordinate (tB) at (2.5,-1);

\draw[boundary] (tC)--(tA);
\draw[boundary] (tA)--(tD);
\draw[boundary] (tD)--(tB);
\draw[boundary] (tB)--(tC);

\node[zero,label=above left:$\tau\widetilde C$] at (tC) {};
\node[one,label=above right:$\tau\widetilde A$] at (tA) {};
\node[zero,label=below right:$\tau\widetilde D$] at (tD) {};
\node[one,label=below left:$\tau\widetilde B$] at (tB) {};
\node at (3.5,1.55) {$\tau\widetilde f$};

\draw[-{Stealth[length=2mm]},line width=.8pt]
    (-1.8,0)--(1.8,0)
    node[midway,above=3pt,font=\scriptsize]
        {$\tau$ reverses orientation}
    node[midway,below=3pt,font=\scriptsize]
        {and swaps colors};
\end{tikzpicture}
\caption{The two lifts of a face and the two choices involved in defining its quotient cross ratio. The arrows show the positive boundary orientations.}
\label{fig:quotient-cross-ratio}
\end{figure}

Figure~\ref{fig:quotient-cross-ratio} illustrates the two lifted faces and their boundary orders.

\begin{proof}
The identity $c\circ\tau=c+1$ follows from Proposition~\ref{prop:cover-coloring}. Starting at the other color-zero vertex replaces $\chi_\beta(A,B;C,D)$ by $\chi_\beta(B,A;D,C)$, which has the same value.

The deck involution reverses the tautological orientation and exchanges the colors. Hence the positively oriented boundary of $\tau\widetilde f$, beginning at a color-zero vertex, is $\tau\widetilde C-\tau\widetilde A-\tau\widetilde D- \tau\widetilde B-\tau\widetilde C$. It gives $\chi_\beta(C,D;A,B)=\chi_\beta(A,B;C,D)$ by the symmetry of $\beta$.
\end{proof}

We define the face weights of an orientation-twisted polarity tiling with nonconjugate labels at adjacent vertices as follows. In the orientable case, fix an orientation of the surface and a two-coloring of the one-skeleton. For each face $f$, read $\chi_\beta(f)$ from its positively oriented boundary beginning at a color-zero vertex. In the nonorientable case, fix a two-coloring of the orientation-cover lift and use Lemma~\ref{lem:quotient-cross-ratio}.

\begin{theorem}[Polarity master theorem]
\label{thm:polarity-master}
Let $T$ be a realized orientation-twisted polarity tiling of a quadrangulation $Q$ of a closed connected surface $S$. If every face except a marked face $f_0$ is coherent, then $f_0$ is coherent.
\end{theorem}

\begin{proof}
In the orientable case, the equality $[\lambda_Q]=w_{1,S}=0$ implies that $Q^{(1)}$ is bipartite. The result is then the usual Fomin--Pylyavskyy cancellation argument.

In either case, the two faces incident with an edge contribute its pairing with opposite exponents, since their induced boundary directions are opposite. In the nonorientable case, compute both cross ratios using the lifts incident with the same lifted edge, as permitted by Lemma~\ref{lem:quotient-cross-ratio}. Thus all edge pairings cancel, giving
\begin{equation}\label{eq:polarity-face-product}
\prod_{f\in F(Q)}\chi_\beta(f)=1.
\end{equation}
If every face except $f_0$ is coherent, this identity forces $\chi_\beta(f_0)=1$, so $f_0$ is coherent.
\end{proof}

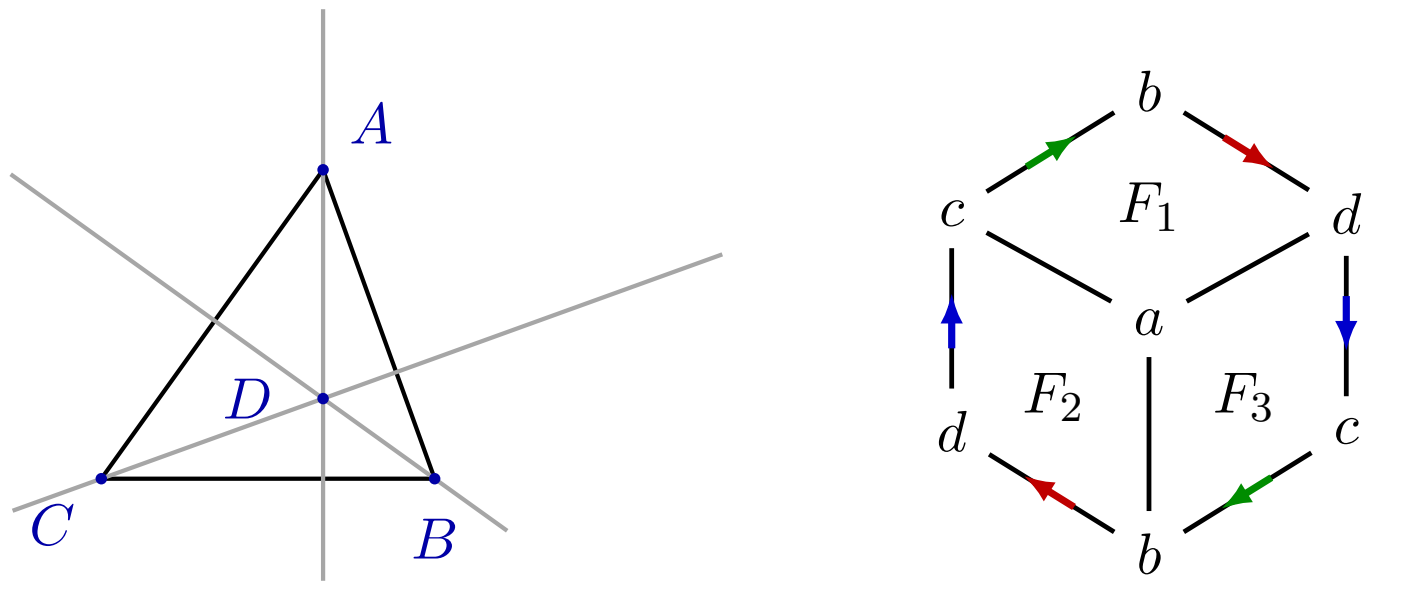
\begin{figure}[!htbp]
\centering
\begin{tikzpicture}[
    scale=.90,
    font=\small,
    geometric/.style={draw=gray!70,line width=.8pt},
    edge/.style={draw=black,line width=.9pt},
    vertex/.style={inner sep=4pt},
    paired/.style={
        -{Latex[length=2.1mm,width=1.5mm]},
        line width=1.35pt
    }
]

\begin{scope}[xshift=-3.8cm]
    \coordinate (A) at (.13,1.12);
    \coordinate (B) at (.95,-1.15);
    \coordinate (C) at (-1.50,-1.15);
    \coordinate (D) at (.13,{-1.15+1.63*.82/2.27});

    \draw[geometric,black] (A)--(B)--(C)--cycle;
    \draw[geometric] (.13,-1.90)--(.13,2.30);
    \draw[geometric]
        ($(B)!-.65!(D)$)--($(B)!3.80!(D)$);
    \draw[geometric]
        ($(C)!-.40!(D)$)--($(C)!2.80!(D)$);

    \foreach \P in {A,B,C,D}{
        \fill[blue!65!black] (\P) circle (1.2pt);
    }
    \node[blue!65!black,above right=2pt] at (A) {$A$};
    \node[blue!65!black,below=4pt]       at (B) {$B$};
    \node[blue!65!black,below left=2pt]  at (C) {$C$};
    \node[blue!65!black,left=6pt]        at (D) {$D$};
\end{scope}

\begin{scope}[xshift=2.4cm]
    \node[vertex] (bt) at (0,1.70) {$b$};
    \node[vertex] (cl) at (-1.45,.80) {$c$};
    \node[vertex] (dl) at (-1.45,-.80) {$d$};
    \node[vertex] (bb) at (0,-1.70) {$b$};
    \node[vertex] (cr) at (1.45,-.80) {$c$};
    \node[vertex] (dr) at (1.45,.80) {$d$};
    \node[vertex] (ac) at (0,0) {$a$};

    \foreach \u/\v in {bt/dr,dr/cr,cr/bb,bb/dl,dl/cl,cl/bt}{
        \draw[edge] (\u)--(\v);
    }
    \draw[edge] (cl)--(ac)--(dr);
    \draw[edge] (ac)--(bb);

    \draw[paired,green!55!black]
        ($(cl)!.38!(bt)$)--($(cl)!.64!(bt)$);
    \draw[paired,green!55!black]
        ($(cr)!.38!(bb)$)--($(cr)!.64!(bb)$);

    \draw[paired,red!75!black]
        ($(bt)!.38!(dr)$)--($(bt)!.64!(dr)$);
    \draw[paired,red!75!black]
        ($(bb)!.38!(dl)$)--($(bb)!.64!(dl)$);

    \draw[paired,blue!80!black]
        ($(dl)!.38!(cl)$)--($(dl)!.64!(cl)$);
    \draw[paired,blue!80!black]
        ($(dr)!.38!(cr)$)--($(dr)!.64!(cr)$);

    \node at (0,.85) {$F_1$};
    \node at (-.70,-.55) {$F_2$};
    \node at (.70,-.55) {$F_3$};
    \node at (0,-2.30) {Hesse tiling on $\mathbb{RP}^2$};
\end{scope}

\end{tikzpicture}
\caption{The Hesse configuration and its polarity tiling on $\mathbb{RP}^2$. Repeated boundary labels and matching colored arrows specify the edge identifications.}
\label{fig:hesse-geometry-tiling}
\end{figure}

\begin{example}[The Hesse polarity tiling]
\label{ex:hesse-polarity-tiling}
The Hesse configuration gives a polarity tiling of $\mathbb{RP}^2$. When the point labels at adjacent vertices are nonconjugate, Theorem~\ref{thm:polarity-master} shows that coherence of all but one face forces coherence of the remaining face. The geometric configuration and its projective-plane tiling are shown in Figure~\ref{fig:hesse-geometry-tiling}.
\end{example}

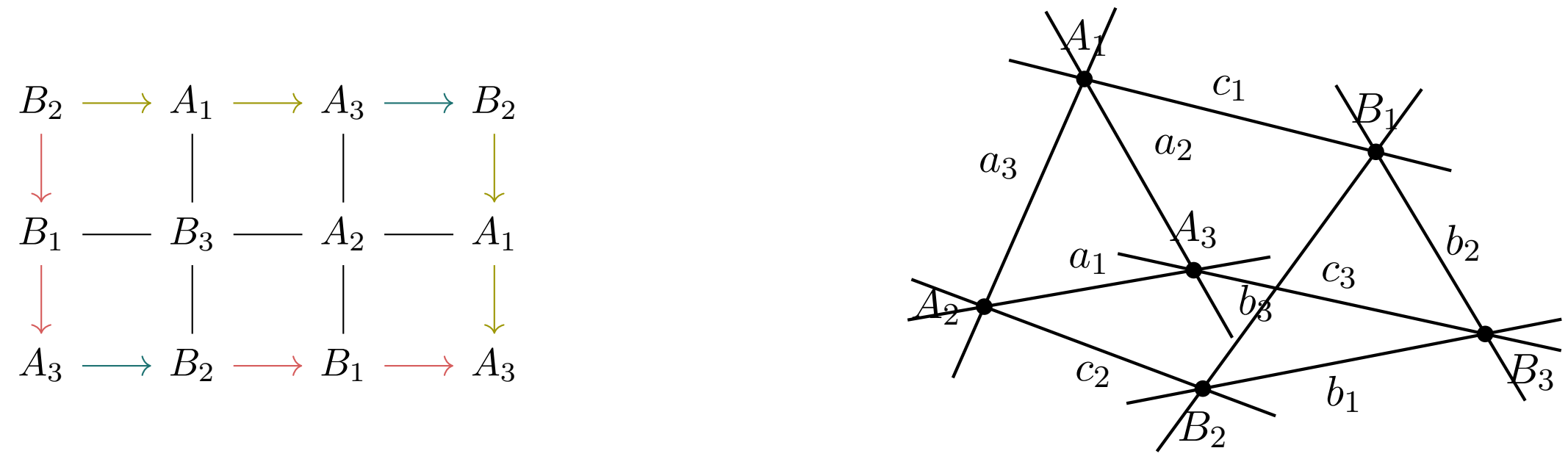
\begin{figure}[!htbp]
  \centering

  \begin{subfigure}[t]{0.46\textwidth}
    \vspace{0pt}
    \centering
    \begin{minipage}[c][5.2cm][c]{\linewidth}
    \centering
    \begin{tikzcd}[column sep=1.8em, row sep=1.8em]
      {B_2} & {A_1} & {A_3} & {B_2} \\
      {B_1} & {B_3} & {A_2} & {A_1} \\
      {A_3} & {B_2} & {B_1} & {A_3}
      \arrow[color={rgb,255:red,157;green,152;blue,7},
        from=1-1, to=1-2]
      \arrow[color={rgb,255:red,214;green,92;blue,92},
        from=1-1, to=2-1]
      \arrow[color={rgb,255:red,157;green,152;blue,7},
        from=1-2, to=1-3]
      \arrow[no head, from=1-2, to=2-2]
      \arrow[color={rgb,255:red,29;green,114;blue,114},
        from=1-3, to=1-4]
      \arrow[no head, from=1-3, to=2-3]
      \arrow[color={rgb,255:red,157;green,152;blue,7},
        from=1-4, to=2-4]
      \arrow[no head, from=2-1, to=2-2]
      \arrow[color={rgb,255:red,214;green,92;blue,92},
        from=2-1, to=3-1]
      \arrow[no head, from=2-2, to=2-3]
      \arrow[no head, from=2-2, to=3-2]
      \arrow[no head, from=2-3, to=2-4]
      \arrow[no head, from=2-3, to=3-3]
      \arrow[color={rgb,255:red,157;green,152;blue,7},
        from=2-4, to=3-4]
      \arrow[color={rgb,255:red,29;green,114;blue,114},
        from=3-1, to=3-2]
      \arrow[color={rgb,255:red,214;green,92;blue,92},
        from=3-2, to=3-3]
      \arrow[color={rgb,255:red,214;green,92;blue,92},
        from=3-3, to=3-4]
    \end{tikzcd}
    \end{minipage}
    \caption{The Klein bottle polarity tiling.}
    \label{fig:klein-bottle-polarity-tiling}
  \end{subfigure}
  \hfill
  \begin{subfigure}[t]{0.46\textwidth}
    \vspace{0pt}
    \centering
    \begin{minipage}[c][5.2cm][c]{\linewidth}
    \centering
    \begin{tikzpicture}[
      scale=0.82,
      every node/.style={font=\small},
      vertex/.style={circle,fill=black,inner sep=1.5pt},
      line/.style={
        thick,
        shorten <=-7mm,
        shorten >=-7mm
      }
    ]

      \coordinate (A1) at (-0.4,3.2);
      \coordinate (A2) at (-1.5,0.7);
      \coordinate (A3) at (0.8,1.1);

      \coordinate (B1) at (2.8,2.4);
      \coordinate (B2) at (0.9,-0.2);
      \coordinate (B3) at (4.0,0.4);

      \draw[line] (A2)--(A3)
        node[midway,above] {$a_1$};

      \draw[line] (A3)--(A1)
        node[midway,above right] {$a_2$};

      \draw[line] (A1)--(A2)
        node[midway,above left] {$a_3$};

      \draw[line] (B2)--(B3)
        node[midway,below] {$b_1$};

      \draw[line] (B3)--(B1)
        node[midway,right] {$b_2$};

      \draw[line] (B1)--(B2)
        node[midway,below left] {$b_3$};

      \draw[line] (A1)--(B1)
        node[midway,above] {$c_1$};

      \draw[line] (A2)--(B2)
        node[midway,below] {$c_2$};

      \draw[line] (A3)--(B3)
        node[midway,above] {$c_3$};

      \node[vertex,label=above:$A_1$] at (A1) {};
      \node[vertex,label=left:$A_2$] at (A2) {};
      \node[vertex,label=above:$A_3$] at (A3) {};

      \node[vertex,label=above:$B_1$] at (B1) {};
      \node[vertex,label=below:$B_2$] at (B2) {};
      \node[vertex,label=below right:$B_3$] at (B3) {};

    \end{tikzpicture}
    \end{minipage}
    \caption{The geometric realization.}
    \label{fig:klein-bottle-configuration}
  \end{subfigure}

  \caption{
  A Klein bottle polarity tiling and its geometric realization. Boundary edges with matching endpoint labels and orientations are identified.
  }
  \label{fig:klein-bottle-example}
\end{figure}

\begin{example}[A Klein bottle polarity tiling]
\label{ex:klein-bottle-tiling}
Let $A_1,A_2,A_3,B_1,B_2,B_3$ be points in the polarized projective plane. For $\{i,j,k\}=\{1,2,3\}$, set $a_i=A_jA_k$, $b_i=B_jB_k$, and $c_i=A_iB_i$. Consider the six conditions $a_i\perp_\beta c_i$ and $b_i\perp_\beta c_i$ for $i=1,2,3$. Under the nondegeneracy assumptions of Theorem~\ref{thm:polarity-master}, any five of these conditions imply the sixth.

Let $C_i=a_i\cap b_i$. Then $c_i\perp_\beta C_i$ for $i\in\{1,2,3\}$. If $C_1,C_2,C_3$ degenerate to a collinear triple, or $c_1,c_2,c_3$ degenerate to a concurrent triple, the same example specializes to the self-dual Desargues configuration. See Figure~\ref{fig:klein-bottle-example}.
\end{example}

To state the boundary version, we extend coherence to polygons with an even number of sides.

\begin{definition}[Coherent polarity polygon]
\label{def:coherent-polarity-polygon}
Let $P=A_1A_2\cdots A_{2n}$ be a polygon with cyclic point labels, where $n\ge1$, and assume that consecutive labels are nonconjugate. Choose nonzero vector representatives $a_i$ of $A_i$, set $a_{2n+1}=a_1$, and define
\[
  \chi_\beta(P)
  :=\frac{\prod_{i=1}^n\beta(a_{2i-1},a_{2i})}
          {\prod_{i=1}^n\beta(a_{2i},a_{2i+1})}.
\]
We call $P$ \emph{coherent} if $\chi_\beta(P)=1$.
\end{definition}

The ratio is independent of the vector representatives. A cyclic shift or reversal preserves it or replaces it by its inverse, so coherence is independent of these choices. For cyclic labels $A,C,B,D$, this recovers~\eqref{eq:fp-polarity-weight}.

\begin{remark}
The same cancellation applies to a compact connected surface with boundary. Each boundary component $B_j$ has an orientable collar, so $\langle w_{1,S},[B_j]\rangle=0$. Since $\ell=w_{1,S}$, its length is even; write it as $2n_j$.

In the orientable case, use the chosen orientation and coloring. In the nonorientable case, choose either lift of $B_j$, give it the boundary orientation induced from the tautological orientation of the cover, and begin at a color-zero vertex. The resulting cyclic labeling determines a polarity $2n_j$-gon $P_j=A_{j,1}C_{j,1}\cdots A_{j,n_j}C_{j,n_j}.$ The choice of lift does not affect $\chi_\beta(P_j)$, and cancellation of the interior edge pairings gives
\[
  \prod_{f\in F(Q)}\chi_\beta(f)
  =
  \prod_{j=1}^{b}\chi_\beta(P_j).
\]
In particular, if $S$ has one boundary component and every face is coherent, then its boundary is a coherent polarity $2n$-gon.

Putting $\ell_{j,i}=C_{j,i}^{\perp}$ identifies this polarity $2n_j$-gon with the alternately labeled polygon $A_{j,1}\ell_{j,1}\cdots A_{j,n_j}\ell_{j,n_j}$ of Fomin and Pylyavskyy. Under the hypotheses of their Proposition~4.2, its coherence is equivalent to tileability by coherent quadrilateral tiles, yielding the corresponding one-boundary statement of their Corollary~4.3 \cite[Proposition~4.2 and Corollary~4.3]{FominPylyavskyy2023}.
\end{remark}

\subsection{From FP Tilings to Polarity Tilings}
\label{subsec:polarity-tilings}
\label{subsec:polarity-conversion-duality}

We now explain how the polarity master theorem extends the FP master theorem, and use this connection to examine projective duality.

\begin{proposition}
\label{prop:polarity-conversion}
Replacing every realized hyperplane label of an FP tiling by its pole preserves face coherence. Conversely, let $T$ be a polarity tiling of a quadrangular complex $Q$, with labeling $\mathcal L:V(Q)\to\Lambda$ and realization $R$. A fixed two-coloring $c:V(Q)\to\Ftwo$ defines an FP tiling on $Q$ with labeling $(c,\mathcal L)$. For each $a\in\Lambda$, realize the point label $(0,a)$ by $R(a)$ and the hyperplane label $(1,a)$ by $R(a)^\perp$. For this fixed coloring, the two conversions are inverse at the level of projective realizations.
\end{proposition}

\begin{proof}
If $C=\ell^\perp$, then $P\notin\ell$ exactly when $\beta(p,c)\ne0$, where $p,c$ represent $P,C$. The equality $\chi_\beta(A,B;\ell^\perp,m^\perp)=\chi(A,B;\ell,m)$ shows that an FP tile is coherent exactly when its polarity tile is. At the level of projective realizations, the conversions are inverse because $(P^\perp)^\perp=P$ and $(\ell^\perp)^\perp=\ell$.
\end{proof}

Using this correspondence, we recover the FP master theorem as a special case of the polarity master theorem.

\begin{corollary}
\label{cor:fp-master}
Let $T$ be a realized FP tiling of a quadrangulation $Q$ of a closed connected oriented surface, and let $T'$ be the polarity tiling on $Q$ obtained by replacing each realized hyperplane label by its pole. Then $T'$ has nonconjugate point realizations at adjacent vertices, and $Q$ satisfies the orientation-twisted condition. For any face $f_0$, the hypothesis that all other faces are coherent and the conclusion that $f_0$ is coherent are each equivalent for $T$ and $T'$. Thus the FP master theorem for $T$ is an instance of Theorem~\ref{thm:polarity-master}.
\end{corollary}

\begin{proof}
The point and line vertices make $Q$ bipartite. Since the surface is oriented, $Q$ satisfies the orientation-twisted condition. Replacing line labels by poles gives nonzero edge pairings, because adjacent FP point and line labels are nonincident.

If all FP faces other than $f_0$ are coherent, they remain coherent after this replacement by Proposition~\ref{prop:polarity-conversion}. Theorem~\ref{thm:polarity-master} then shows that the remaining polarity face is coherent. By the same correspondence, the original FP face is coherent as well.
\end{proof}

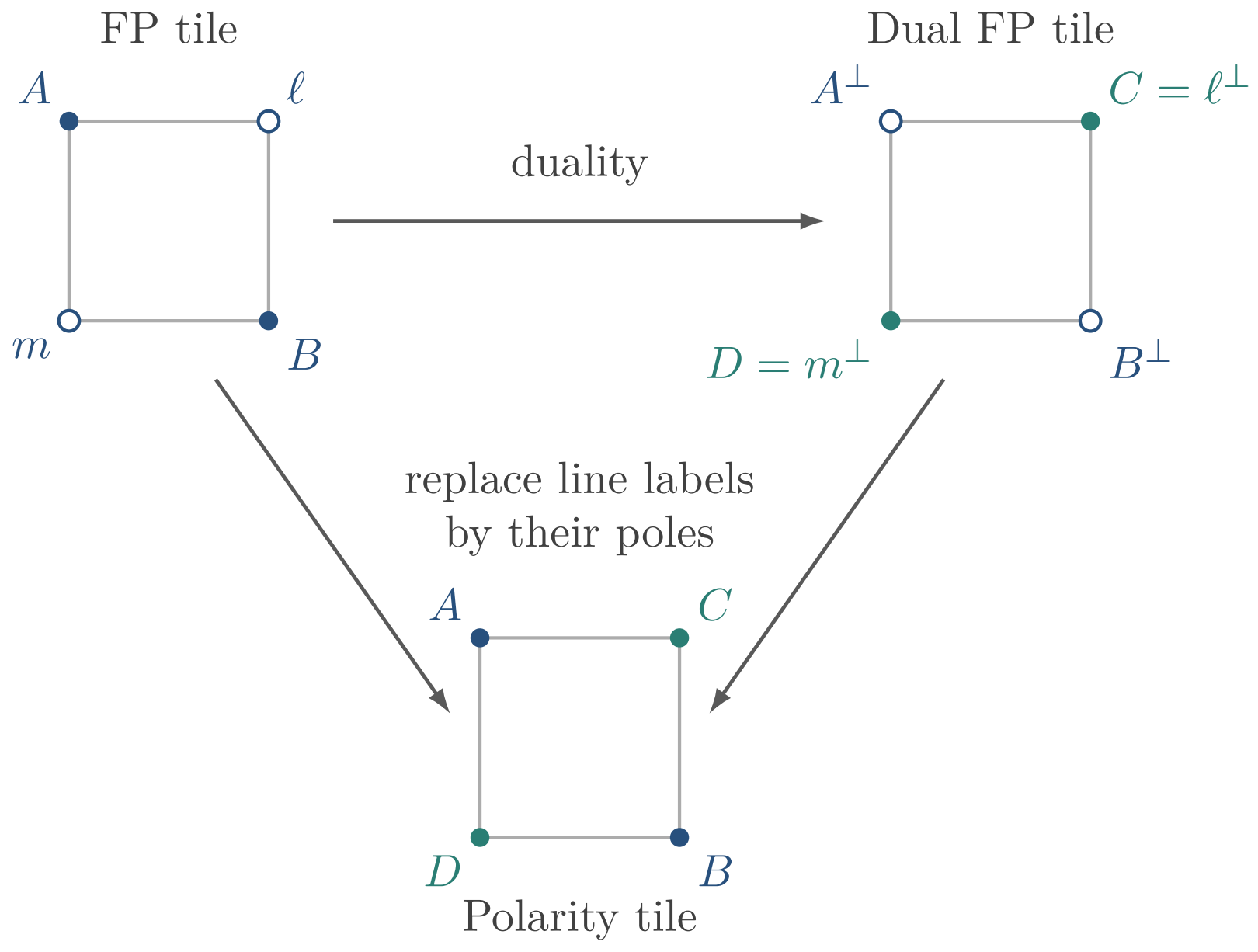
\begin{figure}[htbp]
\centering
\begingroup
\definecolor{dualityblue}{RGB}{40,80,125}
\definecolor{dualityteal}{RGB}{42,126,116}
\begin{tikzpicture}[
    font=\small,
    edge/.style={draw=gray!65,line width=.8pt},
    point/.style={circle,draw=dualityblue,fill=dualityblue,
        minimum size=4.2pt,inner sep=0pt},
    line/.style={circle,draw=dualityblue,fill=white,
        minimum size=5pt,inner sep=0pt,line width=.8pt},
    pole/.style={circle,draw=dualityteal,fill=dualityteal,
        minimum size=4.2pt,inner sep=0pt},
    lab/.style={inner sep=2pt},
    title/.style={text=black!75},
    arrow/.style={-{Latex[length=2.2mm,width=1.5mm]},
        draw=black!65,line width=.9pt}
]
    \begin{scope}[xshift=-3.5cm,yshift=2.5cm]
        \draw[edge] (-.85,.85)--(.85,.85)--(.85,-.85)--(-.85,-.85)--cycle;
        \node[point] at (-.85,.85) {};
        \node[line]  at (.85,.85) {};
        \node[point] at (.85,-.85) {};
        \node[line]  at (-.85,-.85) {};
        \node[lab,above left=3pt,text=dualityblue] at (-.85,.85) {$A$};
        \node[lab,above right=3pt,text=dualityblue] at (.85,.85) {$\ell$};
        \node[lab,below right=3pt,text=dualityblue] at (.85,-.85) {$B$};
        \node[lab,below left=3pt,text=dualityblue] at (-.85,-.85) {$m$};
        \node[title] at (0,1.65) {FP tile};
    \end{scope}

    \draw[arrow] (-2.1,2.5)--(2.1,2.5)
        node[midway,above=5pt,text=black!75] {duality};

    \begin{scope}[xshift=3.5cm,yshift=2.5cm]
        \draw[edge] (-.85,.85)--(.85,.85)--(.85,-.85)--(-.85,-.85)--cycle;
        \node[line] at (-.85,.85) {};
        \node[pole] at (.85,.85) {};
        \node[line] at (.85,-.85) {};
        \node[pole] at (-.85,-.85) {};
        \node[lab,above left=3pt,text=dualityblue] at (-.85,.85) {$A^\perp$};
        \node[lab,above right=3pt,text=dualityteal] at (.85,.85)
            {$C=\ell^\perp$};
        \node[lab,below right=3pt,text=dualityblue] at (.85,-.85) {$B^\perp$};
        \node[lab,below left=3pt,text=dualityteal] at (-.85,-.85)
            {$D=m^\perp$};
        \node[title] at (0,1.65) {Dual FP tile};
    \end{scope}

    \draw[arrow] (-3.1,1.15)--(-1.1,-1.7);
    \draw[arrow] (3.1,1.15)--(1.1,-1.7);
    \node[title,align=center] at (0,.05)
        {replace line labels\\by their poles};

    \begin{scope}[yshift=-1.9cm]
        \draw[edge] (-.85,.85)--(.85,.85)--(.85,-.85)--(-.85,-.85)--cycle;
        \node[point] at (-.85,.85) {};
        \node[pole]  at (.85,.85) {};
        \node[point] at (.85,-.85) {};
        \node[pole]  at (-.85,-.85) {};
        \node[lab,above left=3pt,text=dualityblue] at (-.85,.85) {$A$};
        \node[lab,above right=3pt,text=dualityteal] at (.85,.85) {$C$};
        \node[lab,below right=3pt,text=dualityblue] at (.85,-.85) {$B$};
        \node[lab,below left=3pt,text=dualityteal] at (-.85,-.85) {$D$};
        \node[title] at (0,-1.55) {Polarity tile};
    \end{scope}
\end{tikzpicture}
\endgroup
\caption{Replacing line labels by poles gives the same polarity tile.}
\label{fig:fp-duality-polarity}
\end{figure}

The same conversion sends an FP tiling and its projective dual to the same polarity tiling. Indeed, the cyclic labels $A,\ell,B,m$ and $A^\perp,\ell^\perp,B^\perp,m^\perp$ both become $A,\ell^\perp,B,m^\perp$ after replacing line labels by their poles; see Figure~\ref{fig:fp-duality-polarity}. By Proposition~\ref{prop:polarity-conversion}, a face is coherent if and only if its dual is.

The Pascal and Brianchon fillings below illustrate this relation. The conic polarity exchanges their point and line labels and turns the collinearity conclusion of Pascal's theorem into the concurrence conclusion of Brianchon's theorem.

\subsection{Pascal and Brianchon}

\begin{figure}[!htbp]
\centering
\begin{tikzpicture}[
    scale=2,
    font=\small,
    line width=.8pt,
    point/.style={circle,fill=black,inner sep=1.5pt},
    result/.style={circle,fill=blue!70!black,inner sep=1.7pt}
]
\draw[gray!65] (0,0) circle (1);
\node[gray!70!black] at (145:1.12) {$\mathcal C$};

\coordinate (P1) at (300:1);
\coordinate (P2) at (110:1);
\coordinate (P3) at (225:1);
\coordinate (P4) at (0:1);
\coordinate (P5) at (195:1);
\coordinate (P6) at (50:1);

\draw[red!70!black,name path=s12] (P1)--(P2);
\draw[red!70!black,name path=s45] (P4)--(P5);

\draw[green!45!black,name path=s23] (P2)--(P3);
\draw[green!45!black,name path=s56] (P5)--(P6);

\draw[orange!85!black,name path=s34] (P3)--(P4);
\draw[orange!85!black,name path=s61] (P6)--(P1);

\path[name intersections={of=s12 and s45,by=P7}];
\path[name intersections={of=s23 and s56,by=P8}];
\path[name intersections={of=s34 and s61,by=P9}];

\draw[blue!70!black,dashed,line width=1pt]
    ($(P8)!-.25!(P9)$)--($(P8)!1.25!(P9)$);

\node[point,label=below right:$P_1$] at (P1) {};
\node[point,label=above left:$P_2$] at (P2) {};
\node[point,label=below left:$P_3$] at (P3) {};
\node[point,label=right:$P_4$] at (P4) {};
\node[point,label=left:$P_5$] at (P5) {};
\node[point,label=above right:$P_6$] at (P6) {};

\node[result,label={
    [fill=white,inner sep=1pt]above right:$P_7$
}] at (P7) {};
\node[result,label={
    [fill=white,inner sep=1pt]above left:$P_8$
}] at (P8) {};
\node[result,label={
    [fill=white,inner sep=1pt]below right:$P_9$
}] at (P9) {};

\node[green!45!black,fill=white,inner sep=1pt,above left=2pt]
    at ($(P5)!.72!(P6)$) {$a$};
\node[green!45!black,fill=white,inner sep=1pt,below left=2pt]
    at ($(P2)!.72!(P3)$) {$c$};
\end{tikzpicture}
\caption{The Pascal configuration, with $P_8=a\cap c$. The three blue points $P_7,P_8,P_9$ are collinear.}
\label{fig:pascal}
\end{figure}
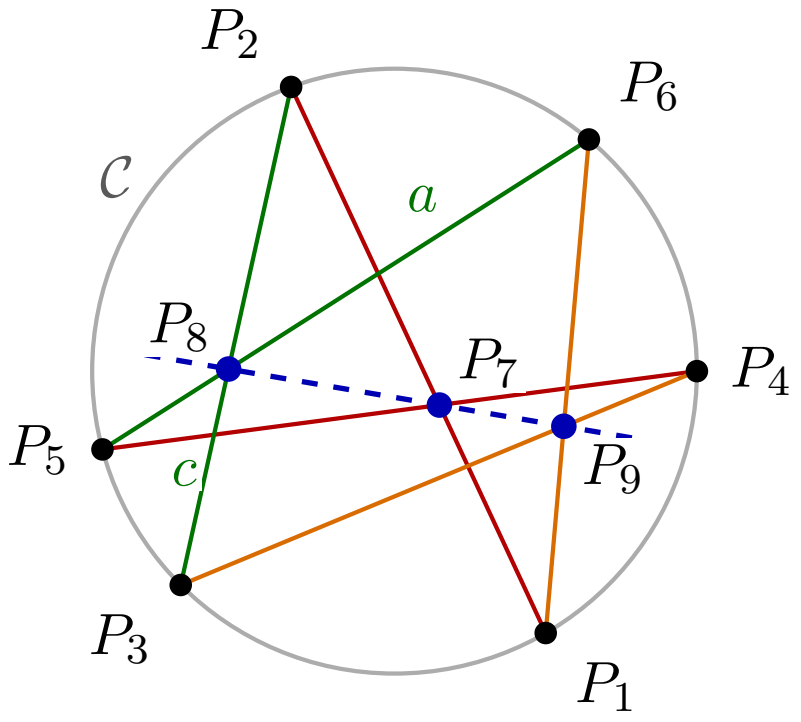

\begin{theorem}[Pascal's theorem]
\label{ex:pascal}

Let $P_1,\ldots,P_6$ be six distinct points on a nonsingular conic $\mathcal C$. Set $a=P_5P_6$ and $c=P_2P_3$, and define $P_7=P_1P_2\cap P_4P_5$, $P_8=a\cap c$, and $P_9=P_3P_4\cap P_6P_1$. Then the three points $P_7,P_8,P_9$ are collinear, as shown in Figure~\ref{fig:pascal}.

\end{theorem}

\begin{figure}[htbp]
\centering
\begin{tikzpicture}[
    scale=.86,
    font=\small,
    line width=.8pt,
    every node/.style={inner sep=4pt}
]
\foreach \shift/\filled in {-3/0,3/1}{
    \begin{scope}[xshift=\shift cm]
        \node (p1t) at (-.95,1.85) {$P_1$};
        \node (c)   at (.95,1.85) {$c$};
        \node (p4r) at (1.85,.95) {$P_4$};
        \node (b)   at (1.85,-.95) {$b$};
        \node (p1b) at (.95,-1.85) {$P_1$};
        \node (a)   at (-.95,-1.85) {$a$};
        \node (p4l) at (-1.85,-.95) {$P_4$};
        \node (d)   at (-1.85,.95) {$d$};

        \draw
            (p1t)--(c)--(p4r)--(b)--
            (p1b)--(a)--(p4l)--(d)--(p1t);

        \ifnum\filled=1
            \node (p7) at (-.58,.20) {$P_7$};
            \node (p9) at (.58,.20) {$P_9$};

            \draw (d)--(p7)--(c);
            \draw (p7)--(a);
            \draw (c)--(p9)--(b);
            \draw (p9)--(a);
        \fi
    \end{scope}
}
\node[font=\Large] at (0,0) {$+$};
\end{tikzpicture}
\caption{The boundary octagon and its quadrilateral filling.}
\label{fig:pascal-filling}
\end{figure}
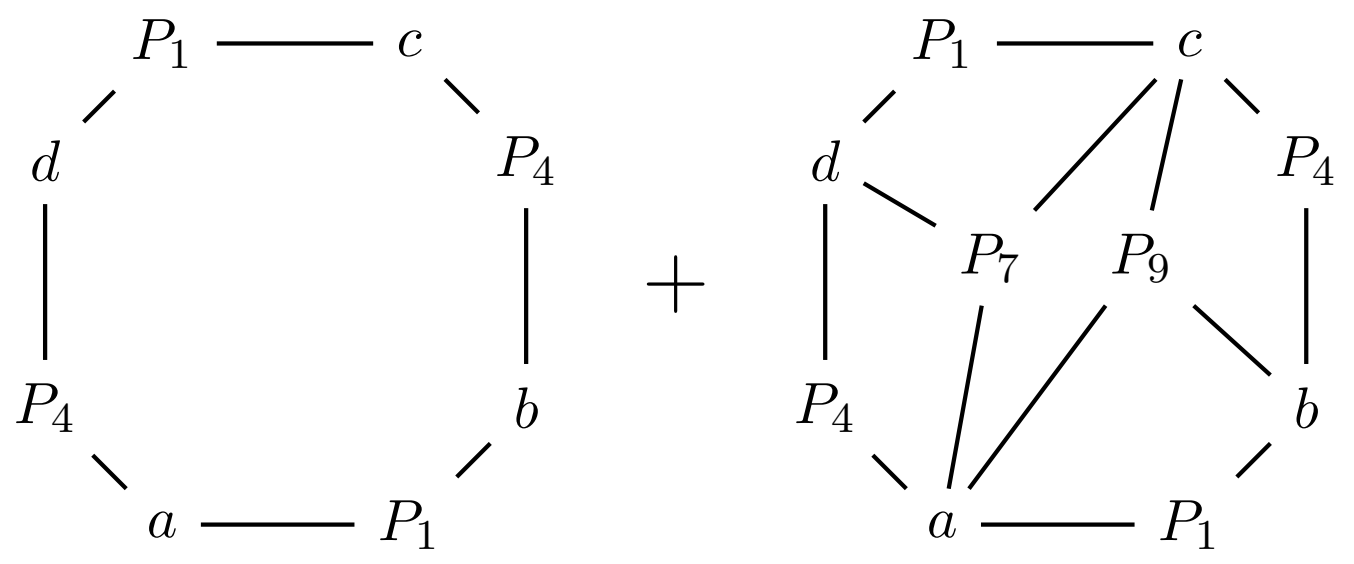

\begin{proof}
Set $b=P_3P_6$ and $d=P_2P_5$, and choose linear forms $a(x),b(x),c(x),d(x)$ defining the four lines. In Figure~\ref{fig:pascal-filling}, the definitions of $P_7$ and $P_9$ make the four outer faces coherent.

No three of $P_2,P_3,P_5,P_6$ are collinear. The quadratic forms vanishing at these four points are therefore spanned by $ac$ and $bd$. Thus $\mathcal C$ has an equation $\lambda a(x)c(x)+\mu b(x)d(x)=0$. Both coefficients are nonzero, since $\mathcal C$ is nonsingular.

Evaluating at $P_1$ and $P_4$ gives $ac/(bd)=-\mu/\lambda$ at both points. These evaluations are nonzero, since no line contains three of the six conic points. Hence
\[
 W=
 \frac{a(P_1)c(P_1)b(P_4)d(P_4)}
      {b(P_1)d(P_1)a(P_4)c(P_4)}
 =1.
\]

Thus the octagon with cyclic boundary labels $P_1,c,P_4,b,P_1,a,P_4,d$ is coherent. Since no three of the lines $a,b,c,d$ are concurrent, the dual form of \cite[Proposition~4.2]{FominPylyavskyy2023} ensures that this octagon admits a tiling by coherent quadrilateral tiles.

Glue this coherent filling to the five-face disk shown in Figure~\ref{fig:pascal-filling} along the corresponding boundary edges. This gives a quadrilateral tiling of the sphere in which every face except possibly the central face of the displayed disk is coherent. By Theorem~\ref{thm:fp-master}, that face is also coherent. Since $a\cap c=P_8$, this means that $P_7,P_8,P_9$ are collinear.
\end{proof}

Under the polarity defined by $\mathcal C$, each $P_i$ becomes the tangent line $\ell_i=P_i^\perp$, for $1\le i\le6$. Set $A=a^\perp$ and $C=c^\perp$. Then $P_8^\perp=(a\cap c)^\perp=AC=\ell_8$. Thus the three collinear points in Pascal's theorem become the three concurrent lines in Brianchon's theorem below.

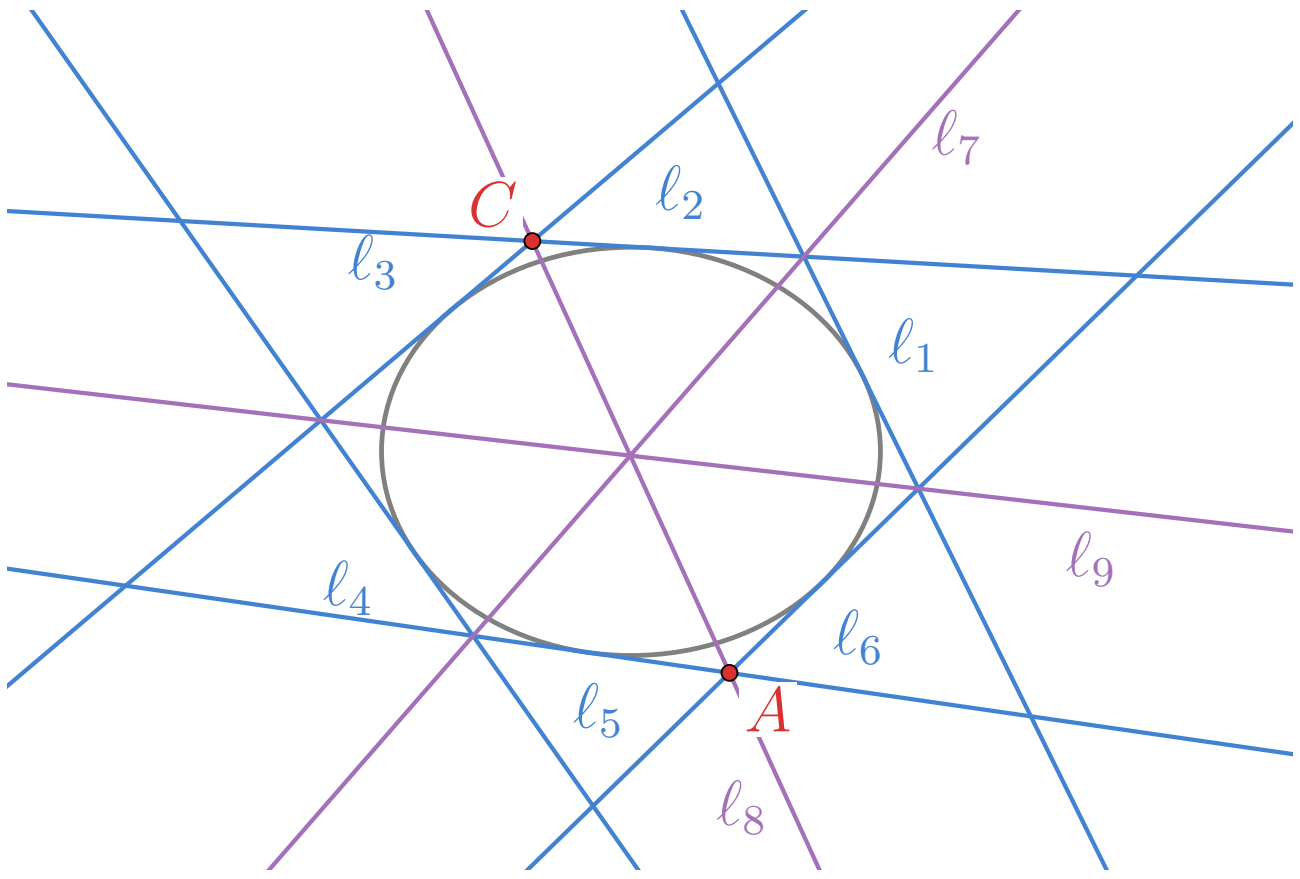
\begin{figure}[htbp]
\centering
\begin{tikzpicture}[
    x=1.65cm,
    y=1.35cm,
    font=\normalsize,
    line width=.8pt,
    tangent/.style={draw=brianchonblue},
    diagonal/.style={draw=brianchonpurple},
    linelabel/.style={inner sep=1pt},
    pointlabel/.style={
        fill=white,
        text=brianchonred,
        inner sep=1pt
    }
]

\foreach \i/\ang in {
    1/22,2/86,3/136,4/210,5/260,6/320
}{
    \coordinate (T\i) at (\ang:1);
    \coordinate (L\i) at
        ($(T\i)+({\ang+90}:5)$);
    \coordinate (R\i) at
        ($(T\i)+({\ang-90}:5)$);
}

\coordinate (V1) at
    (intersection of L1--R1 and L2--R2);
\coordinate (V2) at
    (intersection of L2--R2 and L3--R3);
\coordinate (V3) at
    (intersection of L3--R3 and L4--R4);
\coordinate (V4) at
    (intersection of L4--R4 and L5--R5);
\coordinate (V5) at
    (intersection of L5--R5 and L6--R6);
\coordinate (V6) at
    (intersection of L6--R6 and L1--R1);

\coordinate (A) at (V5);
\coordinate (C) at (V2);
\coordinate (O) at
    (intersection of V1--V4 and V2--V5);

\begin{scope}
    \clip (-2.50,-2.05) rectangle (2.65,2.16);

    \draw[gray,line width=.9pt]
        (0,0) circle (1);

    \foreach \i in {1,...,6}{
        \draw[tangent] (L\i)--(R\i);
    }

    \draw[diagonal]
        ($(V1)!-2!(V4)$)--($(V1)!3!(V4)$);

    \draw[diagonal]
        ($(V2)!-2!(V5)$)--($(V2)!3!(V5)$);

    \draw[diagonal]
        ($(V3)!-2!(V6)$)--($(V3)!3!(V6)$);
\end{scope}

\node[linelabel,text=brianchonblue,above right]
    at (1.00,.35) {$\ell_1$};
\node[linelabel,text=brianchonblue,above]
    at (.20,1.10) {$\ell_2$};
\node[linelabel,text=brianchonblue,above left]
    at (-.90,.76) {$\ell_3$};
\node[linelabel,text=brianchonblue,below left]
    at (-1.00,-.50) {$\ell_4$};
\node[linelabel,text=brianchonblue,below]
    at (-.13,-1.10) {$\ell_5$};
\node[linelabel,text=brianchonblue,below right]
    at (.78,-.74) {$\ell_6$};

\node[linelabel,text=brianchonpurple,right=3pt]
    at ($(O)!1.60!(V1)$) {$\ell_7$};
\node[linelabel,text=brianchonpurple,left=3pt]
    at ($(O)!1.62!(A)$) {$\ell_8$};
\node[linelabel,text=brianchonpurple,below=3pt]
    at ($(O)!1.60!(V6)$) {$\ell_9$};

\foreach \P in {A,C}{
    \filldraw[
        fill=brianchonred,
        draw=black,
        line width=.35pt
    ] (\P) circle (1.5pt);
}

\node[pointlabel,below right=2pt] at (A) {$A$};
\node[pointlabel,above left=2pt]  at (C) {$C$};

\end{tikzpicture}
\caption{The Brianchon configuration, with $\ell_8=AC$. The six blue lines are tangent to the conic, and the three purple lines $\ell_7,\ell_8,\ell_9$ are concurrent.}
\label{fig:brianchon}
\end{figure}

\begin{theorem}[Brianchon's theorem]
\label{ex:brianchon}

Let $k$ be a field with $\operatorname{char} k \ne 2$, and let $\ell_1,\ldots,\ell_6$ be six distinct tangents to a nonsingular conic $\mathcal C \subset \mathbb P^2(k)$. Set $A=\ell_5\cap\ell_6$ and $C=\ell_2\cap\ell_3$ and define $\ell_7=(\ell_1\cap\ell_2)(\ell_4\cap\ell_5)$, $\ell_8=AC$, and $\ell_9=(\ell_3\cap\ell_4)(\ell_6\cap\ell_1)$. Then the three lines $\ell_7,\ell_8,\ell_9$ are concurrent, as shown in Figure~\ref{fig:brianchon}.

\end{theorem}

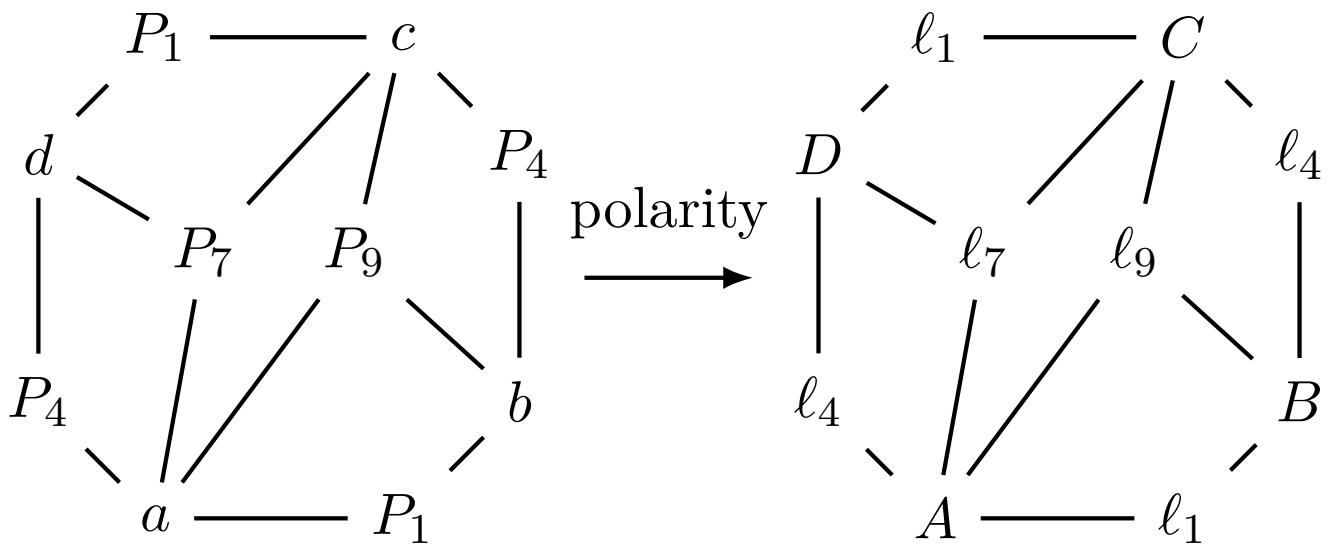
\begin{figure}[htbp]
\centering
\begin{tikzpicture}[
    scale=.86,
    font=\small,
    line width=.8pt,
    every node/.style={inner sep=4pt}
]
\begin{scope}[xshift=-3cm]
    \node (p1t) at (-.95,1.85) {$P_1$};
    \node (c)   at (.95,1.85) {$c$};
    \node (p4r) at (1.85,.95) {$P_4$};
    \node (b)   at (1.85,-.95) {$b$};
    \node (p1b) at (.95,-1.85) {$P_1$};
    \node (a)   at (-.95,-1.85) {$a$};
    \node (p4l) at (-1.85,-.95) {$P_4$};
    \node (d)   at (-1.85,.95) {$d$};

    \node (p7) at (-.58,.20) {$P_7$};
    \node (p9) at (.58,.20) {$P_9$};

    \draw
        (p1t)--(c)--(p4r)--(b)--
        (p1b)--(a)--(p4l)--(d)--(p1t);

    \draw (d)--(p7)--(c);
    \draw (p7)--(a);
    \draw (c)--(p9)--(b);
    \draw (p9)--(a);
\end{scope}

\draw[-{Latex[length=2mm]}] (-.65,0)--(.65,0)
    node[midway,above=3pt] {polarity};

\begin{scope}[xshift=3cm]
    \node (l1t) at (-.95,1.85) {$\ell_1$};
    \node (C)   at (.95,1.85) {$C$};
    \node (l4r) at (1.85,.95) {$\ell_4$};
    \node (B)   at (1.85,-.95) {$B$};
    \node (l1b) at (.95,-1.85) {$\ell_1$};
    \node (A)   at (-.95,-1.85) {$A$};
    \node (l4l) at (-1.85,-.95) {$\ell_4$};
    \node (D)   at (-1.85,.95) {$D$};

    \node (l7) at (-.58,.20) {$\ell_7$};
    \node (l9) at (.58,.20) {$\ell_9$};

    \draw
        (l1t)--(C)--(l4r)--(B)--
        (l1b)--(A)--(l4l)--(D)--(l1t);

    \draw (D)--(l7)--(C);
    \draw (l7)--(A);
    \draw (C)--(l9)--(B);
    \draw (l9)--(A);
\end{scope}
\end{tikzpicture}
\caption{The Pascal and Brianchon fillings related by polarity: $P_i^\perp=\ell_i$, and $a^\perp=A,\ b^\perp=B,\ c^\perp=C,\ d^\perp=D$.}
\label{fig:brianchon-filling}
\end{figure}

\begin{proof}
Apply the conic polarity to the Pascal filling, as shown in Figure~\ref{fig:brianchon-filling}. The points $P_7,P_9$ become the lines $\ell_7,\ell_9$, and $P_8^\perp=(a\cap c)^\perp=AC=\ell_8$. By Pascal's theorem, $P_7,P_8,P_9$ lie on a line $r$. Its pole $r^\perp$ lies on each of $\ell_7,\ell_8,\ell_9$, so these three lines are concurrent.
\end{proof}

\subsection{Two Special Cases of the Cayley--Bacharach Theorem}
\label{subsec:cayley-bacharach}

We give tiling proofs of two special cases of the classical Cayley--Bacharach theorem~\cite{Bacharach1886,EisenbudGreenHarris1996} for plane curves of degrees $3$ and $\ell$. The cubic is a union of three lines; the other curves are unions of lines, with one nonsingular conic allowed. The all-line case uses the FP master theorem and the line-conic configuration case uses our main theorem~\ref{thm:polarity-master}.

\begin{theorem}[The all-line configuration]
\label{prop:cb-lines}
Let $s_1,s_2,s_3\subset\mathbb P^2(k)$ be three nonconcurrent lines, and put
\[
 O_1=s_2\cap s_3,\qquad O_2=s_3\cap s_1,\qquad
 O_3=s_1\cap s_2.
\]
Let $n_1,\ldots,n_\ell$ and $m_1,\ldots,m_\ell$ be two families of lines, distinct within each family, none passing through an $O_i$. Suppose the $3\ell$ points $X_{ij}=s_i\cap n_j$ are distinct. If $\bigcup_a m_a$ contains all but one of these points, then it contains the remaining point.
\end{theorem}

\begin{theorem}[The line--conic configuration]
\label{ex:cb-line-conic}
Let $\operatorname{char}k\ne2$, and retain the $O_1,O_2,O_3$ and the $3\ell$ distinct points $X_{ij}=s_i\cap n_j$ of Theorem~\ref{prop:cb-lines}, with $\ell\ge2$. Let $\mathcal Q$ be a nonsingular conic avoiding the $O_i$, with each $s_i$ meeting $\mathcal Q$ in two distinct points. Let $m_1,\ldots,m_{\ell-2}$ be distinct lines avoiding the $O_i$. If
\[
 D=m_1\cup\cdots\cup m_{\ell-2}\cup\mathcal Q
\]
contains all but one of the $X_{ij}$, it contains the remaining point as well.
\end{theorem}

Throughout this subsection, square brackets denote labels in \emph{cyclic order}. Thus a polarity square $[A,C,B,D]$ is the tile $(A,B;C,D)$ of Definition~\ref{def:polarity-tile}. Distinct vertices may carry the same geometric label.

The proof of Theorem~\ref{prop:cb-lines} is adapted from Zhang's unpublished manuscript~\cite{ZhangTiledGenus}.

\begin{figure}[htbp]
\centering
\includegraphics[width=.72\textwidth]{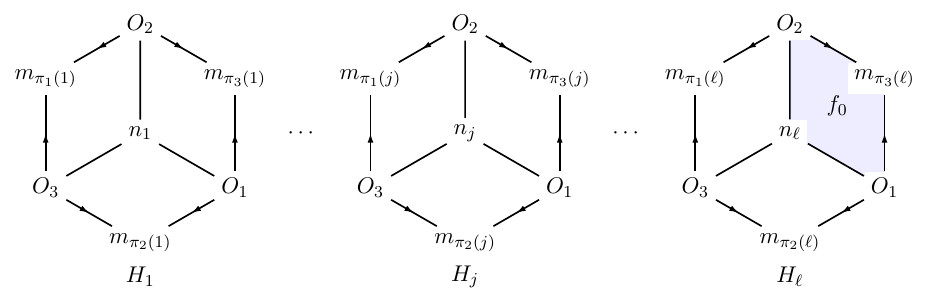}
\par\medskip
\includegraphics[width=.45\textwidth]{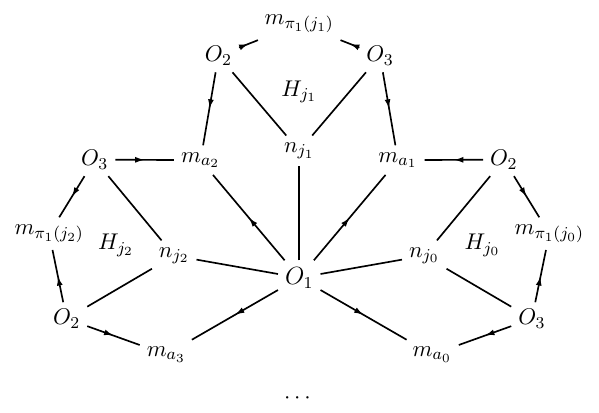}
\caption{The general all-line tiling from the array \eqref{eq:cb-line-array}. Top: one three-square patch $H_j$ for each column. Identify equally labelled boundary edges $O_i m_a$ according to their arrows, giving $3\ell$ faces. The shaded face is the starred condition $s_3\cap n_\ell\in m_{\pi_3(\ell)}$. Bottom: three consecutive patches in a longer cycle around $O_1$. Here $j_{r+1}=(\pi_2^{-1}\pi_3)(j_r)$ and $a_r=\pi_2(j_r)$, so $a_{r+1}=\pi_3(j_r)$. The shared edges $O_1m_{a_1}$ and $O_1m_{a_2}$ have been glued; the dots indicate continuation around the same vertex.}
\label{fig:cb-line-tiling}
\end{figure}

\begin{proof}[Proof of Theorem~\ref{prop:cb-lines}.]
Suppose the remaining point is $X_{3\ell}$. Each of the first two rows consists of $\ell$ distinct points covered by $\ell$ lines $m_a$. Since each $m_a$ meets $s_i$ once, this gives permutations $\pi_1,\pi_2$ with $X_{ij}\in m_{\pi_i(j)}$. In the third row, choose one covering line for each of the $\ell-1$ known points. These choices are distinct, since a line $m_a$ cannot meet $s_3$ at two different points. Assign the unused label to the last position, obtaining a permutation $\pi_3$. This last assignment is formal: its incidence condition is precisely what remains to be proved. The resulting configuration is the general $3\times\ell$ array
\begin{equation}
\label{eq:cb-line-array}
\begin{array}{c|ccccc}
 &n_1&\cdots&n_j&\cdots&n_\ell\\ \hline
 s_1&m_{\pi_1(1)}&\cdots&m_{\pi_1(j)}&\cdots&m_{\pi_1(\ell)}\\
 s_2&m_{\pi_2(1)}&\cdots&m_{\pi_2(j)}&\cdots&m_{\pi_2(\ell)}\\
 s_3&m_{\pi_3(1)}&\cdots&m_{\pi_3(j)}&\cdots&m_{\pi_3(\ell)}^{*}
\end{array}
\end{equation}
where $\pi_1,\pi_2,\pi_3\in S_\ell$. An entry $m_{\pi_i(j)}$ records $s_i\cap n_j\in m_{\pi_i(j)}$; the starred entry is the sole condition not yet assumed.

For each column $j$, form a hexagonal patch $H_j$ by joining a central vertex labelled $n_j$ to the three vertices labelled $O_1,O_2,O_3$. Its three oriented squares are
\begin{equation}
\label{eq:cb-line-squares}
[O_2,n_j,O_3,m_{\pi_1(j)}],\qquad
[O_3,n_j,O_1,m_{\pi_2(j)}],\qquad
[O_1,n_j,O_2,m_{\pi_3(j)}].
\end{equation}
The two occurrences of each $O_i n_j$ are already joined inside $H_j$. The remaining edges lie on the patch boundaries. For every $a\in\{1,\ldots,\ell\}$, their pairings are
\[
\begin{array}{c|c}
\text{boundary edge}&\text{patches containing its two copies}\\ \hline
O_1m_a&H_{\pi_2^{-1}(a)},\ H_{\pi_3^{-1}(a)}\\
O_2m_a&H_{\pi_3^{-1}(a)},\ H_{\pi_1^{-1}(a)}\\
O_3m_a&H_{\pi_1^{-1}(a)},\ H_{\pi_2^{-1}(a)}
\end{array}
\]
Glue each pair of boundary edges by matching their endpoints as indicated by the arrows in Figure~\ref{fig:cb-line-tiling}. The patches thus form a regular quadrangulation of a disjoint union of closed oriented surfaces.

All edge incidences are excluded by the assumptions. Every face except the one assigned to $X_{3\ell}$ is coherent. Apply Theorem~\ref{thm:fp-master} to its connected component. The remaining face is coherent, which gives $X_{3\ell}\in m_{\pi_3(\ell)}$.
\end{proof}

\begin{figure}[htbp]
\centering
\includegraphics[width=.32\textwidth]{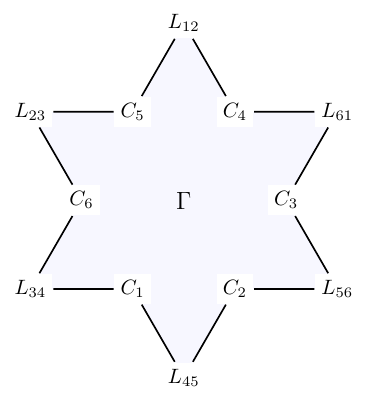}
\caption{The coherent twelve-sided polygon $\Gamma$ of Lemma~\ref{lem:cb-coherent-dodecagon}, with $L_{ij}=(C_iC_j)^\perp$. The interior is left unfilled.}
\label{fig:cb-coherent-dodecagon}
\end{figure}

We now assume $\operatorname{char}k\ne2$ and use the polarity of a nonsingular conic $\mathcal Q$. For a point $P\in\mathcal Q$, its polar $P^\perp$ is the tangent at $P$. We begin with the coherent polygon in Figure~\ref{fig:cb-coherent-dodecagon} determined by six points on $\mathcal Q$.

\begin{lemma}[A coherent conic dodecagon]
\label{lem:cb-coherent-dodecagon}
Let $C_1,\ldots,C_6$ be six distinct points on the nonsingular conic $\mathcal Q$, and put $L_{ij}=(C_iC_j)^\perp$. The polarity polygon
\begin{equation}
\label{eq:cb-coherent-dodecagon}
\Gamma=[C_1,L_{34},C_6,L_{23},C_5,L_{12},
        C_4,L_{61},C_3,L_{56},C_2,L_{45}]
\end{equation}
is coherent. Moreover, $\Gamma$ admits a bipartite disk filling by coherent polarity quadrilaterals.
\end{lemma}

\begin{proof}
Choose coordinates and homogeneous parameters $[u_i:v_i]$ on $\mathcal Q$ so that
\[
 c_i=(u_i^2,u_iv_i,v_i^2),\qquad
 [ij]=u_iv_j-v_iu_j.
\]
The parameters $[u_i:v_i]$ are distinct, so $[ij]\ne0$ whenever $i\ne j$. The chord $C_iC_j$ has the linear equation
\[
 q_{ij}(X,Y,Z)=v_iv_jX-(u_iv_j+v_iu_j)Y+u_iu_jZ=0,
\]
and evaluation at $c_k$ gives $q_{ij}(c_k)=[ki][kj]$. Choose a representative $l_{ij}$ of $L_{ij}$ such that $\beta(\,\cdot\,,l_{ij})=q_{ij}$. Each point on the boundary of $\Gamma$ is distinct from both endpoints of either adjacent chord. Hence all adjacent pairings are nonzero.

By Definition~\ref{def:coherent-polarity-polygon},
\[
\begin{aligned}
\chi_\beta(\Gamma)
 &=\frac{q_{34}(c_1)q_{23}(c_6)q_{12}(c_5)
          q_{61}(c_4)q_{56}(c_3)q_{45}(c_2)}
         {q_{34}(c_6)q_{23}(c_5)q_{12}(c_4)
          q_{61}(c_3)q_{56}(c_2)q_{45}(c_1)}\\[3pt]
 &=\frac{[13][62][51][46][35][24]}
         {[31][26][15][64][53][42]}
 =(-1)^6=1.
\end{aligned}
\]
The second equality cancels the common bracket factors. Thus $\Gamma$ is coherent.

For the filling assertion, replace each $L_{ij}$ by the line $L_{ij}^\perp=C_iC_j$. The resulting FP polygon has point labels $C_1,C_6,C_5,C_4,C_3,C_2$, no three of which are collinear. Its adjacent point and line labels are nonincident, as checked above. Therefore \cite[Proposition~4.2]{FominPylyavskyy2023} gives a filling by coherent quadrilaterals. Although stated over $\mathbb R$ and $\mathbb C$, its inductive construction uses only joins and meets under these same nonincidence conditions, so it works over $k$ as well. Replacing its line labels by their poles gives the required polarity filling, by Proposition~\ref{prop:polarity-conversion}.
\end{proof}

The following observation attaches a column-line label to a prescribed conic point.

\begin{lemma}[The conic attachment square]
\label{lem:cb-conic-collar}
Let $A,B\notin\mathcal Q$ be distinct points on a secant $s$, let $P\in s\cap\mathcal Q$, and let $n$ avoid $A,B$. Then the polarity square $[A,n^\perp,B,P]$ is coherent if and only if $P\in n$. All its adjacent pairs are nonconjugate.
\end{lemma}

\begin{proof}
The tangent $P^\perp$ meets $s$ only at $P$, so it contains neither $A$ nor $B$. If $n^\perp=P$, both assertions about coherence and incidence hold. Otherwise
\[
 (n^\perp P)^\perp=n\cap P^\perp.
\]
This point belongs to $AB=s$ exactly when it is $P$, which is equivalent to $P\in n$.
\end{proof}

Before proving Theorem~\ref{ex:cb-line-conic}, we describe its configuration table and the modification of the all-line tiling. Relabel the unknown point as $X_{3\ell}$, and label the six conic points so that
\[
 s_1\cap\mathcal Q=\{C_6,C_3\},\qquad
 s_2\cap\mathcal Q=\{C_5,C_2\},\qquad
 s_3\cap\mathcal Q=\{C_1,C_4\}.
\]
Replacing two line components by the conic gives the following row lists, whose entries are permuted into the columns $n_1,\ldots,n_\ell$:
\begin{equation}
\label{eq:cb-conic-array}
\begin{array}{c|c}
 &\text{entries in columns }n_1,\ldots,n_\ell\\ \hline
 s_1&\bigl(m_1,\ldots,m_{\ell-2},C_6,C_3\bigr)_{\pi_1}\\
 s_2&\bigl(m_1,\ldots,m_{\ell-2},C_5,C_2\bigr)_{\pi_2}\\
 s_3&\bigl(m_1,\ldots,m_{\ell-2},C_1,C_4\bigr)_{\pi_3}^{*}
\end{array}
\end{equation}
Here $(a_1,\ldots,a_\ell)_\pi$ means $(a_{\pi(1)},\ldots,a_{\pi(\ell)})$, and the star marks the last entry of the third reordered row. An entry $m_a$ in column $n_j$ records $X_{ij}\in m_a$; an entry $C_r$ records $C_r\in n_j$, equivalently $X_{ij}=C_r$. Thus every row contains each actual line once and its two conic points. The starred entry may be a line or a conic point, and records a condition still to be proved. The following lemma chooses the row permutations and the numbering within each conic pair.

\begin{lemma}[Two virtual line labels]
\label{lem:cb-virtual-lines}
Under the hypotheses of Theorem~\ref{ex:cb-line-conic}, the entry assignments in~\eqref{eq:cb-conic-array} can be chosen so that each column contains at most two conic-point entries, including the starred entry if it is assigned to a conic point. The six conic-point entries can then be replaced by two formal labels $\widetilde m_{\ell-1},\widetilde m_\ell$, each occurring exactly once in every row and at most once in every column.
\end{lemma}

\begin{proof}
Each row has $\ell$ component positions: one for each $m_a$ and two for $\mathcal Q$. In either complete row, its $\ell$ distinct grid points use all these positions, so all positions are distinct and both conic points are known grid points. In the third row, at least one known grid point lies on $\mathcal Q$, since the $\ell-2$ lines alone cannot cover its $\ell-1$ known points.

Assign every known point on $\mathcal Q$ to its conic position first. Assign each remaining known point to a covering line $m_a$; different points require different lines. Assign the single unused component position to the starred cell, without assuming its incidence. If both conic points in the third row are known, all six conic entries are actual grid points, the starred entry is a line entry, and every column has at most two conic entries.

Otherwise exactly five conic entries are known grid points, and the starred cell receives the sixth conic position. In this case all $3(\ell-2)$ line positions are known grid points and none lies on $\mathcal Q$. Only the starred column could have three formal conic entries. Suppose it does. Delete that column. The remaining $3(\ell-1)$ grid points consist of all the line positions and three points $A,B,C$ on $\mathcal Q$, one in each row. None of $A,B,C$ lies on an $m_a$. The line $t=AB$ avoids all three $O_i$, since $A$ and $B$ lie on different sides of the triangle, and $t$ is distinct from every $m_a$. The $\ell-1$ lines $m_1,\ldots,m_{\ell-2},t$ cover all the remaining grid points except possibly $C$. Theorem~\ref{prop:cb-lines}, applied to these $\ell-1$ columns, gives $C\in t$. This contradicts the fact that three distinct points of a nonsingular conic are not collinear. Thus the starred column also has at most two conic entries.

Make a bipartite graph with the three rows on one side and the columns on the other, and one edge for each conic entry. Every row vertex has degree two and every column vertex has degree at most two. Each nontrivial component is a path or an even cycle. Color its edges alternately with two colors, and use these colors as the labels $\widetilde m_{\ell-1}$ and $\widetilde m_\ell$. Each row receives each label once, and each column receives each label at most once.
\end{proof}

Number the points within each row so that the two color classes are $\{C_1,C_5,C_6\}$ and $\{C_2,C_3,C_4\}$. Replace their entries by
\[
 C_1,C_5,C_6\longmapsto\widetilde m_{\ell-1},\qquad
 C_2,C_3,C_4\longmapsto\widetilde m_\ell.
\]
With $\widetilde m_a=m_a$ for $a\le\ell-2$, the same row permutations give the virtual-line table
\begin{equation}
\label{eq:cb-virtual-array}
\begin{array}{c|ccccc}
 &n_1&\cdots&n_j&\cdots&n_\ell\\ \hline
 s_1&\widetilde m_{\pi_1(1)}&\cdots&\widetilde m_{\pi_1(j)}&\cdots&\widetilde m_{\pi_1(\ell)}\\
 s_2&\widetilde m_{\pi_2(1)}&\cdots&\widetilde m_{\pi_2(j)}&\cdots&\widetilde m_{\pi_2(\ell)}\\
 s_3&\widetilde m_{\pi_3(1)}&\cdots&\widetilde m_{\pi_3(j)}&\cdots&\widetilde m_{\pi_3(\ell)}^{*}
\end{array}
\end{equation}
Every row is now a permutation of $\ell$ distinct symbols, exactly as in the all-line construction. The two new symbols are used only to organize the tiling; no geometric realizations are assigned to the virtual lines.

Our proof modifies the all-line template as follows: retain only the boundaries of the two virtual three-square patches, removing their interiors, centers, and six spokes. Split the two hexagonal boundaries at their $O_1$-corners, and reconnect the split corners crosswise. This creates a twelve-sided boundary into which we insert the conic filling. We make this operation explicit before applying the master theorem.

The assignment of conic points to virtual labels is therefore
\[
\begin{array}{c|cc}
 &\widetilde m_{\ell-1}&\widetilde m_\ell\\ \hline
 s_1&C_6&C_3\\
 s_2&C_5&C_2\\
 s_3&C_1&C_4
\end{array}
\]
Let $\nu_i$ be the column line assigned to $C_i$, and put $N_i=\nu_i^\perp$. If the starred position is a conic position, the corresponding assignment does not yet assert $C_i\in\nu_i$.

Regroup the all-line squares around their line-label centers. Each virtual label gives a three-square hexagonal patch; write $T_-$ and $T_+$ for its formal center labels. The two retained boundary cycles are
\[
 B_-=[O_1,N_1,O_2,N_6,O_3,N_5],\qquad
 B_+=[O_1,N_4,O_2,N_3,O_3,N_2].
\]
Cut each cycle at its $O_1$-corner into an oriented six-edge path, and join the end of each path to the beginning of the other. The new boundary is
\begin{equation}
\label{eq:cb-conic-column-boundary}
 B=[O_1,N_1,O_2,N_6,O_3,N_5,
    O_1,N_4,O_2,N_3,O_3,N_2].
\end{equation}
Its two $O_1$-vertices remain distinct. This operation preserves all twelve boundary edges, their labels, and their directions; it is performed before the final edge gluings that define the closed surface.

\begin{figure}[htbp]
\centering
\includegraphics[width=.88\textwidth]{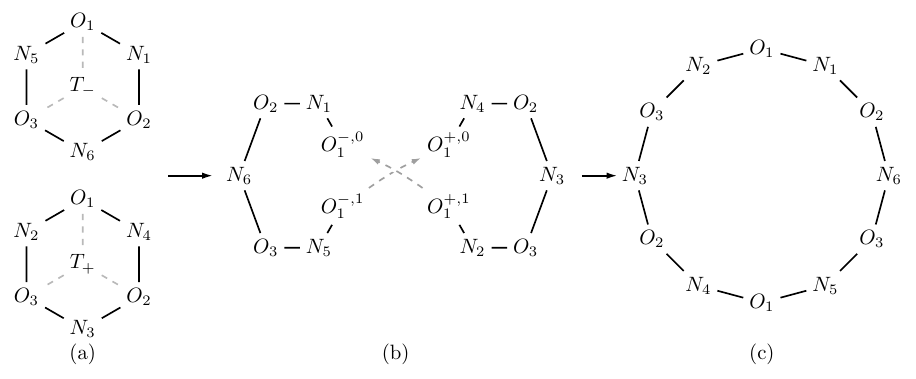}
\caption{Replacing the two virtual-line patches. The temporary centers $T_-,T_+$ correspond to $\widetilde m_{\ell-1},\widetilde m_\ell$. Delete the six dashed spokes and their centers, cut the two boundary cycles at their $O_1$-corners, and reconnect the ends crosswise. The resulting twelve-sided boundary is $B$ in~\eqref{eq:cb-conic-column-boundary}.}
\label{fig:cb-virtual-surgery}
\end{figure}

\begin{lemma}[The conic replacement disk]
\label{lem:cb-conic-disk}
The boundary $B$ admits a bipartite quadrilateral disk filling with six distinguished attachment faces. Each attachment face is coherent exactly when its assigned conic point lies on its column line. All other faces are coherent, and all adjacent labels are nonconjugate.
\end{lemma}

\begin{proof}
Along each boundary path $O_a,N_i,O_b$, attach the square $[O_a,N_i,O_b,C_i]$. By Lemma~\ref{lem:cb-conic-collar}, its coherence is precisely $C_i\in\nu_i$. The six attachment squares leave the inner boundary
\begin{equation}
\label{eq:cb-conic-boundary}
 [O_1,C_1,O_2,C_6,O_3,C_5,
   O_1,C_4,O_2,C_3,O_3,C_2].
\end{equation}

Put $L_{ij}=(C_iC_j)^\perp$, and let $\tau$ act on labels by $C_i\mapsto C_{i+2}$, $O_j\mapsto O_{j+1}$, with subscripts modulo six and three, respectively. Add the following four squares and their images under $\tau$ and $\tau^2$:
\begin{equation}
\label{eq:cb-conic-layers}
\begin{array}{c|ll}
 \text{ears}
 &[O_1,C_1,O_2,L_{12}]
 &[O_1,L_{12},O_3,C_2]\\[2pt]
 \text{expansion}
 &[O_1,L_{34},C_1,L_{45}]
 &[O_1,L_{45},C_2,L_{56}].
\end{array}
\end{equation}
The six ears leave the hexagonal boundary $[O_1,L_{34},O_3,L_{12},O_2,L_{56}]$. The six expansion squares leave exactly the polygon $\Gamma$ of~\eqref{eq:cb-coherent-dodecagon}. Fill $\Gamma$ by coherent polarity quadrilaterals using Lemma~\ref{lem:cb-coherent-dodecagon}.

The ears are coherent by Lemma~\ref{lem:cb-conic-collar}, using the chords $C_1C_2,C_3C_4,C_5C_6$. The first expansion square is coherent because $(L_{34}L_{45})^\perp=C_4\in O_1C_1;$ for the second, the corresponding intersection is $C_5\in O_1C_2$. Their cyclic images give the other four coherence conditions. All adjacent pairings in these layers are nonzero: the indicated chords avoid the adjacent $O_i$, and no chord contains a third one of the six conic points.

This construction is a disk with the prescribed boundary. Color the $O_i$ and the inner copies of the $C_i$ black, and the $N_i$, the $L_{ij}$, and the outer copies of the $C_i$ white. This agrees with the bipartite filling of $\Gamma$. Repeated labels on different vertices remain separate, and every closed quadrilateral is embedded.
\end{proof}

\begin{figure}[htbp]
\centering
\includegraphics[width=.47\textwidth]{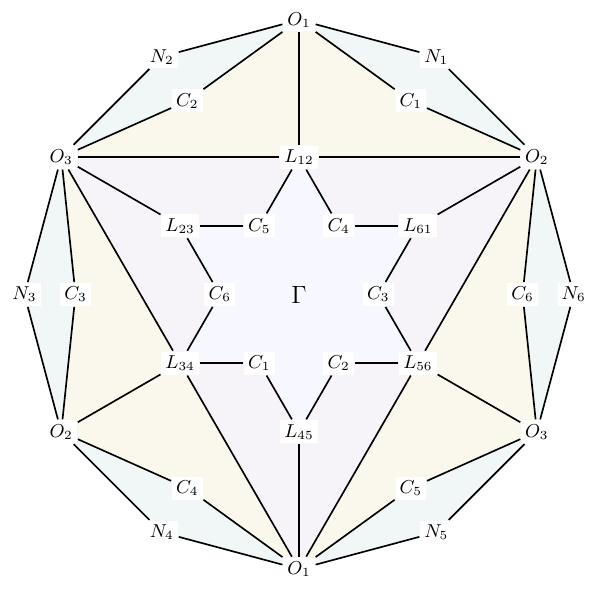}
\caption{The conic replacement disk. The three displayed layers contain six attachment squares, six ears, and six expansion squares. The central coherent polygon $\Gamma$ is filled by Lemma~\ref{lem:cb-coherent-dodecagon}. The six attachment faces record $C_i\in\nu_i$, with $N_i=\nu_i^\perp$; all other faces are coherent.}
\label{fig:cb-conic-filling}
\end{figure}

\begin{proof}[Proof of Theorem~\ref{ex:cb-line-conic}]
Choose the completed configuration table and its two virtual labels by Lemma~\ref{lem:cb-virtual-lines}. For each actual line $m_a$, use the three assigned columns to form its three-square hexagonal patch, replacing each realized line label by its pole. Replace the two virtual patches by the single conic disk just constructed.

The replacement preserves the boundary edge occurrences. Thus each edge $O_i n_j^\perp$ still occurs twice, with opposite boundary orientations. Pair the two occurrences with the same formal indices $(i,j)$, preserving their endpoint labels. As in the all-line proof, the vertex links are unions of circles; different link components give different vertices. We obtain a regular quadrangulation of a disjoint union of closed oriented surfaces, with a compatible bipartite coloring. The pairings do not identify two corners of any square: its $O$-indices are distinct, and its interior vertices do not participate in the boundary identifications. All adjacent labels are nonconjugate, so the hypotheses on the surface and its labels in Theorem~\ref{thm:polarity-master} are satisfied.

Every face except the one corresponding to the starred entry is coherent. If that entry is a line entry, the remaining face expresses the missing line incidence. If it is a conic entry, the remaining attachment face expresses $C_i\in\nu_i$, where $C_i$ is the unused conic point on the third row. This is precisely the assertion that $X_{3\ell}$ is that point of $\mathcal Q$. Apply Theorem~\ref{thm:polarity-master} to the connected component containing this face. Its coherence gives the required incidence in either case.
\end{proof}

\section{Cayley--Klein Geometry Interpretation of Coherence}
\label{sec:ck-coherence}

Projective geometry can describe incidences, but by itself supplies no ruler or notion of a right angle. A polarity adds precisely the quadratic data from which these measurements can be recovered. Let $V=\mathbb R^{n+1}$, $n\ge2$, and let $\beta$ be a nondegenerate symmetric bilinear form. Its zero locus $\mathcal A=\{[x]\in\mathbb P(V_{\mathbb C}):\beta(x,x)=0\}$ is called the \emph{absolute}; the same form defines the polarity $[p]\mapsto p^\perp$. Cayley--Klein geometry measures lengths and angles relative to this fixed quadric. We recall the construction from \cite[Chapter~13]{BadescuCarletti2024}, then interpret a coherent tile as a quadrilateral with orthogonal diagonal directions.

\begin{figure}[htbp]
\centering
\begin{tikzpicture}[scale=1.05,>=Stealth]
  \fill[blue!4] (0,0) circle (1.65);
  \draw[blue!65!black,thick] (0,0) circle (1.65);
  \coordinate (A) at (-.88,-.58);
  \coordinate (B) at (1.04,-.36);
  \coordinate (C) at (.15,1.08);
  \fill[blue!9] (A)--(B)--(C)--cycle;
  \draw[blue!70!black,thick] (A)--(B)--(C)--cycle;
  \foreach \P in {A,B,C} \fill (\P) circle (1.4pt);
  \node[below left] at (A) {$A$};
  \node[below right] at (B) {$B$};
  \node[above] at (C) {$C$};
  \node at (.05,.05) {$\mathbb H^2$};
  \draw[gray!75,->] (2.75,1.12)--(1.36,.93);
  \node[right,align=left] at (2.75,1.12)
    {absolute $\mathcal A$\\ideal boundary};
  \draw[gray!75,->] (2.75,-.8)--(.65,-.405);
  \node[right,align=left] at (2.75,-.8)
    {geodesic segment\\straight in the model};
\end{tikzpicture}
\caption{The Klein disk. Its interior is the hyperbolic plane, its boundary circle is the absolute, and geodesics are straight chords.}
\label{fig:ck-klein-disk}
\end{figure}
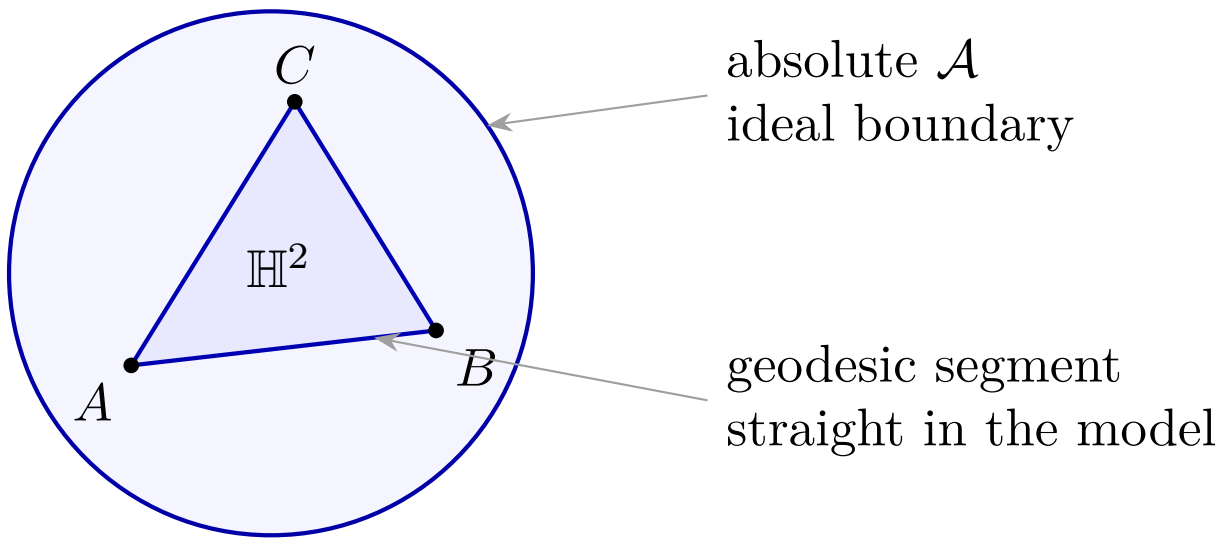

Intuitively, the absolute is the geometry's boundary at infinity. For example, take $\beta(x,y)=-x_0y_0+x_1y_1+x_2y_2$. In the affine chart $[1:x_1:x_2]$, the region $\beta(x,x)<0$ is the open unit disk. This is the \emph{Klein model} of the hyperbolic plane: geodesics are straight chords, and the boundary circle is the absolute (Figure~\ref{fig:ck-klein-disk}). The absolute provides the reference for measuring lengths and angles, as we explain next.

\subsection{Recovering Lengths and Angles}
\label{subsec:ck-metrics}

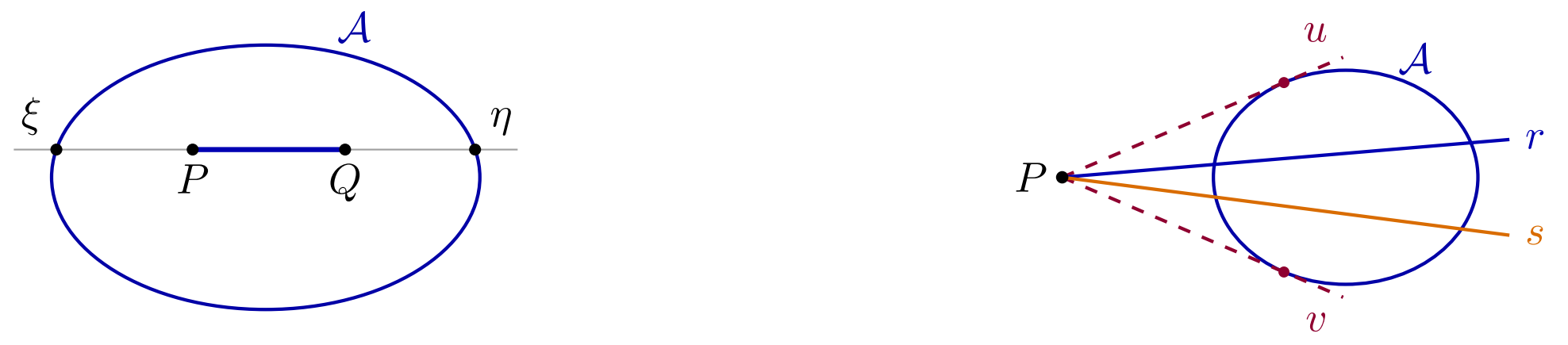
\begin{figure}[htbp]
\centering
\begin{minipage}[t]{.48\textwidth}
\centering
\begin{tikzpicture}[scale=1.05,>=Stealth]
  \path[use as bounding box] (-2.1,-1.3) rectangle (2.15,1.45);
  \draw[blue!65!black,thick] (0,0) ellipse (1.7 and 1.05);
  \node[blue!65!black,above] at (.7,.96) {$\mathcal A$};
  \draw[gray!75] (-2.0,.22)--(2.0,.22);
  \draw[blue!70!black,very thick] (-.58,.22)--(.63,.22);
  \foreach \x in {-1.662,-.58,.63,1.662}
    \fill (\x,.22) circle (1.35pt);
  \node[above left] at (-1.662,.22) {$\xi$};
  \node[above right] at (1.662,.22) {$\eta$};
  \node[below] at (-.58,.22) {$P$};
  \node[below] at (.63,.22) {$Q$};
\end{tikzpicture}

{\small Length: $\rho(P,Q;\xi,\eta)$}
\end{minipage}
\hfill
\begin{minipage}[t]{.48\textwidth}
\centering
\begin{tikzpicture}[scale=1.05,>=Stealth]
  \path[use as bounding box] (-2.1,-1.3) rectangle (2.15,1.45);
  \coordinate (P) at (-1.65,0);
  \coordinate (U) at (.11,.752);
  \coordinate (W) at (.11,-.752);
  \draw[blue!65!black,thick] (.6,0) ellipse (1.05 and .85);
  \node[blue!65!black,above] at (1.15,.71) {$\mathcal A$};
  \draw[purple!75!black,dashed,thick] (P)--(.58,.953)
    node[above left] {$u$};
  \draw[purple!75!black,dashed,thick] (P)--(.58,-.953)
    node[below left] {$v$};
  \draw[blue!70!black,thick] (P)--(1.9,.3) node[right] {$r$};
  \draw[orange!85!black,thick] (P)--(1.9,-.46) node[right] {$s$};
  \fill (P) circle (1.4pt);
  \fill[purple!75!black] (U) circle (1.25pt);
  \fill[purple!75!black] (W) circle (1.25pt);
  \node[left] at (P) {$P$};
\end{tikzpicture}

{\small Angle: $\rho(r,s;u,v)$}
\end{minipage}
\caption{The absolute supplies two reference points on a line and two reference lines in a pencil.}
\label{fig:ck-measurements}
\end{figure}

To measure the distance between two points $P,Q$, extend their joining line to meet the absolute at $\xi,\eta$. The four points have a projectively invariant cross ratio. Our convention, in an affine coordinate on the line, is $\rho(P,Q;\xi,\eta) =\frac{(p-\xi)(q-\eta)}{(p-\eta)(q-\xi)}$. The values at infinity are interpreted by limits. The absolute points serve as the two reference marks of a projective ruler; see the left-hand diagram in Figure~\ref{fig:ck-measurements}. Cross ratios multiply when intervals are concatenated, so $\rho(P,R;\xi,\eta)=\rho(P,Q;\xi,\eta)\rho(Q,R;\xi,\eta)$. Consequently, a continuous additive signed length invariant under the connected group preserving these reference points is locally of the form $\ell(P,Q)=k\operatorname{Log}\rho(P,Q;\xi,\eta)$. The logarithm turns this multiplication into addition. The constant $k$ sets the scale and makes the measurement real on the chosen domain; a logarithm branch, or a covering space, specifies the signed length. Ordinary distances also account for endpoint order and, in elliptic geometry, select the shortest arc.

Angles are recovered by the same construction in a pencil of lines. Let $r,s$ meet at an actual point $P$, and work in the projective plane they span. The two tangent lines $u,v$ from $P$ to the plane section of $\mathcal A$ provide its two reference directions. Taking the cross ratio in this pencil gives, with a suitable order of $u,v$,
\begin{equation}\label{eq:ck-angle}
  \angle(r,s)\equiv\frac{1}{2i}
    \operatorname{Log}\rho(r,s;u,v)\pmod\pi,
  \qquad
  r\perp s\ \Longleftrightarrow\ \rho(r,s;u,v)=-1.
\end{equation}
Thus a right angle is a harmonic relation with the two absolute directions. For an actual elliptic or hyperbolic point, these reference tangents are complex: they are the two isotropic lines through $P$. The right-hand diagram shows the tangency construction schematically, using a real conic to display its incidences.

\begin{example}[Hyperbolic geometry]
Take $\beta$ of signature $(n,1)$, with $n$ positive directions. The actual space is $X_{\mathrm{hyp}}=\{[p]:\beta(p,p)<0\}$, the higher-dimensional Klein ball. Choose representatives $p,q$ in the same component of the negative cone. At unit scale, $d_{\mathrm{hyp}}(P,Q) =\frac12\left|\log\rho(P,Q;\xi,\eta)\right| =\operatorname{arcosh} \frac{-\beta(p,q)}{\sqrt{\beta(p,p)\beta(q,q)}}$. Indeed, normalize $\beta(p,p)=\beta(q,q)=-1$ and write $q=\cosh t\,p+\sinh t\,e$, with $\beta(p,e)=0$, $\beta(e,e)=1$, and $t\ge0$. In the coordinate $z\mapsto[p+ze]$, the two points are $0,\tanh t$, and the absolute points are $1,-1$. Hence $\rho(0,\tanh t;1,-1)=\frac{1+\tanh t}{1-\tanh t}=e^{2t}$, so taking $k=1/2$ in the logarithmic length formula gives distance $t$.
\end{example}

\begin{example}[Elliptic geometry]
If $\beta$ is positive definite, the absolute has no real points, and the actual space is all of $\mathbb P(V)$. It is the unit sphere with antipodal points identified. At unit scale, $d_{\mathrm{ell}}(P,Q) =\frac12\left|\operatorname{Arg}\rho(P,Q;\xi,\eta)\right| =\arccos\frac{|\beta(p,q)|}{\sqrt{\beta(p,p)\beta(q,q)}}$, where $\operatorname{Arg}$ is the principal argument. Choose unit representatives with $\beta(p,q)\ge0$, and write $q=\cos t\,p+\sin t\,e$, where $e$ is a unit vector orthogonal to $p$ and $0\le t\le\pi/2$. In the same line coordinate, $P,Q$ are $0,\tan t$, while $\xi,\eta$ are $-i,i$. Thus $\rho(0,\tan t;-i,i)=\frac{1+i\tan t}{1-i\tan t}=e^{2it}$. Taking $k=1/(2i)$ gives length $t$, with the endpoint case obtained by taking a limit. The shortest-arc convention gives distances in $[0,\pi/2]$.
\end{example}

These constructions also identify the tangent metric. Use the unit sphere in the elliptic case and one sheet of the unit hyperboloid in the hyperbolic case. At a normalized point $p$, $T_pX=p^\perp, g_p(u,v)=\beta(u,v)$. This restriction is positive definite in both models. In an orthonormal tangent basis, the two isotropic directions have slopes $-i,i$; directions of slopes $0,\tan\theta$ have cross ratio $e^{2i\theta}$. Thus \eqref{eq:ck-angle} is precisely the angle of this Riemannian metric. Below, parallel transport is its Levi--Civita parallel transport.

\subsection{Coherence and Orthogonality}
\label{subsec:ck-meaning}

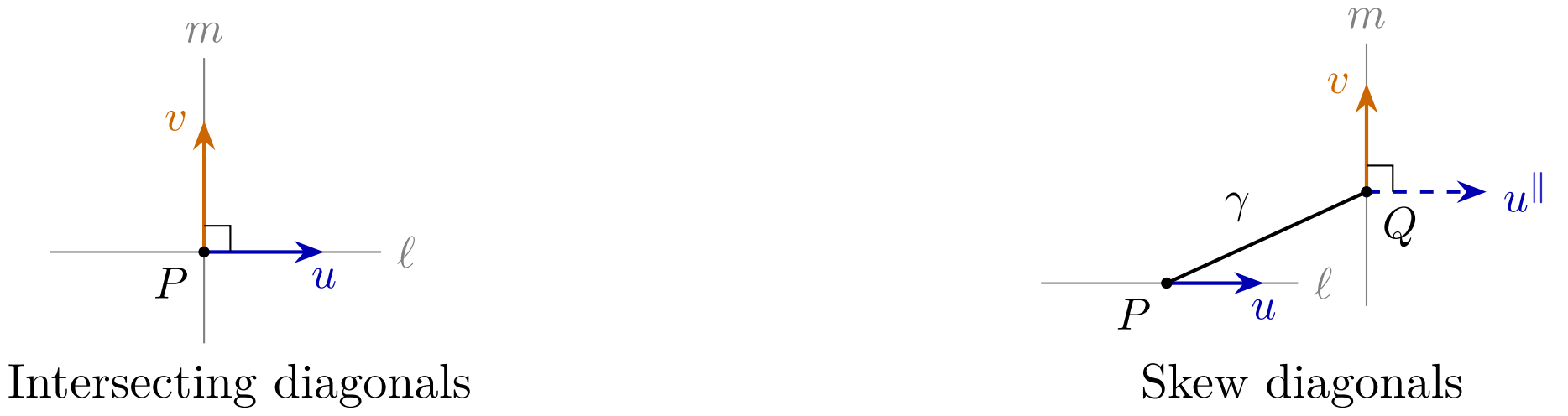
\begin{figure}[htbp]
\centering
\begin{minipage}[t]{.46\textwidth}
\centering
\begin{tikzpicture}[scale=.9,>=Stealth]
  \draw[gray] (-1.35,0)--(1.55,0) node[right] {$\ell$};
  \draw[gray] (0,-.8)--(0,1.7) node[above] {$m$};
  \draw[->,blue!70!black,thick] (0,0)--(1.05,0)
    node[below] {$u$};
  \draw[->,orange!80!black,thick] (0,0)--(0,1.15)
    node[left] {$v$};
  \draw (.23,0)--(.23,.23)--(0,.23);
  \fill (0,0) circle (1.4pt);
  \node[below left] at (0,0) {$P$};
\end{tikzpicture}

{\small Intersecting diagonals}
\end{minipage}
\hfill
\begin{minipage}[t]{.50\textwidth}
\centering
\begin{tikzpicture}[scale=.9,>=Stealth]
  \coordinate (P) at (0,0);
  \coordinate (Q) at (1.75,.8);
  \draw[gray] (-1.1,0)--(1.15,0) node[right] {$\ell$};
  \draw[gray] (1.75,-.2)--(1.75,2.1) node[above] {$m$};
  \draw[->,blue!70!black,thick] (P)--(.85,0)
    node[below] {$u$};
  \draw[->,orange!80!black,thick] (Q)--(1.75,1.75)
    node[left] {$v$};
  \draw[->,blue!70!black,thick,dashed] (Q)--(2.8,.8)
    node[right] {$u^\parallel$};
  \draw[thick] (P)--(Q) node[midway,above left] {$\gamma$};
  \draw (1.98,.8)--(1.98,1.03)--(1.75,1.03);
  \fill (P) circle (1.4pt);
  \fill (Q) circle (1.4pt);
  \node[below left] at (P) {$P$};
  \node[below right] at (Q) {$Q$};
\end{tikzpicture}

{\small Skew diagonals}
\end{minipage}
\caption{Coherence compares the two diagonal directions at the same point. For skew diagonals, transport $u$ along a shortest common perpendicular $\gamma$ before comparing it with $v$.}
\label{fig:ck-coherence}
\end{figure}

We can now read coherence geometrically. For intersecting diagonals, orthogonality is measured at their common point. For disjoint diagonals, first join their closest points by a shortest geodesic segment and transport one tangent direction along it. Comparing the directions at the same endpoint is the curved-space analogue of translating two Euclidean directions to a common origin; see Figure~\ref{fig:ck-coherence}. This description works in both elliptic and hyperbolic geometry.

\begin{theorem}
\label{thm:ck-orthogonality}
Let $\beta$ define an elliptic or hyperbolic space $X$, and let $(A,B;C,D)$ be a polarity tile with labels in $X$. Suppose $A\ne B$ and $C\ne D$, so that both diagonal lines are defined. Write $\ell,m$ for the geodesics with projective closures $AB,CD$. The following are equivalent:
\begin{enumerate}[label=\textup{(\roman*)}]
  \item The tile is polarity coherent: $AB\cap(CD)^\perp\ne\varnothing$.
  \item For homogeneous representatives $a,b,c,d$, $\beta_2(a\wedge b,c\wedge d) :=\beta(a,c)\beta(b,d)-\beta(a,d)\beta(b,c)=0$.
  \item There is a shortest segment joining $\ell$ to $m$ along which parallel transport makes their tangent directions orthogonal. When they intersect, the segment is constant and this says that the diagonals meet at a right angle.
\end{enumerate}
Here $\beta_2$ is the form induced on $\Lambda^2V$; \textup{(ii)} is the algebraic conjugacy condition for the two lines.
\end{theorem}

\begin{proof}
The restricted pairing $\langle a,b\rangle\times\langle c,d\rangle\to\mathbb R$ has a nonzero left kernel exactly when its determinant vanishes. This proves \textup{(i)}$\Longleftrightarrow$\textup{(ii)}.

Suppose first that the two geodesics have closest points. Choose normalized lifts $p,q$ of these points and unit tangent vectors $u,v$ to $\ell,m$. A shortest joining segment is perpendicular to both geodesics at its endpoints; when the segment is constant, $p=q$. In either case, $\beta(u,p)=\beta(u,q)=\beta(v,p)=\beta(v,q)=0$. The lifted segment lies in $\langle p,q\rangle$, so the constant ambient vector $u$ remains tangent and is parallel along it. If $\theta$ is its angle with $v$ after transport, then
\begin{equation}\label{eq:ck-skew-pairing}
  \beta_2(p\wedge u,q\wedge v)
  =\beta(p,q)\beta(u,v)=\beta(p,q)\cos\theta.
\end{equation}
The factor $\beta(p,q)$ is nonzero. Indeed, it cannot vanish for timelike vectors in the hyperbolic model. In the elliptic model, vanishing at a closest pair would make the two entire vector two-planes orthogonal, contradicting the nonconjugacy of adjacent labels in a polarity tile. Thus \eqref{eq:ck-skew-pairing} proves \textup{(ii)}$\Longleftrightarrow$\textup{(iii)}.

The only pairs without closest points are asymptotic hyperbolic geodesics. They share a null vector $k$ at infinity. Choosing timelike vectors $p,q$ on the respective vector two-planes gives $\beta_2(k\wedge p,k\wedge q) =-\beta(k,p)\beta(k,q)\ne0$. Thus both \textup{(ii)} and \textup{(iii)} fail in this remaining case, completing the proof.
\end{proof}

In a plane, condition \textup{(iii)} is simply perpendicular intersection. Indeed, disjoint geodesics with a common perpendicular have parallel tangent directions after transport along it. In higher dimensions the diagonals can be skew, and the transport description gives the corresponding meaning of orthogonality. In hyperbolic three-space, this condition is equivalent to a twist angle of $\pi/2$ between the diagonals; see \cite[Sections~2--3]{tabachnikov-skewers}. Consequently, Theorem~\ref{thm:polarity-master} gives an all-but-one statement for these orthogonal diagonal directions, with its original adjacency and topological hypotheses. The cases $A=B$ or $C=D$ remain coherent by convention, without assigning a direction to a collapsed diagonal.

\subsection{The Euclidean Limit}
\label{subsec:ck-euclidean}

In Euclidean geometry the ideal data move to the hyperplane at infinity. Besides $H_\infty=\{x_0=0\}$, one retains the direction quadric $Q_\infty:\ x_1^2+\cdots+x_n^2=0\quad\text{in }H_\infty$. This quadric has no real points, but its positive definite form measures angles between real directions. In the plane it consists of the two circular points $[0:1:i]$, $[0:1:-i]$. Together with the double hyperplane $x_0^2=0$, these are the Euclidean absolute of \cite[Sections~13.3.4 and~13.4]{BadescuCarletti2024}.

To see how the preceding coherence condition becomes Euclidean orthogonality, consider $\beta_\varepsilon((x_0,x),(y_0,y)) =x_0y_0+\varepsilon\,x\cdot y$. For $\varepsilon>0$ this gives an elliptic geometry; for $\varepsilon<0$ its negative gives a hyperbolic geometry with the same polarity. After lengths are rescaled by $|\varepsilon|^{-1/2}$, distances between fixed affine points converge to their Euclidean distances. On $H_\infty$, cancelling $\varepsilon$ in the absolute equation preserves $Q_\infty$. Direct substitution of $\varepsilon=0$ would lose this direction data.

For affine labels $a,c,b,d$, put $\hat a=(1,a)$, and similarly for the other points. The nondegenerate coherence equation in Theorem~\ref{thm:ck-orthogonality}\textup{(ii)} sets the pairing $\beta_{\varepsilon,2}(\hat a\wedge\hat b,\hat c\wedge\hat d)$ to zero. Its expansion has linear term $\varepsilon(a-b)\cdot(c-d)$ and quadratic term $\varepsilon^2\bigl((a\cdot c)(b\cdot d)-(a\cdot d)(b\cdot c)\bigr)$. Divide out the common factor $\varepsilon$ and take the limit. For $A\ne B$, $C\ne D$, the resulting condition is exactly $(a-b)\cdot(c-d)=0 \Longleftrightarrow AB\perp CD$. Here, also for skew lines, $\perp$ means orthogonality of their direction vectors. Equivalently, $AB\cap H_\infty$ and $CD\cap H_\infty$ are conjugate for $Q_\infty$. Thus orthogonal Euclidean diagonals are the first-order limit of polarity coherence; the rank-one form $\beta_0$ alone would make every face coherent.

The all-but-one conclusion passes to this limit by a short calculation. Label the vertices of an orientation-twisted quadrangulation by arbitrary points of $\mathbb E^n$, and use the face orders in the proof of Theorem~\ref{thm:polarity-master}. Writing $E_f=(a-b)\cdot(c-d)$ for a face in cyclic order $a,c,b,d$, we have $\chi_{\beta_\varepsilon}(f)=1+\varepsilon E_f+O(\varepsilon^2)$. Comparing the coefficients of $\varepsilon$ in the product identity $\prod_f\chi_{\beta_\varepsilon}(f)=1$ gives $\sum_f E_f=0$. The product identity applies for all sufficiently small nonzero $\varepsilon$, since every edge pairing tends to $1$. Hence orthogonality on all but one face forces it on the last; collapsed diagonals contribute zero automatically.

Apply this Euclidean all-but-one rule to the three faces of the Hesse polarity tiling in Example~\ref{ex:hesse-polarity-tiling} (Figure~\ref{fig:hesse-geometry-tiling}). Their diagonal pairs are $AB,CD$, $AC,BD$, and $AD,BC$, giving the following corollary.

\begin{corollary}
\label{cor:ck-tetrahedron}
If two pairs of opposite edges of a Euclidean tetrahedron are perpendicular, then so is the third pair. In symbols, $AB\perp CD,\quad AC\perp BD \Longrightarrow AD\perp BC$.
\end{corollary}

\begin{figure}[htbp]
\centering
\begin{tikzpicture}[scale=1.05]
  \coordinate (A) at (0,2.05);
  \coordinate (B) at (-1.6,-.75);
  \coordinate (C) at (1.35,-.55);
  \coordinate (D) at (1.65,1.15);
  \fill[blue!3] (A)--(B)--(C)--(D)--cycle;
  \draw[orange!85!black,thick,dashed] (B)--(D);
  \draw[blue!70!black,very thick] (A)--(B) (C)--(D);
  \draw[orange!85!black,very thick] (A)--(C);
  \draw[teal!75!black,very thick] (A)--(D) (B)--(C);
  \foreach \P in {A,B,C,D} \fill (\P) circle (1.6pt);
  \node[above] at (A) {$A$};
  \node[below left] at (B) {$B$};
  \node[below right] at (C) {$C$};
  \node[above right] at (D) {$D$};
\end{tikzpicture}

\caption{The three pairs of opposite edges. Each pair has its own color.}
\label{fig:ck-tetrahedron}
\end{figure}
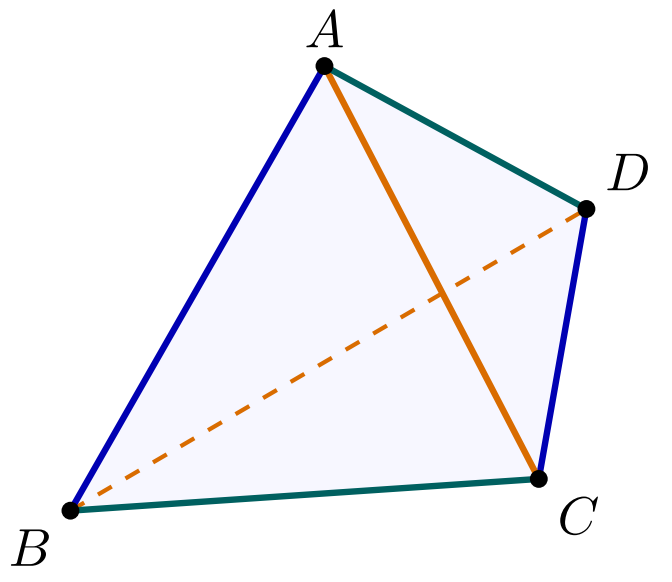

Here perpendicularity refers to the directions of the opposite edges, which are skew. Figure~\ref{fig:ck-tetrahedron} displays the three pairs. In the Euclidean plane, this corollary degenerates into the orthocenter theorem.

\section{Self-Polarity, Lifting, and Descent}
\label{sec:lifting-descent}

The orientation double cover used in the proof of the polarity master theorem also has an incidence-geometric interpretation. Its two-coloring restores point and line types, and the deck involution exchanges these types. A marked face on the original surface becomes a pair of faces on the covering surface, carrying equivalent incidence conditions that are exchanged by duality. The resulting tiling therefore carries a self-dual incidence statement, and the original polarity supplies a self-polar realization. We describe this lift and characterize when a self-polar realization descends to a tiling on a nonorientable surface.

\subsection{Lift to Self-dual Incidence Theorems}

For a tiling $T$ of a quadrangular complex $Q$, we write $\mathcal L$ for the formal vertex labeling. Incidence statements and maps of tilings refer to these formal data, without a chosen realization. In the FP case, a label includes its point or hyperplane type. In particular, the labels $(0,a)$ and $(1,a)$ in Proposition~\ref{prop:polarity-conversion} are distinct formal variables and may be realized independently. Self-duality refers to projective duality between $\mathbb P(V)$ and $\mathbb P(V^*)$, possibly after relabeling, and requires no chosen polarity.

Write $I(f)$ for the coherence condition of a face $f$. A tiling $T$ of $Q$ with a marked face $f$ carries the implication from all other face conditions to $I(f)$. With two distinct marked faces $f,g$, $T$ carries the equivalence of their conditions under the common premises:
\[
  \left(\bigwedge_{h\in F(Q)\setminus\{f,g\}} I(h)\right)
  \quad\Longrightarrow\quad
  \bigl(I(f)\Longleftrightarrow I(g)\bigr).
\]
The same convention applies to polarity tilings. On a closed oriented bipartite tiling, the ordinary master theorem gives both directions of this equivalence by taking either marked face as the conclusion.

\begin{definition}[Maps of tilings]
\label{def:tiling-map}
Let $T,T'$ be FP tilings of quadrangular complexes $Q,Q'$, with labelings $\mathcal L:V(Q)\to\Lambda$ and $\mathcal L':V(Q')\to\Lambda'$. A \emph{map of tilings} $(u,\sigma):T\to T'$ consists of a cellular map $u:Q\to Q'$ and a type-preserving bijection $\sigma:\Lambda\to\Lambda'$ such that $\mathcal L'\circ u=\sigma\circ\mathcal L$. It is an \emph{isomorphism of tilings} if $u$ is a cellular isomorphism.
\end{definition}

\begin{proposition}
\label{prop:color-swapping-duality}
Let $T$ be an FP tiling of a quadrangular complex $Q$ with labeling $\mathcal L$, and let $\tau$ be a cellular involution of $Q$. Suppose a bijection $\sigma$ of the label set exchanges point and hyperplane labels and satisfies $\mathcal L\circ\tau=\sigma\circ\mathcal L$. Then projective duality identifies $I(f)^\vee$ with $I(\tau f)$ for every face $f$, up to the relabeling $\sigma$. Consequently, if $f\ne\tau f$, the statement carried by $T$ with marked faces $f$ and $\tau f$ is self-dual.
\end{proposition}

\begin{proof}
The equality $\mathcal L(\tau v)=\sigma(\mathcal L(v))$ shows that $\tau$ exchanges point and hyperplane types, carrying the dual of each face relation to the relation on its image face, with labels renamed by $\sigma$. It permutes the faces outside $\{f,\tau f\}$ and exchanges the two marked faces. The common premises are therefore preserved, while the two sides of the marked equivalence are interchanged.
\end{proof}

A realization $R$ of an FP tiling $T$ of $Q$ is \emph{self-polar with respect to $\tau$} if $R(\tau v)=R(v)^\perp$ for every vertex $v$, where $\tau$ is a color-swapping cellular involution of $Q$ and $\perp$ is the polarity defined by $\beta$. This imposes an additional condition on the chosen realization.

\begin{proof}[Proof of Theorem~\ref{thm:canonical-lift}]
Let $\pi:\widetilde Q\to Q$ be the orientation-cover lift of the quadrangulation underlying $T$, with deck involution $\tau$. Write $\mathcal L:V(Q)\to\Lambda$ for the labeling of $T$ and $R$ for its realization. By Proposition~\ref{prop:cover-coloring}, $\widetilde Q$ has a two-coloring $c$ satisfying $c\circ\tau=c+1$. Define an FP labeling by
\[
  \widetilde{\mathcal L}(\widetilde v):=
  \bigl(c(\widetilde v),\mathcal L(\pi(\widetilde v))\bigr).
\]
This defines an FP tiling $\widetilde T$ of $\widetilde Q$. Mark its two faces above $f$. Since $\mathcal L\circ\pi$ is deck-invariant, $\tau$ acts on the typed labels by $(b,a)\mapsto(b+1,a)$ and exchanges the marked faces.

Let $\widetilde f$ be either lift of $f$. The lifted surface is closed and orientable, and its quadrangulation is bipartite. Under the common premises, the ordinary master theorem applied with conclusion $\tau\widetilde f$ gives $I(\widetilde f)\Rightarrow I(\tau\widetilde f)$; applying it with conclusion $\widetilde f$ gives the reverse implication. Proposition~\ref{prop:color-swapping-duality} shows that this equivalence is self-dual.

Define a realization of $\widetilde T$ by
\[
  \widetilde R(0,a)=R(a),\qquad
  \widetilde R(1,a)=R(a)^\perp
  \qquad(a\in\Lambda).
\]
Since $\tau$ preserves $a$ and exchanges the types, we have $\widetilde R(\tau\widetilde v)= \widetilde R(\widetilde v)^\perp$. Under polarity conversion, both $R(a)$ and $R(a)^\perp$ give the point $R(a)$. Proposition~\ref{prop:polarity-conversion} therefore identifies each lifted coherence condition for this realization with its counterpart downstairs.
\end{proof}

\begin{figure}[!htbp]
  \centering

  \begin{subfigure}[t]{0.37\textwidth}
    \vspace{0pt}
    \centering
    \begin{minipage}[c][4.1cm][c]{\linewidth}
    \centering
    \begin{tikzcd}[column sep=2.2em, row sep=1.6em]
      & {D_0} && {B_1} \\
      {C_1} && {A_1} \\
      & {A_0} && {C_0} \\
      {B_0} && {D_1}
      \arrow[no head, from=1-2, to=1-4]
      \arrow[no head, from=1-2, to=3-2]
      \arrow[no head, from=1-4, to=3-4]
      \arrow[no head, from=2-1, to=1-2]
      \arrow[dashed, no head, from=2-1, to=2-3]
      \arrow[no head, from=2-1, to=4-1]
      \arrow[dashed, no head, from=2-3, to=1-4]
      \arrow[dashed, no head, from=2-3, to=4-3]
      \arrow[no head, from=3-2, to=3-4]
      \arrow[no head, from=3-2, to=4-1]
      \arrow[no head, from=3-4, to=4-3]
      \arrow[no head, from=4-1, to=4-3]
    \end{tikzcd}
    \end{minipage}
    \caption{The Hesse hemicube.}
    \label{fig:hesse-hemicube}
  \end{subfigure}%
  \hspace{0.015\textwidth}%
  \begin{minipage}[t]{0.12\textwidth}
    \vspace{0pt}
    \begin{minipage}[c][4.1cm][c]{\linewidth}
      \centering
      \begin{tikzpicture}
        \draw[-{Latex[length=2.4mm,width=1.5mm]},line width=0.65pt]
          (0,0) -- (1.25,0) node[midway,above=2pt] {\small lift};
      \end{tikzpicture}
    \end{minipage}
  \end{minipage}%
  \hspace{0.015\textwidth}%
  \begin{subfigure}[t]{0.37\textwidth}
    \vspace{0pt}
    \centering
    \begin{minipage}[c][4.1cm][c]{\linewidth}
    \centering
    \begin{tikzcd}[column sep=2.2em, row sep=1.6em]
      & {A_2} && {a_2} \\
      {b} && {C} \\
      & {c} && {B} \\
      {A_1} && {a_1}
      \arrow[no head, from=1-2, to=1-4]
      \arrow[no head, from=1-2, to=3-2]
      \arrow[no head, from=1-4, to=3-4]
      \arrow[no head, from=2-1, to=1-2]
      \arrow[no head, from=2-1, to=2-3]
      \arrow[no head, from=2-1, to=4-1]
      \arrow[no head, from=2-3, to=1-4]
      \arrow[no head, from=2-3, to=4-3]
      \arrow[no head, from=3-2, to=3-4]
      \arrow[no head, from=3-2, to=4-1]
      \arrow[no head, from=3-4, to=4-3]
      \arrow[no head, from=4-1, to=4-3]
    \end{tikzcd}
    \end{minipage}
    \caption{The Desargues cube.}
    \label{fig:desargues-cube}
  \end{subfigure}

  \caption{The orientation-cover lift from the Hesse hemicube to the Desargues cube. The left-hand drawing represents the hemicube by identifying antipodal cells of the cube.}
  \label{fig:hesse-desargues-lift}
\end{figure}
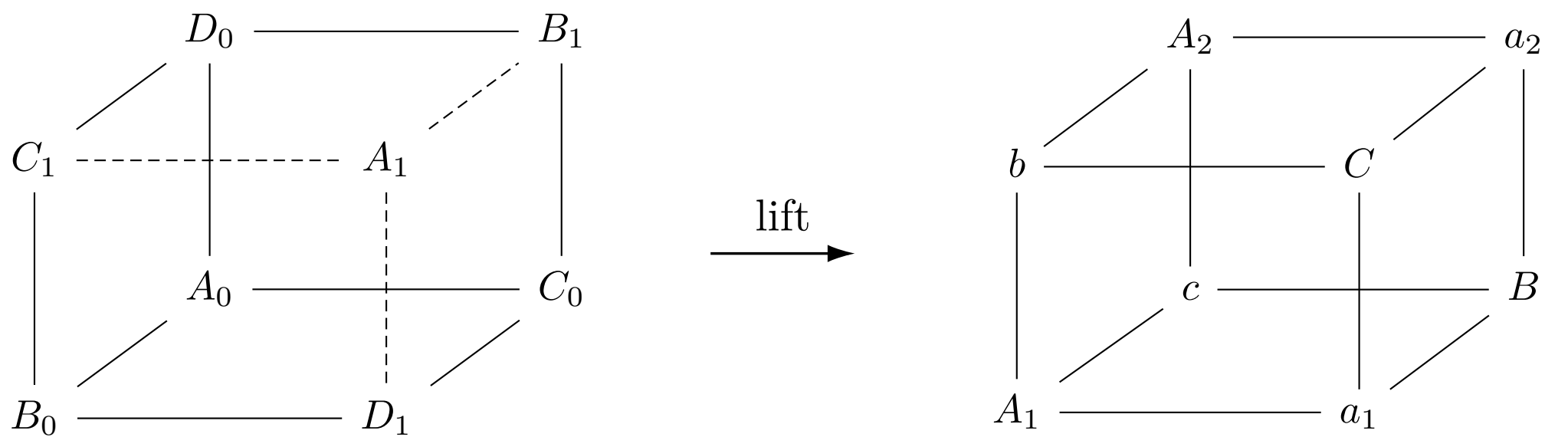

\begin{example}[The Hesse--Desargues lift]
\label{ex:hesse-desargues-lift}
The hemicube quadrangulation of $\mathbb{RP}^2$ lifts to the cube quadrangulation of $S^2$. At the level of incidence statements, this realizes the classical passage from Hesse's theorem to Desargues' theorem \cite{Brown1960}.

If $f$ is the conclusion face of the hemicube, its two lifts $\widetilde f$ and $\tau\widetilde f$ are antipodal cube faces. The canonical lift is therefore an FP tiling $\widetilde T$ of the cube quadrangulation $\widetilde Q$, with marked faces $\widetilde f$ and $\tau\widetilde f$. Under the four remaining face conditions, the conditions on these two faces are equivalent and exchanged by projective duality. The original Hesse polarity gives a self-polar realization of this Desargues tiling, with polar labels at antipodal vertices. Thus the one-boundary Hesse tiling lifts to the two-boundary Desargues tiling; see Figure~\ref{fig:hesse-desargues-lift}.
\end{example}

\FloatBarrier

\subsection{Descent of Self-polar Realizations}

For descent, the polarity must be compatible with a cellular involution of the chosen tiling. A local condition on edges ensures both freeness of the action and regularity of the quotient.

\begin{definition}[Admissible quadrangulation involution]
\label{def:admissible}
Let $\widetilde Q$ be a connected bipartite regular quadrangulation of a closed orientable surface $\widetilde S$, with two-coloring $c$. A cellular involution $\tau$ is \emph{admissible} if
\begin{enumerate}[label=\textup{(\roman*)}]
  \item it reverses the orientation of $\widetilde S$;
  \item it exchanges the colors;
  \item no vertex $v$ is adjacent to $\tau v$.
\end{enumerate}
\end{definition}

\begin{lemma}
\label{lem:quotient-regularity}
A color-swapping cellular involution $\tau$ of a bipartite regular quadrangulation of a closed surface acts freely with a regular quadrangulation quotient if and only if no vertex is adjacent to its image under $\tau$.
\end{lemma}

\begin{proof}
Suppose no vertex is adjacent to its image. If $\overline f\cap\tau(\overline f)\ne\varnothing$ for a face $f$, this intersection contains a common vertex $v$. Both $v$ and $\tau v$ are then vertices of $f$. They are distinct and have opposite colors, so they are adjacent on the quadrilateral boundary, a contradiction. Hence every closed face is disjoint from its image.

Every point of the surface belongs to a closed face, so the action is free. The quotient is a closed surface, and the quotient projection embeds every closed cell. Its face images are therefore disks bounded by four distinct edges and vertices, giving a regular quadrangulation.

Conversely, if an edge joins $v$ to $\tau v$, its endpoints are identified in the quotient. Under a free action its image is a loop, contrary to regularity.
\end{proof}

On $S^2$, a free involution is necessarily orientation-reversing: otherwise its quotient would be a closed orientable surface with Euler characteristic $2/2=1$, which is impossible. Thus, on the sphere, conditions \textup{(ii)} and \textup{(iii)} already imply \textup{(i)}, and the quotient is $\mathbb{RP}^2$.

\begin{theorem}[Self-polar descent]
\label{thm:cover-characterization}
A connected bipartite quadrangulation $\widetilde Q$ of a closed orientable surface $\widetilde S$ is the orientation double cover of an orientation-twisted quadrangulation of a nonorientable surface if and only if it admits an admissible involution $\tau$. For an FP tiling $\widetilde T$ of $\widetilde Q$, every realization self-polar with respect to such a $\tau$ descends to a realization of a polarity tiling $T$ on the quotient $Q=\widetilde Q/\langle\tau\rangle$, preserving face coherence.
\end{theorem}

\begin{proof}
For an orientation cover, the deck involution reverses orientation, and color exchange follows from Proposition~\ref{prop:cover-coloring}. No edge can join $\widetilde v$ to $\tau\widetilde v$, since its projection would be a loop in the regular quadrangulation downstairs. Thus the deck involution is admissible.

Conversely, Lemma~\ref{lem:quotient-regularity} gives a free action with a regular quadrangulation quotient $Q:=\widetilde Q/\langle\tau\rangle$. The projection is a two-sheeted covering $\widetilde S\to S:=\widetilde S/\langle\tau\rangle$. Since $\tau$ reverses orientation, $S$ is nonorientable and this is its orientation cover. To check that $Q$ is orientation-twisted, lift a closed edge path $\gamma$ from a vertex $\widetilde v$. An orientation-preserving path lifts to a closed path, which has even length in the bipartite graph. If $\gamma$ reverses orientation, its lift ends at $\tau\widetilde v$, which has a different color from $\widetilde v$. Each edge of $\gamma$ changes the color once along the lift, so $\gamma$ has odd length. Hence $[\lambda]=w_{1,S}$.

Let $\widetilde T$ be an FP tiling of $\widetilde Q$, with labeling $\widetilde{\mathcal L}$ and a realization $\widetilde R$ self-polar with respect to $\tau$. Each vertex orbit has a unique color-zero representative $\widetilde v$. Define the quotient labeling and realization by
\[
  \mathcal L([\widetilde v]):=
  \widetilde{\mathcal L}(\widetilde v),\qquad
  R([\widetilde v]):=\widetilde R(\widetilde v)
  \qquad(c(\widetilde v)=0).
\]
The quotient labels are point labels, and $R$ is well defined on these formal labels. They define a polarity tiling $T$ of $Q$. Self-polarity gives $\widetilde R(\tau\widetilde v)=\widetilde R(\widetilde v)^\perp$, so the two objects assigned to each orbit convert to the same point. Proposition~\ref{prop:polarity-conversion} identifies the coherence of each face in $\widetilde T$ with that of its quotient face in $T$. The marked faces $\widetilde f$ and $\tau\widetilde f$ of $\widetilde T$ descend to a single marked face of $T$. Lifting the realization of $T$ with the original two-coloring recovers $\widetilde R$ at every vertex.
\end{proof}

The quotient has nonorientable genus $g+1$ if $\widetilde S$ has orientable genus $g$, and half as many vertices, edges, and faces. Two admissible involutions of $\widetilde Q$ yield isomorphic quotient quadrangulations if and only if they are conjugate in $\Aut_{\mathrm{map}}(\widetilde Q)$.

\section{Hermitian Polarity over Division Rings}
\label{sec:division-rings}

The preceding master theorem relies on cancellation of edge pairings. Over a noncommutative division ring, the order of multiplication matters, and this cancellation argument no longer applies directly. Hermitian polarity introduces a further change: reversing a pairing involves the coefficient involution, even when the division ring is a field.

Izosimov's approach expresses local coherence as trivial facial holonomy \cite{Izosimov2026}. This suggests studying when flatness on all faces except one forces flatness on the remaining face. The gauge correspondence makes this question accessible through the fundamental group and the boundary word of the punctured surface. We use this perspective to study Hermitian polarity, where both the coefficient involution and the orientability of the surface affect the result.

\subsection{Gauge Connections and Boundary Holonomy}

We first recall the gauge-theoretic mechanism used by Izosimov. Let $\Gamma$ be a connected graph, let $\vec E(\Gamma)$ be its set of oriented edges, and write $\bar e$ for the reversal of $e$. For a group $G$, a \emph{$G$-connection} is a map $U:\vec E(\Gamma)\to G$ satisfying $U(\bar e)=U(e)^{-1}$. Its \emph{holonomy} along an edge path $\gamma=e_1\cdots e_r$ is $\Hol_U(\gamma)=U(e_1)\cdots U(e_r)$. A vertex function $q:V(\Gamma)\to G$ acts by
\[
  (q\cdot U)(e)
  =q(t(e))U(e)q(h(e))^{-1}.
\]
Connections in the same orbit are \emph{gauge-equivalent}.

Suppose now that $\Gamma$ is a map on a compact connected surface $S$, so the closures of the complementary components are the faces of a cellular decomposition. A connection is \emph{flat} if its holonomy around every face is trivial. Changing the base vertex of a facial boundary conjugates its holonomy, while reversing the boundary orientation inverts it, so flatness is independent of these choices. This is Izosimov's terminology for connections on maps \cite[Definitions~3.1--3.7 and Remark~3.8]{Izosimov2026}.

The following formulation combines results proved by Izosimov in Propositions~3.5, 3.9, and~3.10 \cite{Izosimov2026}.

\begin{theorem}[Gauge correspondence]
\label{thm:gauge-representations}
Fix $v_0\in V(\Gamma)$. Holonomy induces natural bijections
\[
  \frac{\{\text{$G$-connections on }\Gamma\}}
       {\{\text{gauge equivalence}\}}
  \cong
  \Hom\bigl(\pi_1(\Gamma,v_0),G\bigr)/G
\]
and
\[
  \frac{\{\text{flat $G$-connections on }\Gamma\}}
       {\{\text{gauge equivalence}\}}
  \cong
  \Hom\bigl(\pi_1(S,v_0),G\bigr)/G,
\]
where $G$ acts on representations by conjugation.
\end{theorem}

The gauge correspondence translates this flatness implication into a condition on boundary holonomy.

\begin{theorem}
\label{thm:group-valued-excision}
Let $Q$ be a regular quadrangulation of a closed connected surface $S$, with a marked face $f_0$, and let $G$ be a group. Put $Y:=S\setminus\operatorname{Int}(f_0)$. Choose a boundary vertex $x_0\in V(Q)\cap\partial f_0$ and an orientation of the boundary edge loop of $f_0$, and let $\delta_{f_0}\in\pi_1(Y,x_0)$ be its class. The following conditions are equivalent:
\begin{enumerate}[label=\textup{(\roman*)}]
  \item Every $G$-connection on $Q^{(1)}$ that is flat on every face $f\ne f_0$ is also flat on $f_0$.
  \item Every homomorphism $\rho:\pi_1(Y,x_0)\to G$ satisfies $\rho(\delta_{f_0})=1$.
\end{enumerate}
\end{theorem}

Changing the direction of the boundary loop inverts $\delta_{f_0}$, while changing the base vertex conjugates it. Consequently, condition \textup{(ii)} is independent of these choices; the base vertex will henceforth be suppressed.

\begin{proof}
Assume \textup{(i)}, and let $\rho:\pi_1(Y,x_0)\to G$ be a homomorphism. The graph $Q^{(1)}\subset Y$ is a map whose faces are precisely the faces of $Q$ distinct from $f_0$. By Izosimov's gauge correspondence, Theorem~\ref{thm:gauge-representations}, there is a flat $G$-connection $U$ on this map whose holonomy representation is conjugate to $\rho$. Regarded as a connection on $Q^{(1)}$, it is flat on every face other than $f_0$, and hence also on $f_0$ by \textup{(i)}. Thus $\rho(\delta_{f_0})=1$, since conjugation does not affect whether an element is trivial.

Conversely, let $U$ be a $G$-connection on $Q^{(1)}$ that is flat on every face distinct from $f_0$. Since $U$ is flat on the map $Q^{(1)}\subset Y$, its graph holonomy factors through a homomorphism $\rho_U:\pi_1(Y,x_0)\to G$ by \cite[Proposition~3.9]{Izosimov2026}. The holonomy of $U$ around $f_0$ is $\rho_U(\delta_{f_0})$, up to inversion and conjugation, and is therefore trivial by \textup{(ii)}. This proves \textup{(i)}.
\end{proof}

Because face coherence is equivalent to flatness of an associated connection, these conditions are sufficient for the corresponding master theorem. The converse need not hold, since geometric realizations may produce only a subclass of all $G$-connections.

\subsection{The Hermitian Polarity Master Theorem}
\label{subsec:hermitian-master}

Izosimov proves that flatness on all but one face forces flatness on the remaining face in two cases: on the sphere, for any coefficient group $G$, and on any closed connected orientable surface, for an abelian coefficient group $G$ \cite[Corollaries~3.11--3.12]{Izosimov2026}. For nonorientable surfaces, however, the boundary word is a product of squares rather than commutators, and abelianity alone does not make these squares trivial. We therefore keep track of orientation reversal in the coefficient group, using a semidirect product that also encodes Hermitian symmetry.

Use the conventions for $\mathrm{D},V,\sigma,h$ from Section~\ref{subsec:projective-polarities}. The map $\widehat\sigma(g):=\sigma(g)^{-1}$ is an involutive automorphism of $\mathrm{D}^\times$. Put
\begin{equation}\label{eq:semidirect}
  \mathrm{D}^\times_\sigma:=\mathrm{D}^\times\rtimes_{\widehat\sigma}\Ftwo,
  \qquad
  (g,a)(k,b)=\bigl(g\widehat\sigma^a(k),a+b\bigr).
\end{equation}
Let $p:\mathrm{D}^\times_\sigma\to\Ftwo$ be the projection onto the second coordinate.

For a polarity tiling with realization $R$, choose vector representatives of its realized point labels and write $H(v,w)$ for their $h$-pairing at vertices $v,w$. We likewise write $H(A,C)$ when the vertices are denoted by their point labels.

\begin{definition}[Hermitian connection]\label{def:hermitian-connection}
On a polarity tiling whose point labels at adjacent vertices are nonconjugate, the \emph{Hermitian connection} is given on each oriented edge by
\[
  U_h(v,w):=\bigl(H(v,w),1\bigr)\in \mathrm{D}^\times_\sigma.
\]
\end{definition}

Hermitian symmetry and $(g,1)^{-1}=(\sigma(g),1)$ give $U_h(w,v)=U_h(v,w)^{-1}$. Rescaling a representative at $v$ by $a_v\in \mathrm{D}^\times$ acts by the gauge transformation $q(v)=(\sigma(a_v),0)$.

For a face with boundary $A-C-B-D-A$ whose consecutive point labels are nonconjugate, set
\begin{equation}\label{eq:chi-h}
  \chi_h(A,B;C,D)
  :=H(A,C)H(B,C)^{-1}H(B,D)H(A,D)^{-1}.
\end{equation}
The face is \emph{$h$-coherent} when $\chi_h(A,B;C,D)=1$. Thus an $h$-coherent face has nonconjugate consecutive labels, just as in the symmetric bilinear setting.

\begin{proposition}
\label{prop:hermitian-connection}
The following properties hold.
\begin{enumerate}[label=\textup{(\arabic*)}]
  \item For a face with boundary $A-C-B-D-A$,
  \begin{equation}\label{eq:hermitian-cross-ratio}
    \Hol_{U_h}(\partial f)=\bigl(\chi_h(A,B;C,D),0\bigr).
  \end{equation}
  \item Every edge path $\gamma$ satisfies
  \[
    p\bigl(\Hol_{U_h}(\gamma)\bigr)=|\gamma|\pmod2.
  \]
  For an orientation-twisted tiling, this equals $\langle w_{1,S},[\gamma]\rangle$ on closed paths.
  \item For every $q\in \mathrm{D}^\times$,
  \begin{equation}\label{eq:hermitian-odd-square}
    (q,1)^2=\bigl(q\sigma(q)^{-1},0\bigr).
  \end{equation}
  Thus every odd element has square equal to the identity when $\sigma=\id$.
\end{enumerate}
\end{proposition}

\begin{proof}
For \textup{(1)}, multiplication around the face uses $\widehat\sigma(H(C,B))=H(B,C)^{-1}$ and $\widehat\sigma(H(D,A))=H(A,D)^{-1}$, giving \eqref{eq:hermitian-cross-ratio}. For independent opposite pairs, the geometric condition is equivalent to a nonzero right kernel of
\[
  \begin{pmatrix}H(A,C)&H(A,D)\\H(B,C)&H(B,D)\end{pmatrix}.
\]
Eliminating the first variable gives $H(B,D)=H(B,C)H(A,C)^{-1}H(A,D)$, equivalently $\chi_h=1$. For \textup{(2)}, path parity follows by adding the second coordinates; on closed paths, orientation-twistedness gives $[\lambda]=w_{1,S}$. Property \textup{(3)} follows directly from \eqref{eq:semidirect}.
\end{proof}

\begin{remark}[Reduction on orientable surfaces]
For an orientation-twisted tiling on an orientable surface, closed-path holonomy lies in $\mathrm{D}^\times\times\{0\}$. To obtain a connection taking values in this subgroup on every edge, choose a two-coloring $c$ and apply the gauge transformation $q(v)=(1,c(v))$. Then
\[
  (q\cdot U_h)(v,w)=
  \begin{cases}
    (H(v,w),0),&c(v)=0,\\
    (H(w,v)^{-1},0),&c(v)= 1,
  \end{cases}
\]
which is equivalent to the setting of Fomin and Pylyavskyy\cite[Theorem 2.6]{FominPylyavskyy2023}.
\end{remark}

We consider the polarity implication associated with an orientation-twisted quadrangulation $Q$ with a marked face $f_0$: coherence on every other face implies coherence on $f_0$. Universal validity means validity for every assignment of projective points to the vertices of $Q$ with nonconjugate adjacent points.

\begin{theorem}[Hermitian polarity master theorem]
\label{thm:hermitian-polarity-master}
If a polarity implication is represented by a tiling of the sphere, then it holds for arbitrary $\mathrm{D}$ and $\sigma$. If the tiling is on a closed orientable surface of positive genus, then it holds for any field $\mathrm{D}$ and arbitrary $\sigma$. If the tiling is on a closed nonorientable surface, then it holds for any field $\mathrm{D}$ with $\sigma=\id$.
\end{theorem}

\begin{proof}
Let $Q$ be the quadrangulation of $S$ underlying the tiling, with a marked face $f_0$, and fix a point realization for which every face other than $f_0$ is $h$-coherent. Every edge belongs to two distinct faces, at least one of which is a premise face. Since consecutive labels of an $h$-coherent face are nonconjugate, all edge pairings are nonzero and $U_h$ is defined on the whole one-skeleton. Put $Y=S\setminus\Int(f_0)$, and let $i:Y\hookrightarrow S$ be the inclusion. By Proposition~\ref{prop:hermitian-connection}, $U_h$ is flat on every face of $Y$. By Theorem~\ref{thm:gauge-representations}, it therefore determines a representation $\rho:\pi_1(Y)\to \mathrm{D}^\times_\sigma$. The path-parity formula and the orientation-twisted identity $[\lambda]=w_{1,S}$ show that $p(\rho(\gamma))=\langle i^*w_{1,S},[\gamma]\rangle$ for every $\gamma\in\pi_1(Y)$. Thus we apply the boundary-holonomy criterion to representations with the prescribed parity $p\circ\rho=i^*w_{1,S}$. The holonomy around $f_0$ is $\rho(\delta_{f_0})$, up to inversion and conjugation, so it remains to prove that this element is trivial. We use the standard boundary words for punctured surfaces \cite[Section~1.2]{Hatcher2002}.

If $S=S^2$, then $Y$ is a disk and $\delta_{f_0}=1$. This proves the spherical case.

Let $S=\Sigma_g$ with $g\ge1$. Since $S$ is orientable, $w_{1,S}=0$, so the image of $\rho$ is contained in the even subgroup $\mathrm{D}^\times\times\{0\}$. Write $\rho(a_i)=(\alpha_i,0)$ and $\rho(b_i)=(\beta_i,0)$. Since $\delta_{f_0}=\prod_{i=1}^g[a_i,b_i]$, the orientable boundary word gives
\[
  \rho(\delta_{f_0})
  =
  \left(\prod_{i=1}^g[\alpha_i,\beta_i],0\right).
\]
When $\mathrm{D}$ is commutative, every commutator in $\mathrm{D}^\times$ is trivial, which proves the orientable case. Notice that no condition on $\sigma$ is needed: only the even subgroup occurs.

Finally, let $S=N_k$. The standard generators $x_j$ are orientation-reversing, hence $\rho(x_j)=(q_j,1)$ for some $q_j\in \mathrm{D}^\times$, and $\delta_{f_0}=\prod_{j=1}^k x_j^2$. From \eqref{eq:hermitian-odd-square}, each square $\rho(x_j)^2$ is the identity when $\sigma=\id$. The nonorientable boundary word therefore gives $\rho(\delta_{f_0})=1$, proving the nonorientable case. Since the identity map can be an anti-automorphism only on a commutative division ring, commutativity here is automatic from $\sigma=\id$.

In all three cases $U_h$ is flat around $f_0$. Its face holonomy is therefore trivial, and Proposition~\ref{prop:hermitian-connection} implies that $f_0$ is $h$-coherent.
\end{proof}

\begin{remark}[The field with three elements]
\label{rem:hermitian-f3}
Over $\mathrm{D}=\mathbb F_3$, the conclusion of Theorem~\ref{thm:hermitian-polarity-master} holds for any quadrangulation $Q$ of a closed connected surface, without the orientation-twisted condition. In particular, no bipartiteness assumption is needed. To see this, fix a realization with nonconjugate adjacent points. The involution of $\mathbb F_3$ is necessarily the identity, so $h$ is symmetric bilinear. Choose vector representatives $x_v$ of the realized points, and for each edge $e=uv$ put
\[
 s_e:=h(x_u,x_v)\in\mathbb F_3^\times=\{1,-1\}.
\]
Symmetry makes $s_e$ independent of the direction of the edge, and $s_e^{-1}=s_e$. Thus \eqref{eq:chi-h} gives, for every face $f$,
\[
 \chi_h(f)=\prod_{e\in\partial f}s_e,
\]
independently of the initial vertex and the direction in which its boundary is read. Since each edge belongs to exactly two faces,
\[
 \prod_{f\in F(Q)}\chi_h(f)
   =\prod_{e\in E(Q)}s_e^2=1.
\]
Consequently, coherence of all but one face forces coherence of the remaining face.

The same identity $s^2=1$ underlies the observation of Pylyavskyy and Skopenkov for triangulated closed nonorientable surfaces \cite[discussion following Proposition~4.21]{PylyavskyySkopenkov2026}. On a bipartite quadrangulation, replacing the realized points at vertices of one color by their polar hyperplanes identifies the face conditions with those of the FP construction, by Proposition~\ref{prop:polarity-conversion}. The argument above also applies to quadrangulations whose one-skeleton is not bipartite.
\end{remark}

Outside the coefficient cases of Theorem~\ref{thm:hermitian-polarity-master}, even an orientation-twisted tiling can fail to give a valid implication because the boundary holonomy need not be trivial.

\begin{example}[Failure of the complex Hermitian Hesse analogue]\label{ex:hermitian-hesse}
Let
\[
  \mathrm{D}=\mathbb C,\qquad \sigma(z)=\overline z,\qquad
  h(x,y)=\overline x^{\,T}y,
\]
and take
\[
  A=(1,0,0),\quad B=(1,1,0),\quad
  C=(1,i,1),\quad D=(1,i,-1-i).
\]
All six pairings between distinct points are nonzero. For the three opposite-line conditions in the Hesse configuration of Example~\ref{ex:hesse-desargues-lift}, the relevant cross ratios are
\[
\begin{aligned}
 \chi_h(A,B;C,D)
 &=1\,(1+i)^{-1}(1+i)\,1^{-1}=1,\\
 \chi_h(A,C;B,D)
 &=1\,(1-i)^{-1}(1-i)\,1^{-1}=1,\\
 \chi_h(A,D;B,C)
 &=1\,(1-i)^{-1}(1+i)\,1^{-1}=i\ne1.
\end{aligned}
\]
Thus $AB\perp_h CD$ and $AC\perp_h BD$, but $AD\not\perp_h BC$. This is an explicit $h$-realizable obstruction, not merely an abstract representation-theoretic one. Equivalently, deleting the third hemicube face gives a M\"obius band with $\delta_{f_0}=x^2$. For the odd core generator $x=[A\,B\,C\,A]$, direct multiplication gives
\[
  \rho(x)
  =(1,1)(1+i,1)(1,1)
  =\bigl((1+i)/2,1\bigr).
\]
Writing $p=(1+i)/2$, we obtain
\[
  \rho(\delta_{f_0})=(p\overline p^{-1},0)=(i,0)\ne1.
\]
\end{example}

\subsection{Realization of Hermitian connections}
\label{subsec:hermitian-realization}

To prove the converse in Theorem~\ref{thm:hermitian-extensions-intro} and turn holonomy obstructions into geometric counterexamples, we need to realize prescribed connections using a fixed Hermitian form. Izosimov's realization argument chooses point vectors and hyperplane covectors independently \cite[Theorem~5.17 and its second proof]{Izosimov2026}. In the polarity setting, we instead choose the point vectors successively, prescribing their pairings with previously chosen neighbors. The following theorem gives sufficient dimension bounds for this construction in terms of graph degeneracy.

\begin{definition}[Graph degeneracy]\label{def:graph-degeneracy}
The \emph{degeneracy} of a finite graph $\Gamma$ is
\[
  d(\Gamma):=
  \max_{\substack{\Gamma'\subseteq\Gamma\\V(\Gamma')\ne\varnothing}}
  \min_{v\in V(\Gamma')}\deg_{\Gamma'}(v).
\]
Equivalently, $d(\Gamma)$ is the least integer $d$ for which the vertices admit an ordering in which each vertex has at most $d$ earlier neighbors.
\end{definition}

\begin{theorem}[Hermitian realization]\label{thm:hermitian-realization}
Let $\Gamma$ be a finite connected simple graph with $n\ge2$ vertices and degeneracy $d$. Let $h$ be a nondegenerate $\sigma$-Hermitian form on a right $\mathrm{D}$-vector space $V$, and put $r=\dim\mathbb P(V)$, so that $\dim_{\mathrm{D}} V=r+1$. Assume that
\[
  \begin{cases}
    r\ge2d-1,&\text{if $\mathrm{D}$ is infinite},\\[2pt]
    r\ge2d-2+\lceil\log_s n\rceil,&\text{if $|\mathrm{D}|=s<\infty$}.
  \end{cases}
\]
For every $\mathrm{D}^\times_\sigma$-connection $U$ on $\Gamma$ with $p(U(e))=1$ on every oriented edge, there are nonzero vectors $a_v\in V$, with pairwise distinct projective points $[a_v]$, such that
\[
  U(v,w)=\bigl(h(a_v,a_w),1\bigr)
  \qquad\text{on every oriented edge }(v,w).
\]
In particular, the resulting polarity realization has nonconjugate point labels at adjacent vertices and Hermitian connection exactly $U$.
\end{theorem}

\begin{proof}
Choose an ordering $v_1,\ldots,v_n$ with at most $d$ earlier neighbors at each vertex, and put
\[
  I_k:=\{j<k:v_j\sim v_k\},\qquad N:=r+1.
\]
Write $U(v_j,v_k)=(u_{jk},1)$. The connection identity gives $u_{kj}=\sigma(u_{jk})$, so prescribing one direction of each pairing also prescribes its reverse.

We inductively choose $a_1,\ldots,a_i$ with the following properties:
\begin{enumerate}[label=\textup{(\arabic*)}]
  \item The projective points $[a_1],\ldots,[a_i]$ are pairwise distinct.
  \item For every edge with both endpoints already chosen, $h(a_j,a_k)=u_{jk}$.
  \item For each $k>i$, the vectors $\{a_j:j\in I_k\cap\{1,\ldots,i\}\}$ are right-linearly independent.
\end{enumerate}
Suppose that vectors through $a_{i-1}$ have been chosen. We seek $a_i$ satisfying the $t:=|I_i|\le d$ equations $h(a_j,a_i)=u_{ji}$ for $j\in I_i$. The right-linear functionals $h(a_j,-)$, $j\in I_i$, are left-linearly independent. Indeed, a relation $\sum_{j\in I_i}c_jh(a_j,-)=0$ gives $h\left(\sum_{j\in I_i}a_j\sigma(c_j),-\right)=0.$ Nondegeneracy of $h$ and the induction hypothesis imply that every $c_j$ is zero. Thus these equations have an affine solution space $\mathcal A_i$ of dimension $N-t$. If $t>0$, its nonzero right-hand sides give $0\notin\mathcal A_i$; if $t=0$, then $\mathcal A_i=V$.

To preserve distinctness and independence, avoid
\[
  \mathcal B_i:=
  \left(\bigcup_{j<i}a_j\mathrm{D}\right)
  \cup
  \left(\bigcup_{\substack{k>i\\v_i\sim v_k}}
    \Span_{\mathrm{D}}\{a_j:j\in I_k\cap\{1,\ldots,i-1\}\}\right).
\]
This is a union of at most $n-1$ linear subspaces, each of dimension at most $b:=\max\{1,d-1\}$. For the second family, $i\in I_k$ leaves at most $d-1$ earlier neighbors. Also $0\in\mathcal B_i$: for $i=1$, connectedness provides a later neighbor, and for $i>1$ the first family is nonempty.

Suppose first that $\mathrm{D}$ is infinite. The bound gives $N\ge2d$, hence $\dim\mathcal A_i=N-t\ge d\ge b$. If $t>0$, no component of $\mathcal B_i$ contains $\mathcal A_i$: equality of dimensions would force the affine space to be a linear subspace containing zero. If $t=0$, then $\dim\mathcal A_i=N>b$. In either case, its intersections with the components are proper affine subspaces. Choose a direction outside their direction subspaces \cite[Lemma~4.7]{Izosimov2026}; a line in this direction meets each intersection in at most one point. Since $\mathrm{D}$ is infinite, $\mathcal A_i\setminus\mathcal B_i\ne\varnothing$.

Now suppose that $|\mathrm{D}|=s<\infty$, and put $m:=\lceil\log_s n\rceil$. Then $N\ge2d-1+m$. If $t>0$, each nonempty intersection of $\mathcal A_i$ with a component of $\mathcal B_i$ has dimension at most $b-1$, since $0\notin\mathcal A_i$. As $b-1\le d-1$,
\[
  |\mathcal A_i\cap\mathcal B_i|
  \le(n-1)s^{b-1}
  <s^{m+d-1}
  \le s^{N-d}
  \le|\mathcal A_i|.
\]
If $t=0$, then $\mathcal A_i=V$, and $b\le2d-1$ gives
\[
  |\mathcal B_i|\le(n-1)s^b<s^{m+b}\le s^N=|\mathcal A_i|.
\]
Thus again we may choose $a_i\in\mathcal A_i\setminus\mathcal B_i$. This vector is nonzero, satisfies all prescribed pairings, and preserves the three induction properties. The induction yields the required realization.
\end{proof}

\begin{remark}\label{rem:hermitian-realization}
\begin{enumerate}[label=\textup{(\arabic*)}]
  \item The proof of Theorem~\ref{thm:hermitian-realization} uses only the graph, without any two-cells or two-coloring, and is therefore independent of orientability.
  \item If $\mathcal L$ is not injective, identify equally labeled vertices and apply the construction to the resulting graph of labels, provided it is loopless and the prescribed edge values descend to it. This identification need not preserve a quadrangulation. For example, identifying repeated labels in the cube lift of the Hesse tiling and merging parallel edges gives $K_4$.
  \item Parallel edges impose the same pairing equation and may be merged whenever their prescribed connection values agree after the edges are oriented in the same direction.
  \item Beyond compatibility of the reverse-edge values, the inductive argument uses only linearity and nondegeneracy of the pairing. The same method therefore applies to analogous master theorems whose realization conditions are prescribed values of a nondegenerate pairing, such as the M\"obius geometry described by Fomin and Pylyavskyy \cite[Real cross-ratios]{FominPylyavskyy2023}.
\end{enumerate}
\end{remark}

\begin{proof}[Proof of Theorem~\ref{thm:hermitian-extensions-intro}]
The sufficient cases follow from Theorem~\ref{thm:hermitian-polarity-master}. For necessity, let $Q$ be an orientation-twisted quadrangulation of a closed connected surface $S$, with a marked face $f_0$ and a simple one-skeleton, and fix a nondegenerate $\sigma$-Hermitian form $h$ on $V$. Assume that $r=\dim\mathbb P(V)$ satisfies the bound in Theorem~\ref{thm:hermitian-realization} for $Q^{(1)}$. Put $Y=S\setminus\Int(f_0)$, and let $i:Y\hookrightarrow S$ be the inclusion. Outside the three listed cases, we construct a representation
\[
  \rho:\pi_1(Y)\longrightarrow \mathrm{D}^\times_\sigma,
  \qquad p\circ\rho=i^*w_{1,S},
  \qquad \rho(\delta_{f_0})\ne1.
\]
We use the same standard generators and boundary words as in the proof of Theorem~\ref{thm:hermitian-polarity-master}. If $S\cong\Sigma_g$ with $g\ge1$ and $\mathrm{D}$ is noncommutative, choose $\alpha,\beta\in \mathrm{D}^\times$ with $[\alpha,\beta]\ne1$. Send $a_1,b_1$ to $(\alpha,0),(\beta,0)$, respectively, and all remaining generators to the identity. Then $\rho(\delta_{f_0})=\bigl([\alpha,\beta],0\bigr)\ne1.$ If $S\cong N_k$ and $\sigma\ne\id$, choose $q\in \mathrm{D}^\times$ with $\sigma(q)\ne q$. Send $x_1$ to $(q,1)$ and all remaining generators to $(1,1)$. Then $\rho(\delta_{f_0})=\bigl(q\sigma(q)^{-1},0\bigr)\ne1.$ Both representations have the required parity character.

By Theorem~\ref{thm:gauge-representations}, $\rho$ is represented, up to conjugation, by a connection $U_0$ on $Q^{(1)}$ that is flat on every face other than $f_0$. Since $Q$ is orientation-twisted, the cochains $\alpha_0(e):=p(U_0(e))$ and $\lambda(e):=1$ have the same sums along every closed edge path. Hence $\alpha_0+\lambda=\mathrm dc$ for a vertex function $c:V(Q)\to\Ftwo$. The gauge transformation $q_0(v):=(1,c(v))$ gives a connection $U:=q_0\cdot U_0$ with $p(U(e))=1$ on every oriented edge. Gauge transformations preserve facial flatness, so $U$ is flat on every premise face and has nontrivial holonomy around $f_0$.

Theorem~\ref{thm:hermitian-realization} realizes $U$ using the fixed form $h$, with pairwise distinct projective points. By Proposition~\ref{prop:hermitian-connection}, the resulting injective labeling makes every premise face $h$-coherent and $f_0$ noncoherent.
\end{proof}

\section{Tiling Surface}
\label{sec:tiling-surface}

The preceding sections start with a tiling and study the incidence implication it carries. We now consider the inverse problem: given a formal polarity presentation $I:Q_1,\ldots,Q_m\Longrightarrow Q_0$, does it admit a tiling proof, and how can one construct a proof of least Euler genus? We consider orientation-twisted tiling proofs, with one conclusion face and all other faces labeled by the given premises; premise faces may be repeated or omitted, but no auxiliary labels are introduced.

We adopt Zhang's framework of fixed-form implications and elementary surface-tiling proofs \cite[Sections~2.3 and~3]{ZhangTiledGenus}. We extend it to polarity presentations by allowing nonbipartite label graphs and nonorientable source surfaces, with path parity required to agree with the orientation character.

\subsection{Premise Complexes and Tiling Proofs}
\label{subsec:premise-complexes}

The following definition adapts Zhang's fixed-form FP implications \cite[Definitions~2.8--2.9]{ZhangTiledGenus} to formal boundary labels in a finite alphabet. It applies both to point--hyperplane presentations and to polarity presentations.

\begin{definition}[Formal presentation]
\label{def:formal-presentation}
Let $\mathcal A$ be a finite alphabet. A \emph{formal square} over $\mathcal A$ is an equivalence class of quadruples $(a_0,a_1,a_2,a_3)\in\mathcal A^4$ under cyclic permutation and reversal, with $a_j\ne a_{j+1}$ for all $j$ modulo four. We write $\partial Q=a_0a_1a_2a_3a_0$ for a representative boundary word. A \emph{formal presentation} over $\mathcal A$ is a finite ordered list
\[
 I=(Q_0,Q_1,\ldots,Q_m),\qquad m\ge0,
\]
of formal squares such that every symbol of $\mathcal A$ occurs in the list and the four labels of $Q_0$ are pairwise distinct. The squares $Q_1,\ldots,Q_m$ are called the \emph{premise squares}, and $Q_0$ is called the \emph{conclusion square}. We also denote this list by
\[
 I:Q_1,\ldots,Q_m\Longrightarrow Q_0.
\]
Repetitions among the $Q_i$, and repeated opposite labels in a premise square, are allowed.

A \emph{formal polarity presentation} uses all symbols of $\mathcal A$ as point variables. A \emph{formal point--hyperplane presentation} additionally includes a partition $\mathcal A=\mathcal P\sqcup\mathcal H$ into point and hyperplane variables, such that each boundary word alternates between the two types.
\end{definition}

These are formal data: no surface, geometric realization, or gluing of edges is part of a presentation. The squares are therefore not assumed to be faces of a preassigned tiling, and no orientability condition is imposed at this stage. The arrow is notation for the ordered list with its distinguished conclusion.

For the geometric interpretation, fix the field, projective space, and, in the polarity case, the nondegenerate symmetric bilinear form. Use the realizations and coherence conditions introduced above (Section~\ref{subsec:topological-conventions}, Theorem~\ref{thm:fp-master}, and Definition~\ref{def:polarity-labeling}). The nonincidence or nonconjugacy requirement applies to consecutive labels in every $Q_i$, including $Q_0$. Repeated formal labels have the same image; distinct labels need not have distinct images. The associated incidence statement asserts that, for every such realization, coherence of all realized premise squares implies coherence of the realized conclusion square. When $m=0$, it asserts coherence of $Q_0$ in every realization. The truth of this statement is not assumed as part of the formal data.

Following Zhang \cite[Definition~3.1]{ZhangTiledGenus}, we associate to $I$ a premise complex, a premise group, and a conclusion class. Here the label graph need not be bipartite.

\begin{definition}[Premise complex]
\label{def:premise-complex}
Let $\Gamma_I$ have the occurring labels as vertices and one undirected edge for each unordered pair of consecutive labels in any $\partial Q_i$, including $\partial Q_0$. Attach one two-cell $D_i$ along each premise word, and set
\[
 X_I=\Gamma_I\cup_{\partial Q_1}D_1\cup\cdots
                  \cup_{\partial Q_m}D_m.
\]
This is the \emph{premise complex} of $I$. No two-cell is attached for the conclusion $Q_0$.
\end{definition}

\begin{figure}[htbp]
\centering
\begin{tikzpicture}[
    x=1cm,
    y=1cm,
    premise/.style={
        draw=black,
        line width=0.9pt,
        fill=blue!7
    },
    conclusion/.style={
        draw=red!75!black,
        line width=1.6pt,
        dash pattern=on 4pt off 2.5pt
    },
    vertex/.style={
        circle,
        fill=black,
        inner sep=1.35pt
    }
]
    \filldraw[premise] (0,0) rectangle (1.6,1.35);
    \filldraw[premise] (1.6,0) rectangle (3.2,1.35);
    \filldraw[premise] (0,1.35) rectangle (1.6,2.70);
    \draw (1.6,1.35)--(3.2,1.35)
        --(3.2,2.70)--(1.6,2.70)--cycle;
    \draw[conclusion] (1.6,1.35)--(3.2,1.35)
        --(3.2,2.70)--(1.6,2.70)--cycle;
    \foreach \x in {0,1.6,3.2}
        \foreach \y in {0,1.35,2.70}
            \node[vertex] at (\x,\y) {};
    \node at (0.8,0.675) {$D_1$};
    \node at (2.4,0.675) {$D_2$};
    \node at (0.8,2.025) {$D_3$};
    \node[text=red!75!black] at (2.4,2.025) {$\partial Q_0$};
\end{tikzpicture}
\caption{A schematic premise complex with three attached premise cells. The dashed square marks the conclusion boundary; its interior is not a two-cell of $X_I$.}
\label{fig:premise-complex}
\end{figure}
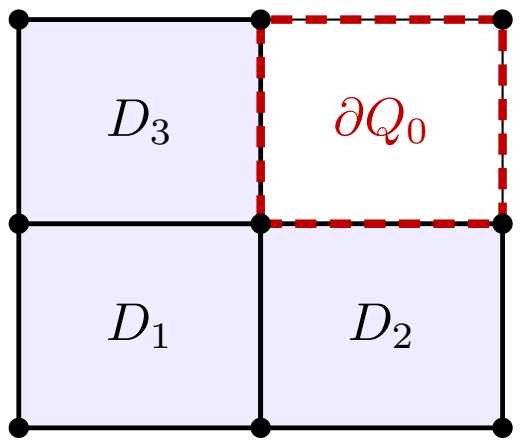

Figure~\ref{fig:premise-complex} illustrates the distinction between the attached premise cells and the retained conclusion boundary. The drawing is schematic: a general premise complex need not be a surface, and repeated labels can make its attaching maps noninjective.

Let $X_I^0$ be the component containing $\partial Q_0$, and choose a basepoint $x_0\in\partial Q_0$. The \emph{premise group} is $G_I=\pi_1(X_I^0,x_0).$ The conclusion boundary determines an element up to conjugacy and inversion. We denote this unoriented conjugacy class by $\delta_I$ and call it the \emph{conclusion class}.

All attached cells have boundary length four, so the parity construction applies to $X_I^0$. Write $\lambda_I=\lambda_{X_I^0}, \varepsilon_I=\varepsilon_{X_I^0}.$ The conclusion also has four sides, hence $\varepsilon_I(\delta_I)=0$. In a point--hyperplane presentation, $\Gamma_I$ is bipartite and $\varepsilon_I=0$. In a polarity presentation the character may be nonzero.

A tiling proof can be recorded by removing its conclusion face and mapping the remaining faces to the permitted premise cells. We adapt Zhang's surface-tiling diagrams \cite[Definition~3.2]{ZhangTiledGenus} to allow nonorientable source surfaces, and add the orientation-twisted condition below. We use finite square-cell structures: each face attaching boundary has four edge occurrences.

\begin{definition}[Surface-tiling diagram]
A \emph{surface-tiling diagram} for $I$ is a cellular map $F:Y\to X_I^0$, where $Y$ is a compact connected surface, orientable or nonorientable, with a finite square-cell structure and one boundary component, satisfying the following conditions:
\begin{enumerate}[label=\textup{(\roman*)},nosep]
\item Each edge maps homeomorphically to an edge, and each open face maps homeomorphically to an open premise cell.
\item The four boundary edges read $\partial Q_0$, up to cyclic shift and reversal.
\end{enumerate}
Every face of $Y$ must map to one of the premise cells $D_1,\ldots,D_m$; a given premise cell may be used more than once or omitted. The diagram is \emph{orientation-twisted} if
\[
 F^*[\lambda_I]=w_{1,Y},
 \qquad\text{equivalently}\qquad
 \varepsilon_I\circ F_*=w_{1,Y}.
\]
\end{definition}

An orientation-twisted surface-tiling diagram gives a tiling proof of $I$ of the kind discussed above: $F$ specifies the formal labels and premise squares, and filling the boundary with $Q_0$ supplies the conclusion face. The orientation-twisted condition is preserved. This extends Zhang's elementary surface-tiling proofs \cite[Definition~2.10]{ZhangTiledGenus}.

\FloatBarrier

\subsection{Parity Factorizations and Tiling Criteria}

In the ordinary bipartite setting, Zhang \cite[Theorem~1.9]{ZhangTiledGenus} proves that the least orientable genus of a tiling proof is
\begin{equation}\label{eq:zhang-genus}
  g_{\mathrm{or}}(I)=\cl_{G_I}(\delta_I).
\end{equation}
Here $\cl$ denotes commutator length. For polarity presentations, the orientation-twisted condition imposes parity restrictions on the factors, and nonorientable surfaces lead to products of squares.

The parity character $\varepsilon_I$ defined above records the length of a closed edge path modulo two. A group element is called \emph{even} or \emph{odd} according as its image under $\varepsilon_I$ is zero or one.

\begin{definition}[Square length and commutator length]
For $\gamma\in\ker\varepsilon_I$, define its \emph{odd square length} and \emph{even commutator length}, respectively, by
\begin{align}
  \slodd_{G_I}(\gamma)
  &:={}
  \min\left\{k\ge1\ \middle|\
  \gamma=z_1^2\cdots z_k^2,
  \ \varepsilon_I(z_j)=1\text{ for all }j\right\},
  \label{eq:odd-square-length}\\
  \cleven_{G_I}(\gamma)
  &:={}
  \min\left\{g\ge0\ \middle|\
  \gamma=[a_1,b_1]\cdots[a_g,b_g],
  \ \varepsilon_I(a_i)=\varepsilon_I(b_i)=0\text{ for all }i\right\}.
  \label{eq:even-commutator-length}
\end{align}
The empty commutator product is allowed. A minimum is $\infty$ when no such factorization exists.
\end{definition}

Conjugation preserves the parity of each factor, and inversion reverses the factor order using $(z^2)^{-1}=(z^{-1})^2$ and $[a,b]^{-1}=[b,a]$. Thus both lengths are well-defined on $\delta_I$.

For a surface-tiling diagram $F:Y\to X_I^0$, recall that the orientation-twisted condition is
\[
  F^*[\lambda_I]=w_{1,Y},
  \qquad\text{equivalently}\qquad
  \varepsilon_I\circ F_*=w_{1,Y}.
\]
Write $\Sigma_{g,1}$ for the orientable genus-$g$ surface with one boundary component, where $g\ge0$, and $N_{k,1}$ for the nonorientable genus-$k$ surface with one boundary component, where $k\ge1$.

The following theorem extends Zhang's commutator criterion \cite[Lemma~3.4]{ZhangTiledGenus} by imposing parity restrictions on the factors and including nonorientable surfaces.

\begin{theorem}\label{thm:parity-tiling-criterion}
For every formal presentation $I$, the following statements hold.
\begin{enumerate}[label=\textup{(\roman*)}]
  \item For $h\ge0$, an orientation-twisted orientable surface-tiling diagram for $I$ of genus at most $h$ exists if and only if
  \begin{equation}\label{eq:bounded-even-commutator-length}
    \cleven_{G_I}(\delta_I)\le h.
  \end{equation}
  \item For $h\ge1$, an orientation-twisted nonorientable surface-tiling diagram for $I$ of genus at most $h$ exists if and only if
  \begin{equation}\label{eq:bounded-odd-square-length}
    \slodd_{G_I}(\delta_I)\le h.
  \end{equation}
\end{enumerate}
\end{theorem}

\begin{proof}
Suppose first that a diagram has source $\Sigma_{r,1}$, where $0\le r\le h$. For standard generators $a_i,b_i$, the boundary relation is $[\partial\Sigma_{r,1}]=\prod_{i=1}^{r}[a_i,b_i]$, up to conjugacy and inversion. Since $w_{1,\Sigma_{r,1}}=0$, orientation-twistedness gives $\varepsilon_I(F_*a_i)=\varepsilon_I(F_*b_i)=0$. Applying $F_*$ to the boundary relation yields an even-commutator factorization of $\delta_I$, and therefore $\cleven_{G_I}(\delta_I)\le r\le h$.

Conversely, put $r=\cleven_{G_I}(\delta_I)\le h$ and choose a factorization $\delta_I=\prod_{i=1}^{r}[a_i,b_i]$ with every $a_i,b_i$ even. Choose based edge-loop representatives $\alpha_i,\beta_i$ of $a_i,b_i$, respectively, and let $c_0$ represent the conclusion boundary. If $r>0$, minimality ensures that every commutator is nontrivial, so reduced representatives of its two factors are nonconstant. For $r=0$, there are no paired arcs. The loop $W_{\mathrm{or}} =\left(\prod_{i=1}^{r}[\alpha_i,\beta_i]\right)c_0^{-1}$ is null-homotopic in $X_I^0$, so choose a van Kampen disk over the premise cells with boundary word $W_{\mathrm{or}}$. Glue each $\alpha_i$-arc to its inverse arc with opposite boundary directions, and likewise for $\beta_i$. Corresponding points have the same image in $X_I^0$, so the disk map descends. These are the standard handle identifications, giving a diagram on $\Sigma_{r,1}$ with remaining boundary $c_0^{-1}$. The standard generators have even images, so $\varepsilon_I\circ F_*=0=w_{1,\Sigma_{r,1}}$. This is Zhang's handle-gluing construction \cite[Lemma~3.4 and Theorem~3.6]{ZhangTiledGenus}.

Now let the source be $N_{r,1}$, where $1\le r\le h$, with standard one-sided generators $x_1,\ldots,x_r$. They satisfy $w_{1,N_{r,1}}(x_j)=1$ and give the boundary relation $[\partial N_{r,1}]=x_1^2\cdots x_r^2$, again up to conjugacy and inversion. An orientation-twisted diagram therefore gives $\delta_I=\prod_{j=1}^{r}F_*(x_j)^2$ with $\varepsilon_I(F_*x_j)=1$, so $\slodd_{G_I}(\delta_I)\le r\le h$.

Conversely, put $r=\slodd_{G_I}(\delta_I)\le h$ and choose a factorization $\delta_I=\prod_{j=1}^{r}z_j^2$ with every $z_j$ odd. Choose based edge loops $\zeta_j$ with $[\zeta_j]=z_j$, and take a van Kampen disk over $X_I^0$ for $W_{\mathrm{non}} =\left(\prod_{j=1}^{r}\zeta_j^2\right)c_0^{-1}.$ For each $j$, glue the two copies of $\zeta_j$ with the same boundary direction, matching corresponding points of the path. The map descends, and these are the standard crosscap identifications. Every $\zeta_j$ has odd length and is nonconstant, so the quotient is $N_{r,1}$, with remaining boundary $c_0^{-1}$. On every standard generator, both $\varepsilon_I\circ F_*$ and $w_{1,N_{r,1}}$ take the value one. Thus the resulting diagram is orientation-twisted and has nonorientable genus $r\le h$.

The conclusion boundary $c_0$ is a simple four-cycle in $\Gamma_I$, because its four formal labels are distinct. Since a graph has no two-cells, $[c_0]\ne0$ in $H_1(\Gamma_I;\Ftwo)$. Every commutator and every square has zero class in this homology group, so both $W_{\mathrm{or}}$ and $W_{\mathrm{non}}$ have class $[c_0]\ne0$. Neither loop is null-homotopic in the graph alone, so each disk contains a two-cell of $X_I^0$, necessarily a premise cell. This adapts Zhang's nondegeneracy argument \cite[Lemma~3.3]{ZhangTiledGenus} to $\Ftwo$-homology.

No new two-cells are introduced by the identifications, and each open face maps homeomorphically to an open premise cell. Reversal of the conclusion boundary is allowed. Capping it with $Q_0$ gives the corresponding closed tiling proof, and the parity compatibility extends over the cap as observed above.
\end{proof}

\begin{remark}\label{rem:exact-nonorientable-genus}
For a diagram on the specified surface $N_{h,1}$, the same construction requires an exact $h$-term factorization
\begin{equation}\label{eq:exact-odd-square-factorization}
  \delta_I=\prod_{j=1}^{h}z_j^2,
  \qquad \varepsilon_I(z_j)=1\quad(1\le j\le h).
\end{equation}
An odd-square factorization can be enlarged by two terms, by appending $u^2(u^{-1})^2$ for an odd factor $u$, but need not admit an enlargement by one term. Thus existence at genus at most $h$ does not by itself give existence at genus exactly $h$.
\end{remark}

\subsection{Euler Genus}

We extend Zhang's elementary genus \cite[Definition~2.10]{ZhangTiledGenus} to orientation-twisted tiling proofs, measuring the supporting surface by its Euler genus. Let $\eg(I)$ be the least Euler genus of a closed surface supporting an orientation-twisted tiling proof of $I$, with value $\infty$ if no such proof exists. The Euler genus of a closed surface is $2-\chi$. For a closed orientable surface of genus $g$, it is $2g$; for a closed nonorientable surface of genus $k$, it is $k$.

The following definition helps us to study the Euler genus of a tiling surface uniformly.

\begin{definition}[Twisted length]
  For $\gamma\in\ker\varepsilon_I$, define its \emph{twisted length} by
  \[
    \tl_{G_I}(\gamma)
    :=\min\left\{
      \slodd_{G_I}(\gamma),
      2\cleven_{G_I}(\gamma)
    \right\}.
  \]
\end{definition}

With this notation, Theorem~\ref{thm:surface-formula} states $\eg(I)=\tl_{G_I}(\delta_I)$. It extends Zhang's genus formula; the proof also extends his construction of an optimal tiling \cite[Theorems~1.9 and~3.6]{ZhangTiledGenus}.

\begin{proof}[Proof of Theorem~\ref{thm:surface-formula}]
Delete the conclusion square from a tiling proof. If its closed surface is orientable of genus $g$, Theorem~\ref{thm:parity-tiling-criterion}\textup{(i)} gives $\cleven_{G_I}(\delta_I)\le g$, while the Euler genus is $2g$.

If the surface is nonorientable of genus $k$, Theorem~\ref{thm:parity-tiling-criterion}\textup{(ii)} gives $\slodd_{G_I}(\delta_I)\le k$, while the Euler genus is $k$. Hence the Euler genus of any tiling proof is at least the right-hand side of \eqref{eq:euler-genus}.

If both factorization lengths are infinite, the forward implications of Theorem~\ref{thm:parity-tiling-criterion} exclude every tiling proof, so both sides of \eqref{eq:euler-genus} are infinite. Otherwise, a finite shortest even-commutator factorization of length $g$ produces a diagram on $\Sigma_{g,1}$, and hence a closed proof of Euler genus $2g$; a finite shortest odd-square factorization of length $k$ produces a diagram on $N_{k,1}$, and hence a closed proof of Euler genus $k$. The empty commutator factorization gives the sphere construction when $g=0$. Choosing a branch that realizes the minimum proves the reverse inequality and the construction statement. If the resulting surface had smaller Euler genus, the forward direction of the criterion would give a smaller value of the minimum in \eqref{eq:euler-genus}, a contradiction.
\end{proof}

\begin{remark}[Recovery of Zhang's formula]
In the ordinary bipartite setting every closed edge path is even, so $\varepsilon_I=0$, $\cleven=\cl$, and $\slodd=\infty$. Formula~\eqref{eq:euler-genus} then becomes $\eg(I)=2\cl_{G_I}(\delta_I)=2g_{\mathrm{or}}(I)$, which is precisely the Euler-genus form of \eqref{eq:zhang-genus}.
\end{remark}

\begin{example}[The Hesse presentation]\label{ex:hesse-length}
For the hemicube presentation, $X_I$ is the complement of the open conclusion face in the hemicube and hence a M\"obius band. Thus
\[
  G_I\cong\mathbb Z=\langle x\rangle,
  \qquad \varepsilon_I(x)=1,
  \qquad \delta_I=x^2.
\]
It follows that $\slodd_{G_I}(\delta_I)=1$. Since $\ker\varepsilon_I=\langle x^2\rangle$ is abelian and $x^2\ne1$, one has $\cleven_{G_I}(\delta_I)=\infty$. Therefore $\eg(I)=1$: an optimal surface-tiling diagram is the M\"obius band, and capping its conclusion boundary gives $\mathbb{RP}^2$. No orientable tiling proof exists. Its canonical orientation lift is an annulus whose two conclusion boundaries are exchanged by $\tau$, as in Theorem~\ref{thm:canonical-lift}.
\end{example}

\bibliographystyle{plain}
\IfFileExists{ref.bib}{\bibliography{ref}}{\bibliography{ref}}

\end{document}